\documentclass[12pt,twoside]{report}

\usepackage{amsfonts}       %
\usepackage[leqno]{amsmath} 
\usepackage{amssymb}        %
\usepackage{amstext}        %

\usepackage{amsthm}         
\usepackage{amscd}          
\usepackage{amsxtra}
\usepackage{makeidx}        
\makeindex

\usepackage{latexsym}
\usepackage[dvips]{epsfig}
\usepackage[all,dvips]{xy}  
\usepackage{multicol}
\usepackage{enumerate}      

\usepackage{color}

\newtheorem{meinzaehler}{ist quatsch}[chapter]
\newtheorem{thm}[meinzaehler]{Theorem}
\newtheorem{Def}[meinzaehler]{Definition}
\newtheorem{prop}[meinzaehler]{Proposition}
\newtheorem{lemma}[meinzaehler]{Lemma}

\newtheorem{expl}[meinzaehler]{Example}

\newtheorem{remark}[meinzaehler]{Remark}

\newcommand{\Numerierung}
    {\refstepcounter{meinzaehler}{\bf(\arabic{chapter}.\arabic{meinzaehler}) }
    }
\newcommand{\meinchapter}[1]{\refstepcounter{chapter}\section*{\arabic{chapter} $\ $ #1}
    \addcontentsline{toc}{chapter}{\arabic{chapter}$\ $ #1}
    \markboth{\scshape \arabic{chapter}. #1}{}
}

\newcommand{\Absatz}{\par\bigskip\par}
\newcommand{\absatz}{\par\medskip\par}
\newcommand{\klabsatz}{\par\smallskip\par}

\newcommand{\Seitenumbruch}{\clearpage}

\newcommand{\A}{\mathbb{A}}
\newcommand{\N}{\mathbb{N}}
\newcommand{\Z}{\mathbb{Z}}
\newcommand{\Q}{\mathbb{Q}}
\newcommand{\R}{\mathbb{R}}
\newcommand{\C}{\mathbb{C}}
\newcommand{\F}{\mathbb{F}}
\newcommand{\G}{\mathbb{G}}

\newcommand{\GG}{\mathbb{G}}

\newcommand{\QQ}{\mathbb{Q}}
\newcommand{\RR}{\mathbb{R}}

\newcommand{\ZZ}{\mathbb{Z}}

\newcommand{\CCC}{\mathcal{C}}

\newcommand{\LLL}{\mathcal{L}}
\newcommand{\MMM}{\mathcal{M}}
\newcommand{\NNN}{\mathcal{N}}
\newcommand{\PPP}{\mathcal{P}}

\newcommand{\VVV}{\mathcal{V}}
\newcommand{\WWW}{\mathcal{W}}

\newcommand{\ZZZ}{\mathcal{Z}}

\newcommand{\aaa}{\mathfrak{a}}
\newcommand{\bbb}{\mathfrak{b}}

\newcommand{\gGg}{\mathfrak{g}}
\newcommand{\hhh}{\mathfrak{h}}
\newcommand{\kkk}{\mathfrak{k}}
\newcommand{\lLl}{\mathfrak{l}}
\newcommand{\mmm}{\mathfrak{m}}
\newcommand{\nnn}{\mathfrak{n}}

\newcommand{\ppp}{\mathfrak{p}}
\newcommand{\qqq}{\mathfrak{q}}

\newcommand{\SsS}{\mathfrak{S}}
\newcommand{\ttt}{\mathfrak{t}}
\newcommand{\uuu}{\mathfrak{u}}

\newcommand{\Groth}{{\mathcal G}\bf\it{ro}}

\newcommand{\tran}{{}^t}

\newcommand{\bs}{\backslash}
\newfont{\mycyr}{wncyr10 scaled 1200}

\newcommand{\GL}{\text{GL}}
\newcommand{\SL}{\text{SL}}

\newcommand{\SO}{\text{SO}}
\newcommand{\GSO}{\text{GSO}}

\newcommand{\Sp}{\text{Sp}}
\newcommand{\SP}{\text{Sp}}
\newcommand{\GSp}{\text{GSp}}
\newcommand{\PGL}{\text{PGL}}
\newcommand{\Spin}{\text{Spin}}
\newcommand{\GSPIN}{\text{GSpin}}

\newcommand{\set}[1]{\left\{#1\right\}} 

{%
        \begin{list}{}{\setlength{\parsep}{0pt}\setlength{\itemsep}{0pt}
                    \setlength{\labelsep}{8pt}\setlength{\labelwidth}{15pt}
                      }
}{\end{list}}

\definecolor{hellgrau}{gray}{0.9}
\definecolor{mittelgrau}{gray}{0.7}
\definecolor{mittelgrau2}{gray}{0.705}
\definecolor{dunkelgrau}{gray}{0.5}
\definecolor{dunkelgrau2}{gray}{0.51}

\begin{document}

\pagestyle{headings}

\markboth{}{}
\title{Siegel Modular Varieties and the Eisenstein Cohomology of $\PGL_{2g+1}$}
\author{ {Uwe Weselmann}
}
\date{\today}
\maketitle
\begin{abstract}
We use the twisted topological trace formula developed in \cite{Uwe} to understand liftings from symplectic to general linear groups. We analyse the lift from $\SP_{2g}$ to $\PGL_{2g+1}$ over the ground field $\Q$ in further detail, and we get a description of the image of this lift of the $L^2$ cohomology of $\SP_{2g}$ (which is related to the intersection cohomology of the Shimura variety attached to $\GSp_{2g}$) in terms of the Eisenstein cohomology of the general linear group, whose building constituents are cuspidal representations of Levi groups. This description may be used to understand endoscopic and CAP-representations of the symplectic group. 
\end{abstract}

\thispagestyle{empty}
\Seitenumbruch 
\pagenumbering{roman}\setcounter{page}{1}
\markboth{}{}
\tableofcontents\thispagestyle{empty}        \Seitenumbruch
\markboth{}{}\thispagestyle{empty}
\pagenumbering{arabic} \setcounter{page}{1}
\meinchapter{Introduction}
This paper is a sequel of our paper \cite{Uwe}, where we developed a twisted topological trace formula and tried to understand liftings from symplectic to general linear groups. Here we want to analyse the lift from $\SP_{2g}$ to $\PGL_{2g+1}$ over the ground field $\Q$ in further detail, and we get a description of the image of this lift for the $L^2$ cohomology of $\SP_{2g}$ (which is
related to the intersection cohomology of the Shimura variety attached to $\GSp_{2g}$) in terms of the Eisenstein cohomology of
the general linear group, whose building constituents are cuspidal representations of Levi groups.
This description may be used to understand endoscopic and CAP-representations of the symplectic group. Roughly speaking we do not classify cohomology classes for the symplectic group with respect to their cuspidal support on 
$\Sp_{2g}$ in the sense of \cite{FrankeSchwermer} but in terms of their cuspidal support on their lifts to 
$\PGL_{2g+1}$. In view of the strong multiplicity one theorem on $\GL_{2g+1}$ this gives a classification of all "cohomological" packets (equivalence classes of automorphic representations which are isomorphic at all but finitely many places) in the discrete spectrum of $\Sp_{2g}$. 

\absatz
Let us describe what we will do in more detail, starting with  the heart of our arguments:

\absatz
In chapter \ref{Ch Franke} we use the fundamental work of Franke \cite{Franke} to write the cohomology of $G$ as a sum of parabolically induced $L^2$ cohomologies attached to the different Levi subgroups of $G$ and with respect to suitable coefficient systems on these Levis. For $G=\PGL_{2g+1}$ the index  sum may be rewritten in terms of partitions of the set $J=\{-g,\ldots,g\}$. In chapter \ref{Ch Partitions} we compare these sums for the goups in question  using an induction procedure and get that the
$L^2$-cohomology attached to $G=\SP_{2g}$ lifts to a sum, indexed by  $\eta$-invariant partitions of $J$, of the parabollically induced $L^2$-cohomology
 attached to suitable Levi factors (which are products of general linear groups) and to suitable coefficient systems (Theorem \ref{diskretes Lifting theorem}). Here a partition $J=\bigcup_{i=1}^r J_i$ is called $\eta$-invariant if we have
 $-J_i=J_i$ for $i=1,\ldots,r$.
\absatz Now one can go further: the discrete  spectrum of general linear groups is well understood by the work of
   Moeglin and Waldspurger \cite{MoglinWald}: each automorpic representation in the discrete spectrum is a Langlands quotient
   $MW(\pi^{(i)},n_i)$
   of some parabollically induced representation, in  which a cuspidal $\pi^{(i)}$ appears $n_i$ times with various
   twists by powers of the modulus character. The extended version of the strong multiplicity for $\GL_n$ (\cite{JacquetSmult})
   then implies that an $\eta$-invariant automorphic representation contributing to the lifted cohomology on $\PGL_{2g+1}$ comes
   from $\eta$-invariant cuspidal representations $\pi^{(i)}$ on some smaller linear groups $G^{(i)}$.
   By the work of Cogdell, Kim, Piatetski-Shapiro, Shahidi and Soudry (\cite{Cogdell01}, \cite{Cogdell04}, \cite{Soudry Paris}, \cite
   {Cogdell}) one knows that $\pi^{(i)}$ (resp. the twist of $\pi^{(i)}$ by a quadratic character in the case $G^{(i)}=\GL_{2\gamma_i+1}$)  is a semi-weak lift of a globally generic cuspidal representation $\pi_1^{(i)}$ on some quasi-split classical group $G_1^{(i)}$.
 \absatz
   In case $G^{(i)}=\GL_{2\gamma_i+1}$ one knows that $G_1^{(i)}=\SP_{2\gamma_i}$. In case $G^{(i)}=\GL_{2\gamma_i}$
   one knows that $G_1^{(i)}=\SO_{2\gamma_i+1}$ if $L(s,\pi^{(i)},\Lambda^2)$ has a pole at $s=1$, but
   $G_1^{(i)}=\SO_{2\gamma_i}^d$ if $L(s,\pi^{(i)},Sym^2)$ has a pole at $s=1$. Here the discriminant  of the quadratic form
   describing the quasi split $\SO_{2\gamma_i}^d$ is related via class field theory to the central character $\omega_i$
   of $\pi^{(i)}$, which is quadratic. If $\omega_i\ne 1$ we are always in the case of an even orthogonal group.
   But since the notion of a semi-weak lift means that the lift is compatible with the local Langlands isomorphism at the archimedean prime, we do not have to look for poles of $L$-functions: the fact that $\pi^{(i)}$ is cuspidal implies that
   $\pi^{(i)}_\infty$ is tempered and then the Langlands data are completely determined by the coefficient system.
   But by lemma \ref{Unterscheidungslemma}
    already the coefficient system decides, whether we are in the odd or in the even orthogonal case.
  \absatz
   If one considers this paper  as a contribution to "Langlands-functoriality for the cohomology of arithmetic groups", one
   has to think about the question, what  Langlands-functoriality for coefficient systems $V_\chi$ should mean.
   Here $V_\chi$ is the highest weight module for some dominant weight $\chi\in X^*(T)$ on a maximal torus $T\subset G$.
    We will see, that $\chi+\delta_G\in X^*(T)\otimes \C=X_*(T)\otimes \hat T$ behaves functorial, where $\delta_G$ is as usual the half sum of the positive roots and where $\hat T$ is a maximal torus in the dual group: The action of the center of the universal enveloping algebra $\ZZZ(\gGg)$ on $V_\chi$ is given by $\chi+\delta_G$, and by Wigner's lemma this is up to sign the
   action of $\ZZZ(\gGg)$ on the real representations $\pi_\infty$ contributing to cohomology. Now this action is related to the Langlands parameters of $\pi_\infty$, which behave functorial.
   In the case of $\GL_{2g+1}$ this means that an $\eta$-invariant coefficient system is already described by the (characteristic) set $S$  of all (integral) entries of the vector  $\chi+\delta_G\in X^*(T)=\Z^{2g+1}$, since these entries form a strictly increasing sequence. Here $\eta$-invariance means $-S=S$.
   \absatz
   In order to include the construction of M\oe glin-Waldspurger to describe the discrete spectrum of general linear groups into the notations we define characteristic sets in chapter \ref{Ch Characteristic Sets} as finite non empty subsets $S$ of $\frac 12\Z$, such that $-S=S$ and $s-t\in\Z$ for all $s,t\in S$. These are of type $C_g$ if $0\in S$,
   of type $B_g$ if $S\subset \frac 12+\Z$ and finally of type $D_g$ if $S\subset \Z$, but $0\notin S$.
   For $n\ge 1$ and a characteristic set $S$ we define $MW(S,n)=\{s-\frac{n+1}2+i| 1\le i\le n\}$ if $S$ is $n$-admissible in the sense of definition \ref{Characteristic Sets}(b).
   \absatz
   After this short set theoretic beginning we recall in chapter \ref{Ch Endosopic groups} the notion of
   (stable) twisted endoscopy, introduce the groups involved in the game and relate the coefficient systems to characteristic sets.
   In chapter \ref{Ch Real} we analyse the real representations contributing to cohomology for a given coefficient system and their Langlands parameters. Furthermore we calculate the Lefschetz numbers for the $\eta$-action on the cohomology
   of the real representations.
   \absatz
   In chapter \ref{Ch Statement} we formulate the main theorem about the structure of lifts (Theorem \ref{Maintheorem}): To an
   irreducible automorphic representation $\tau$ of $\SP_{2g}(\A)$ in the discrete spectrum, which is cohomological with respect to some coefficent system $V_\chi$ with characteristic set $S=S_\chi$ one associates a family of octupels:
   $$(X_{\gamma_i}, G^{(i)}, G^{(i)}_1, S_i,n_i,d_i,\pi^{(i)},\pi_1^{(i)})_{i=1,\ldots,r},$$
   which describes the lift of $\tau$ to $\PGL_{2g+1}$ in a way we already discussed.
   One should not be afraid of this description by octupels: to get all possible such families one should at first
   look for partitions of $S$ into integral characteristic sets: $S=\bigcup_{i=1}^r \tilde S_i$.
   One of them is of type $C_{\gamma_i}$, the others are of type $B_{\gamma_i}$.
   \klabsatz
   Then one should
   look for integers $n_i$ such that $\tilde S_i=MW(S_i,n_i)$ where $S_i$ is an $n_i$-admissible characteristic set.
   Observe that for every $i$ one  always can take $n_i=1$ since $MW(S,1)=S$.
   Let $X_{\gamma_i}$ be the type of $S_i$, and let $d_i\in \Q^*/(\Q^*)^2$ be a family of square classes such that
   $d_i=1$ for type $B_{\gamma_i}$, $(-1)^{\gamma_i}d_i>0$ for type $D_{\gamma_i}$ and $\Pi_{i=1}^r d_i=1$.
     Then the $S_i$ and $d_i$ determine the classical groups $G_1^{(i)}$ and the linear groups
   $G^{(i)}$ in question. Finally $(\pi^{(i)},\pi_1^{(i)})$ is a pair of irreducible globally generic cuspidal
   automorphic representations for these two groups, which have the correct Langlands parameters at the archimedean prime and such that $\pi^{(i)}$  has central character $\chi_{d_i}$ and is a semi weak lift of $\pi_1^{(i)}$.
   The proof of the main theorem and its converse is given in the later chapter \ref{Ch Proof} and includes the arguments
   already sketched.
   \absatz
   In the last chapter \ref{Ch Examples} we describe  the elementary pairs  $(\pi^{(i)},\pi_1^{(i)})$ in the case of
   characteristic sets of type $C_1,D_1,B_1$ and $D_2$, we analyse the case $g=2$ in further detail, and give hints to understand
   the weights appearing in the recent work of Bergstr\"om, Faber and van der Geer \cite{BFGeer} in case $g=3$, explain how  the Ikeda lift fills into our language and finally try to give
   a  relation to elliptic endoscopic groups for $\GSp_{2g}$ as appearing in the work of Morel \cite{Morel}.
   \absatz
   Finally we have to mention some technical computations: We will calculate and compare the factors $\alpha_\infty(\gamma_0,h_\infty)$ from the topological trace formula (restated in chapter
\ref{Ch Restatement}) in chapter \ref{Ch alpha}, after we carried out in chapter \ref{Ch Chevalley} a comparison of volumes of the real compact forms of $\SP_{2g}$
and $\SO_{2g+1}$ with respect to measures attached to Chevalley basis. This will make the comparison of topological trace formulas more
explicit, and we get a more precise meaning  of the statement that the  cohomology attached to $\PGL_{2g+1}$ is a lift of the cohomology attached to  $\SP_{2g}$. Here we mean by "cohomology attached to an algebraic group $G$" the alternating sum of the cohomologies of locally symmetric spaces attached to $G$ with coefficients in some local system coming from a finite dimensional algebraic representation $V_\chi$ of $G$.

   \absatz
   It should be noted that their may be some overlap of our work  with the forthcoming book of Arthur \cite{Arthur}: His results depend on the stabilization of the twisted analytic trace formula, which seems not to be established yet,  but give much more precise statements about the local and global packets of (automorphic) representations of $\SP_{2g}(\A)$. But for applications to cohomology and for a
   first understanding of the lift to $\PGL_{2g+1}$ our approach seems to be more direct.

\meinchapter{Characteristic sets}\label{Ch Characteristic Sets}

\begin{Def}\label{Characteristic Sets}
      \begin{itemize}\item[(a)] A characteristic set is a  finite subset $S\ne\emptyset$ of $\frac 12 \Z$, such that $-s\in S$ for all $s\in S$ and $s-t\in \Z$ for all
      $s,t\in \Z$
      \item[(b)] For $n\in \N$ a characteristic set is called $n$-admissible if we have $|s-t|\ge n$ for all $s,t\in S$ with $s\ne t$.
      \item[(c)] For $n\in \N$ we call $E_n=\set{-\frac{n+1}2+i| \; 1\le i\le n}$  the elementary characteristic set of order $n$.
      \item[(d)] If $S$ is an $n$-admissible characteristic set we put $$MW(S,n)=\{s+e\;|\;s\in S,e\in E_n\}.$$
       \end{itemize}
\end{Def}
\begin{lemma} Let  $n\in\N$ and $S$ be an $n$-admissible characteristic set.
   \begin{itemize}
     \item[(a)] Each $E_n$ and each $MW(S,n)$ is a characteristic set.
     \item[(b)] The map $\alpha: S\times E_n\to MW(S,n), \; (s,e)\mapsto s+e$ is a bijection.
   \end{itemize}
\end{lemma}
\qed

\begin{remark} \em It is obvious that each characteristic set $S$
is of exactly one of the following types for some $g\in\N$:
     \begin{itemize}
     \item[(a)] $S$ is of type $B_g$ ($g\ge 1$), if $S\subset \Z+\frac 12$ and $\#(S)=2g$;
     \item[(b)] $S$ is of type $C_g$ ($g\ge 0$), if $0\in S$ and $\#(S)=2g+1$;
     \item[(c)] $S$ is of type $D_g$ ($g\ge 1$), if $S\subset \Z, 0\notin S$ and $\#(S)=2g$.
   \end{itemize}
   The characteristic sets of type $C_g$ and $D_g$ are called {\it integral}.
   Let $S$ be an $n$-admissible characteristic set. If $n$ is odd, then $MW(S,n)$ is of the same type (but with different index) as $S$.
   If $n>0$ is an even integer, then $MW(S,n)$ is integral if and only if $S$ is not integral.
\end{remark}

\meinchapter{Endoscopic groups and coefficient modules}\label{Ch Endosopic groups}

\Numerierung {\sc Split Groups with automorphism.}\label{splitmitauto}
Let $G/F$ be a connected reductive split group over a field $F$, and let  $(B,T,\{X_\alpha\}_{\alpha\in \Delta})$ be a splitting, i.e.  $T$ is a maximal split torus
inside an $F$-rational Borel $B$, $\Delta=\Delta_G=\Delta(G,B,T)\subset \Phi(G,T)\subset X^*(T)=Hom(T,\GG_m)$ denotes the set
 of simple roots inside the system of roots, and
the $X_\alpha$ for the simple roots $\alpha\in\Delta$ are a system (nailing)  of  generators of the root spaces $\gGg_\alpha$
in the Lie algebra.
Let  $\eta\in Aut(G)$ be an automorphism of $G$ of finite order $l$,  which fixes the
splitting. We denote by
$\tilde G=G\rtimes \langle \eta\rangle$
the (nonconnected) semidirect product of $G$ with $\eta$.
 $\eta$ acts on the cocharacter module $X_*(T)=Hom(\GG_m,T)$ via
$X_*(T)\ni \alpha^\vee \mapsto \eta\circ\alpha^\vee$ and on  $X^*(T)$ via $\alpha\mapsto \alpha\circ
\eta^{-1}$.
\Absatz
\Numerierung{\sc The dual group.}
Let $\hat G=\hat G(\C)$ be the dual group of $G$. By de\-fi\-nition $\hat G$ has a splitting
$(\hat B,\hat T,\{\hat X_{\hat\alpha}\})$ such that we have identifications
$X^*(\hat T)=X_*(T), \enskip X_*(\hat T)= X^*(T)$, which identify the (simple) roots
$\hat\alpha\in X^*(\hat T)$ with the (simple) coroots $\alpha^\vee \in X_*(T)$ and the (simple) coroots
$\hat\alpha^\vee\in X_*(\hat T)$ with the (simple) roots $\alpha \in X^*(T)$. There exists a
unique automorphism $\hat\eta$ of $\hat G$ which stabilizes $(\hat B,\hat T,\{\hat X_{\hat\alpha}\})$
and induces on $(X_*(\hat T),X^*(\hat T))$ the same automorphism as $\eta$ on $(X^*(T),X_*(T))$.
\Absatz
For some $\hat s\in \hat T$ let $\hat H=(\hat G^{\hat\eta\hat s})^\circ$ be the connected component of the subgroup
 of elements in $\hat G$ which are fixed under $\hat\eta\circ int(\hat s)$.
 It is a reductive group with a splitting satisfying
 $\hat B_H=\hat B^{\hat\eta\hat s}$,
$\hat T_H=\hat T^{\hat\eta}$ and where in the case $\hat s=1$ the family $\{\hat X_{\hat\beta}\}_{\beta\in\Delta_{\hat H}}$
is described in \cite[5.3 and 5.4]{Uwe}.

\begin{Def}[ $\eta$-endoscopic group]
In the above situation a connected reductive split group  $H/F$ will be called
an $\eta$-endoscopic group for $(G,\eta,\hat s)$  if its dual group together with the
splitting is isomorphic to the above $(\hat H,\hat B_H,\hat T_H,\{X_\beta\}_{\beta\in\Delta_{\hat H}})$.
In the case $\hat s=1$ we call $H$ a stable $\eta$-endoscopic group for $(G,\eta)$.
\end{Def}
 For a maximal split torus $T_H\subset H$ we have:
\begin{align}
 \label{RelationX^*} X^*(T_H)\;=\;X_*(\hat T_H)
\quad &=\quad X_*(\hat T)^{\hat\eta}\; = \; X^*(T)^\eta, \qquad\text{ and}\\
\label{RelationX_*}
 X_*(T_H)\;= \;X^*(\hat T_H)\quad &= \quad (X^*(\hat T)_{\hat\eta})_{free}\;=\;(X_*(T)_\eta)_{free} \qquad
\end{align} is the maximal free abelian quotient of the coinvariant
module $X^*(\hat T)_{\hat\eta}=X_*(T)_\eta$.

\Absatz\Numerierung {\sc The classical groups.}
  We use the following { notations}:
  \klabsatz\quad $diag(a_1,\ldots,a_n)\in \GL_n$ denotes the diagonal
  matrix $(\delta_{i,j}\cdot a_i)_{ij}$ and
  \klabsatz\quad $antidiag(a_1,\ldots,a_n)\in \GL_n$ the antidiagonal matrix
  $(\delta_{i,n+1-j}\cdot a_i)_{ij}$. Consider
  \begin{gather*}
    J=J_n=(\delta_{i,n+1-j}(-1)^{i-1})_{1\le i,j\le n}=antidiag(1,-1,\ldots,(-1)^{n-1})\in \GL_n(F).
  \end{gather*}
  and its modification $J'_{2g}=antidiag(1,-1,1,\ldots,(-1)^{g-1},(-1)^{g-1},\ldots,1,-1,1)$.

  Since $\tran J_n=(-1)^{n-1}\cdot J_n$  and since $J'_{2g}$ is symmetric we can define the
  \begin{align*} \text{standard symplectic group} \quad\qquad\Sp_{2g}\quad &= \quad \Sp(J_{2g})\\
    \text{standard split odd orthogonal  group}\qquad \SO_{2g+1}\quad &= \quad \SO(J_{2g+1}).\\
    \text{standard split even orthogonal  group}\quad\qquad \SO_{2g}\quad &= \quad \SO(J'_{2g}).
  \end{align*}
  We consider the groups $\GL_n,\SL_n,\PGL_n,\SP_{2g},\SO_n$ with the splittings consisting
  of the diagonal torus, the Borel consisting of upper triangular matrices and the standard nailing.
  The following map defines an involution of $\GL_n,\SL_n$ and $\PGL_n$:
  \begin{align}\label{eta} \eta=\eta_n: g\mapsto J_n\cdot \tran g^{-1}\cdot J_n^{-1}.
  \end{align}
\begin{expl}[$A_{2g}\leftrightarrow C_g$]\label{BspPGL2n+1Sp2n}
   \em
  The group  $G=\PGL_{2g+1}$ with the involution $\eta=\eta_{2g+1}$  has the dual group
  $\hat G=\SL_{2g+1}(\C)$ with involution $\hat\eta=\eta_{2g+1}$. Then
 $H=\SP_{2g}$ is a stable $\eta$-endoscopic group, since its dual is $\hat H=\SO_{2g+1}(\C)=\hat G^{\hat \eta}$.

\end{expl}
\Absatz
\begin{expl}[$A_{2g-1}\leftrightarrow B_g$]\label{BspGL2nGl1GSpin2n+1}
  \em
  The group $G=\GL_{2g}\times \G_m$ has the involution
  $ \eta:(\gamma,a)\mapsto (\eta_{2g}(\gamma),\det(\gamma)\cdot a)$ and the
   dual $\hat \eta \in Aut(\hat G)$ satisfies
  $\hat\eta(c,b)=(\eta_{2g}(c)\cdot b, b).$  Then $\hat H=\hat G^{\hat \eta}=\GSp_{2g}(\C)$ is the dual of
  the stable $\eta$-endoscopic group $H=\GSPIN_{2g+1}=\left(\G_m\times \Spin_{2g+1}\right)/\mu_2$.
  \absatz
  If we consider the involution $\eta_{2g}$ on the group $G=\GL_{2g}$ we get the stable $\eta$-endoscopic group
  $H=\SO_{2g+1}$ with dual group $\hat H=\hat G^{\hat \eta}=\Sp_{2g}(\C)$.
  \end{expl}

 \begin{expl}[$A_{2g-1}\leftrightarrow D_g$]\label{BspGL2nGl1GSpin2n}
  \em
  For the group $G=\GL_{2g}$, the involution
  $\eta:=\eta_{2g}$ and the element $\hat s=diag(1,\ldots,1,-1,\ldots,-1)\in \hat T$ we get $H=\SO_{2g}$ as
  an $\eta$-endoscopic group, since its dual is
    $ \hat H=(\hat G^{\hat \eta\hat s})^\circ=\SO_{2g}(\C)$.
   \absatz
   For the group $G=\GL_{2g}\times \G_m$, the involution
  $\eta$ as in example \ref{BspGL2nGl1GSpin2n+1} and the element $\hat s=(diag(1,\ldots,1,-1,\ldots,-1);1)\in \hat T$
  we get $H=\GSPIN_{2g}$ as
  an $\eta$-endoscopic group, since its dual is
    $ \hat H=(\hat G^{\hat \eta\hat s})^\circ=\GSO_{2g}(\C)$ (see \ref{GSpin2g} below for details).
\end{expl}

\Absatz\Numerierung{\sc Irreducible finite dimensional representations.}
The isomorphism classes of finite dimensional algebraic representations  of $G/\bar \Q$ are in one to one correspondence
with the positive Weyl chamber
\begin{align}\label{posWeyl}
X^*(T)^+\quad =\quad \left\{ \chi\in X^*(T)| \;\langle\chi,\alpha^\vee\rangle \ge 0\text{ for all }\alpha\in \Delta\right\}
\end{align}
in such a way that $\chi\in X^*(T)^+$ corresponds to the irreducible representation $V_\chi=(V_\chi,\rho_\chi)$ of highest weight $\chi$.
We say that an irreducible finite dimensional representations
$\tilde V$ of $\tilde G$ is $\eta$-extended if it remains irreducible after restriction to $G$ and if $\eta$ acts as identity on the one dimensional space of highest weight vectors.
Since we have ${}^\eta V_\chi= V_{\eta(\chi)}$ (where ${}^\eta V_\chi=(V_\chi,\rho_\chi\circ \eta)$) we get a bijection between the isomorphism classes of $\eta$-extended irreducible finite dimensional representations
$\tilde V_\chi$ of $\tilde G$
and the set $X^*(T)^+\cap X^*(T)^\eta$. Denote by $\tilde V_\chi$ the $\eta$-extended representation which restricts to $V_\chi$.

\klabsatz
Let $\rho=\rho_G=\frac 12\sum_{\alpha>0}\alpha \in X^*(T)\otimes \Q$ be half the sum of the positive roots.

\begin{lemma}\label{Repkorrespondenz} Let $H$ be a stable $\eta$-endoscopic group for $(G,\eta)$ and let $T_H,T$ be maximal split tori.
Using  the identification $X^*(T_H)=X^*(T)^\eta$ we have:
\begin{itemize}
\item[(a)]  $X^*(T_H)^+=X^*(T)^\eta\cap X^*(T)^+$
\item[(b)]   $\rho_G=\rho_H$
\item[(c)]
the maps $\chi\mapsto V_{H,\chi},\quad\chi\mapsto \tilde V_\chi,$ give  bijections between
   \begin{itemize}
   \item the set $X^*(T_H)^+=X^*(T)^\eta\cap X^*(T)^+$

   \item the isomorphism classes of  finite dimensional irr. representations of $H$

   \item  the isom. classes of $\eta$-extended finite dim. irr. representations of $\tilde G$.

   \end{itemize}
   \end{itemize}
   \end{lemma}
 Proof: (a) The simple coroots of $H$ are the images of the coroots of $G$ under the natural projection
  $p:X_*(T)\mapsto X_*(T_H)$ (\cite[Prop. 5.4.]{Uwe}).
  If we interpret coroots as linear forms on the character group, then $p$ has to be interpreted as the
  restriction of linear forms from $X^*(T)$ to $X^*(T_H)=X^*(T)^\eta$. In this sense the simple coroots of $H$ are the
  restrictions  of the simple coroots of $G$, and the claim is a consequence of the definition (\ref{posWeyl}).
\klabsatz
 (b) Assume without loss of generality that $G$ is simple.
   Recall (\cite[5.3.]{Uwe}) that the roots $\beta\in \Phi(H,T_H)$ are of the form $\beta=c(\alpha)\cdot S_\eta(\alpha)$,
 where  $S_\eta(\alpha)$ denotes the sum over all elements in the $\eta$-orbit of $\alpha\in \Phi(G,T)\subset X^*(T)$ and where
 $c(\alpha)\in \{1,2\}$. Here $\alpha$ runs through a set of representatives for those $\eta$-orbits
 in $\Phi(G,T)$, whose elements are not "long" in the sense that in the case of a root system of type $A_{2g}$, $ord(\eta)=2$ they cannot be written in the form $\alpha_0=\alpha_1+\eta(\alpha_1)$ with $\alpha_1\in \Phi(G,T)$ (call such $\alpha_1$ short).

 We have $c(\alpha)=1$ if $\alpha$ is not "short" and we have   $c(\alpha_1)=2$ if $\alpha_1$ is "short". In this latter case we have
 $\beta=2\cdot S_\eta(\alpha_1)= \alpha_0 +S_\eta(\alpha_1)$. These identities imply the claim $\rho_G=\rho_H$.
 \klabsatz
 (c) is a consequence of  (a) and the remarks preceding the lemma.
\qed
\Absatz\Numerierung {\sc Linear groups.}
If  $\eta$ is the  involution on $G=\GL_n,\PGL_n$ resp. $\GL_n\times \G_m$ as defined above,
  we can make things more explicit:
  For $G=\GL_n$ we denote by $(a_1,a_2,\ldots,a_n)\in \ZZ^n\simeq X^*(T)$ the character
  $\chi:diag(t_1,\ldots,t_n)\mapsto t_1^{a_1}\cdots t_n^{a_n}$. Similarly we denote
 a character of the group $\G_m^n\times \G_m$ by $(a_1,a_2,\ldots,a_n;a_0)\in \ZZ^n\times \ZZ$.
  The characters of the diagonal torus in $\PGL_n$ may be described by the set $\set{(a_1,\ldots,a_n)\in \ZZ^n |\sum_{i=1}^n a_i=0}$.
  The simple roots are of the form $\alpha_i=e_i-e_{i+1}\in X^*(T)$ for the standard basis $e_i$ of $\Z^n$.
  Since the simple coroots are of the same form $\alpha_i^\vee=e_i-e_{i+1}$,
  the positive Weyl chamber $X^*(T)^+$ is given by the inequalities $a_1\ge a_2\ge \ldots \ge a_n$. The
  automorphism $\eta$ acts by $(a_1,a_2,\ldots,a_n)\mapsto (-a_n,-a_{n-1},\ldots,-a_1)$ (cases $\GL_n$ and $\PGL_n$)
  respectively by $(a_1,\ldots,a_n;a_0)\mapsto (a_0-a_n,\ldots,a_0-a_1;a_0)$
  in the case $G=\GL_n\times \GL_1$.
  We have $\rho=
  \left(\frac{n-1}{2},\frac{n-3}{2},\ldots,\frac{1-n}{2}\right)$ for $\GL_n$ and $\PGL_n$ resp.
   $\rho=\left(\frac{n-1}{2},\frac{n-3}{2},\ldots,\frac{1-n}{2};0\right)$ for $\GL_n\times \GL_1$.

   \Absatz
   If $\chi+\rho_G=(b_1,\ldots,b_n)$  in the  case $\GL_n$ and $\PGL_n$
   we call $S_\chi=\{b_1,\ldots,b_n\}\subset \frac 12\Z$ the {\bf describing set} of $\chi\in X^*(T)^+$.
   If $\chi+\rho=(b_1,\ldots,b_n;b_0)$ for $G=\GL_n\times \GL_1$ we call
   $S_\chi=\left\{ b_1-\frac {b_0}2,\ldots,b_n-\frac{b_0}2\right\}\subset \frac 12\Z$ the {\bf describing set} of $\chi\in X^*(T)^+$
   and $w(\chi)=b_0$ the {\bf weight} of $\chi$.
   \absatz
   Since the entries of $\chi+\rho$ are in a strictly decreasing order it is clear that $\chi\in X^*(T)^+$ is
   determined uniquely by the set $S_\chi$ in the cases $\GL_n$ and $\PGL_n$.
   In the case $\GL_n\times \GL_1$ the character $\chi$ is determined uniquely by the pair $(S_\chi,w(\chi))$.

   \begin{lemma}\label{Klassifikation C} Let $G=\PGL_{2g+1}$ and $H=\SP_{2g}$ be the stable $\eta$-endoscopic group.
   The map $  \chi \mapsto S_\chi$ gives a  bijection between
   \begin{itemize}
   \item the set $X^*(T_H)^+=X^*(T)^\eta\cap X^*(T)^+$ and

   \item  the set of characteristic sets of type $C_g$.
   \end{itemize}
   \end{lemma}
   Proof: In view of the description of the $\eta$-action on $X^*(T)$ given above and of
   $\rho_G=(g,\ldots,1,0,-1,\ldots,-g)$, it is clear that
   the describing set $S_\chi$ is a characteristic set of type $C_g$.
   It is easy to see, that the map is a bijection.
   \qed
   \begin{lemma}\label{Klassifikation B} Let $G=\GL_{2g}\times \G_m$ and $H=\GSPIN_{2g+1}$ its stable
   $\eta$-endoscopic group. The map
   $ \chi \mapsto (S_\chi,w(\chi))$ gives a  bijection between
   \begin{itemize}
   \item the set $X^*(T_H)^+=X^*(T)^\eta\cap X^*(T)^+$ and

   \item  the set of pairs $(S,w)$ where $w\in \Z$ and $S$ is a characteristic set of type $B_g$ for $w$ even
   and $S$ is a characteristic set of type $D_g$ for $w$ odd.
   \end{itemize}
   \end{lemma}
   We omit the easy proof.
   \qed
   \Absatz\Numerierung {\sc Even GSO.} Let $T_{GSO}=\{(diag(t_1,\ldots,t_{2g}),t_0)\;|\; t_it_{2g+1-i}=t_0\}$
   be the  diagonal torus in $H=\GSO_{2g}$. Using the identification $T_{GSO}=\G_m^{g+1}$ given by the map
   $(t_0;t_1,\ldots,t_g)\mapsto diag(t_1,\ldots,t_g, t_0t_g^{-1},\ldots, t_0t_1^{-1})$
   we get an isomorphism $X^*(T_{GSO})=\Z\times \Z^g$. Thus
   $(n_0;n_1,\ldots,n_g)$ denotes the character $(diag(t_1,\ldots,t_{2g}),t_0)\;\mapsto \Pi_{i=0}^g t_i^{n_i}$.
  The  $\alpha_i^\vee=e_i-e_{i+1}$ for $1\le i\le g-1$ and $\alpha_g^\vee=e_{g-1}-(e_0-e_g)$ are the simple coroots with
  respect to some splitting of the form $(T_{GSO},B_H,\set{ X_\alpha})$, where $B_H$ denotes the Borel of upper triangular matrices.
  Thus the
  positive Weyl chamber  $X^*(T_{GSO})^+$ is described by the conditions $n_1\ge \ldots\ge n_{g-1}\ge n_g\ge n_0-n_{g-1}$.
   Let $\sigma\in O_{2g}$ be the block diagonal matrix $(E_{g-1},\begin{pmatrix}0& 1\cr 1&0\cr\end{pmatrix},E_{g-1})$. Then
   the conjugation by $\sigma$ is an outer automorphism of $H$ which preserves some splitting as above.
   It acts on $X^*(T_{GSO})$ by $\sigma^*: (n_1,\ldots,n_g;n_0)\mapsto (n_1,\ldots,n_{g-1},n_0-n_g;n_0)$.
   A character $\chi=(n_0;n_1,\ldots,n_g)\in X^*(T_{GSO})^+$ is called $\sigma$-positive resp. $\sigma$-semipositive
   if $n_g-\frac {n_0}2 >0$ resp. $n_g-\frac {n_0}2 \ge 0$.
   It is clear, that each $\chi\in X^*(T_{GSO})^+$ fulfilles exactly one of the following three conditions:
   $\sigma^*\chi=\chi$, $\chi$ is $\sigma$-positive, $\sigma^*\chi$ is $\sigma$-positive.

   \Absatz\Numerierung \label{GSpin2g}{\sc General Spin groups.} Under the identification $X^*(T_{GSO})= \Z^{g+1}=\oplus_{i=0}^g \Z e_i$ and the dual
   identification $X_*(T_{GSO})= \Z^{g+1}$ the system of roots resp. coroots can be written as
   \begin{align*} \Phi(GSO)\; &= \; \set{ \pm(e_i-e_j)| 1\le i<j}\cup  \set{ \pm(e_i+e_j-e_0)| 1\le i<j}
    \; \subset X^*(T_{GSO})\\
     \Phi(GSO)^\vee\; &= \;\Phi_{SO}\;=\; \set{ \pm e_i \pm e_j| 1\le i<j\le g}    \;\subset \; X_*(T_{GSO}).\\
    \end{align*}
    If $\GSPIN_{2g}$ denotes the dual group of $\GSO_{2g}$ we get a root datum, where the roles of $\Phi$ and $\Phi^{\vee}$ are interchanged,
    i.e. $\Phi(\GSPIN)^{\vee}= \Phi(GSO)$ and $\Phi(\GSPIN)=\Phi_{SO}$.
    The  projection $\Spin_{2g}\to \SO_{2g}$ induces the following inclusion on the level of cocharacter groups
    of maximal tori:
    \begin{align*}X_*(T_{Spin})\;=\; &\left\{ \sum_{i=1}^g a_i e_i\in \Z^g \; \left| \; \sum_{i=1}^g a_i \equiv 0
     \;mod \; 2\right.\right\}\quad \hookrightarrow \quad X_*(T_{SO})=\bigoplus_{i=1}^g \Z e_i. \quad \\
    \end{align*}
    Then we may define an isogeny $\pi: \G_m\times \Spin_{2g}\to \GSPIN_{2g}$ by presribing its effect $\pi_*$ on the level
    of cocharacter groups of maximal tori:
    \begin{align*} X_*(\G_m)\times X_*(T_{Spin})\quad &\to \quad X_*(T_{GSpin})=\Z^{g+1},\\
     \quad \left(a_0,e_0,\sum_{i=1}^g a_i e_i\right)\quad &\mapsto \quad
       a_0e_0+\sum_{i=1}^g a_i e_i -\frac 12\cdot \left(\sum_{i=1}^g a_i \right)\cdot e_0.\\
    \end{align*}
    In fact it is easy to check that $\pi_*$ respects the coroots and its dual $\pi^*$ respects the roots.
    From $Hom(\ker(\pi),\G_m)={\rm coker}(\pi^*)\cong \Z/2\Z$ we deduce $\ker(\pi)\cong\mu_2$, and a more careful analysis implies that $\mu_2$ is embedded diagonally in $\G_m\times \Spin_{2g}$. Thus the dual group $\GSPIN_{2g}$ of $\GSO_{2g}$ is in fact what is usually called $\GSPIN_{2g}$.
With
    \begin{align*} X^*(T_{Spin})\;=  
     \;\left\{\left.\sum_{i=1}^g a_i e_i \right| a_i\in \frac 12\Z,\quad a_1\equiv\ldots\equiv a_g\; mod \;\Z\right\}\;\supset\; X^*(T_{SO})=\bigoplus_{i=1}^g \Z e_i
     \end{align*}
    the induced map $\pi^*$
   on the level of character groups of maximal tori is given by:
   \begin{align*}
     \pi^* :\quad X^*(T_{GSpin})=\Z^{g+1}\quad &\to \quad X^*(\G_m)\times X^*(T_{Spin}) =\Z \times X^*(T_{Spin})\\
      (m_0,m_1,\ldots, m_g)\quad&\mapsto \quad \left( m_0; m_1-\frac{m_0}2,\ldots,m_g-\frac{m_0}2\right).\\
   \end{align*}
   The inclusion $i: X^*(T_{GSpin})=X_*(T_{GSO})\hookrightarrow X_*(T_{\hat G})=X^*(T_G)$ where
   $T_G$ is the maximal split torus in $G=\GL_{2g}\times\G_m$, is given by the formula
   $$ (m_0,m_1,\ldots,m_g)\quad\mapsto \quad (m_1,\ldots,m_g,m_0-m_g,\ldots,m_0-m_1; m_0).$$
   The half sum of positive roots for $G=\GL_{2g}\times\G_m$ and for $\GSPIN$ is given by
   \begin{align*}
    \delta_G\quad &=\quad \left( g- \frac 12 ,\ldots,\frac 12,-\frac 12,\ldots,\frac 12 - g;0\right)\in X^*(T_G)\otimes \Q,\\
    \delta_{GSpin}&=\quad (0; g-1,\ldots,1,0)\in \Z^{g+1}=X^*(T_{GSpin}).
   \end{align*}
   Note that the composed map:
   $$\iota: X^*(T_G)^\eta \otimes \Q\to (X^*(\G_m)\times X^*(T_{Spin}))\otimes \Q,\quad \chi\mapsto \pi^*(i^{-1}(\chi+\delta_G)-\delta_{GSpin})$$
   is of the form
   $$(m_1,\ldots,m_g,m_0-m_g,\ldots,m_0-m_1; m_0) \mapsto \left(m_0; m_1-\frac{m_0-1}2,\ldots,m_g-\frac{m_0-1}2\right)\!,$$
   and thus maps the lattice $X^*(T_G)^\eta$ to the lattice $X^*(\G_m)\times X^*(T_{Spin})$.
   For $\chi\in X^*(T_G)^\eta$ with odd $m_0$ we have  $\iota(\chi)\in \Z^{g+1}= X^*(\G_m)\times X^*(T_{SO})$, so that $\iota(\chi)$ descends to a character
   of the maximal torus of $\G_m\times \SO_{2g}$. Let

    $X^*(T_H)^{++}=X^*(T_G)^\eta\cap X^*(T_G)^+$ and $X^*(T_H)^{++odd}=\set {\chi\in X^*(T_H)^{++}|
   \; w(\chi)\text{ odd}}$.
   For   $\chi\in X^*(T_H)^{++odd}$ let $V_{H,\chi}$ be the finite dimensional
   irreducible representation of $H=\G_m\times {\SO}_{2g}$ with highest weight $\iota(\chi)$.

   \begin{lemma}\label{Klassifikation D} Let $G={\GL}_{2g}\times \G_m$  and  $H=\G_m\times {\SO}_{2g}$.
   \klabsatz
   (a) $X^*(T_H)^{++}$ is the set of all $\sigma$-semipositive elements in $X^*(T_H)^+$, and
   $X^*(T_H)^{++odd}$ may also be described as the set of all $\sigma$-positive elements in $X^*(T_H)^+$ with odd $w(\chi)$.
   \klabsatz
   (b) The maps
   $\chi\mapsto V_{H,\chi},\quad\chi\mapsto \tilde V_\chi,\quad  \chi \mapsto (S_\chi,w(\chi))$ give  bijections between
   \begin{itemize}
   \item the set $X^*(T_H)^{++odd}$ of $\sigma$-positive dominant characters with odd $w(\chi)$,

   \item the isomorphism classes of  finite dimensional irreducible representations of $H$ with $\sigma$-positive highest weight on $\SO_{2g}$ and odd weight on $\G_m$,

   \item  the isomorphism classes of $\eta$-extended finite dimensional irreducible representations of $\tilde G$  which restrict to an odd character of $\G_m$,

   \item  the set of pairs $(S,w)$ where $S$ is a characteristic set of type $D_g$ and $w\in \Z$ is odd.
   \end{itemize}
   \end{lemma}
   Proof: (a) $\chi\in X^*(T_H)^+$  means $n_1\ge \ldots\ge n_{g-1}\ge n_g\ge n_0-n_{g-1}$ and $\chi$ is $\sigma$-semipositive
    if $n_g\ge n_0-n_{g}$. Both conditions are satisfied if and only if $\chi\in X^*(T_G)^\eta\cap X^*(T_G)^+$. For odd $w(\chi)$ each $\sigma$-semipositive
   $\chi$ is already $\sigma$-positive.
   \absatz
   (b) may be deduced  easily from the definitions and the considerations above.
   \qed

\meinchapter{Real representations}\label{Ch Real}

   \Absatz\Numerierung \label{Weilgruppenhom} {\sc Weil group homomorphisms}
   Let $ W_\R\;=\; \C^*\cup \C^*\sigma$ be the Weil group of $\R$ with the multiplication rules $\sigma\cdot c=\bar c\cdot \sigma$ for $c\in \C^*$ and
   $\sigma^2=-1\in \C^*$ .
   Recall that the group of continuous homomorphisms from $\C^*$ to $\C^*$ is isomorphic to
      $$\Xi\;=\;\set{(p,q)\in \C^2\;|\quad  p-q\in \Z} $$ in such a way that $(p,q)$ corresponds to the homomorphism
   $$ z\; \mapsto\; z^{(p,q)}:= z^p\cdot \bar z^q:= z^{p-q}\cdot (z\bar z)^q.$$
   Similarly the group of continuous homomorphisms from $\C^*$ to a complex torus $\hat T$ is isomorphic to
   $$X_*(\hat T)\otimes_\Z \Xi\;=\; \left\{ (\mu,\nu)\in (X_*(\hat T)\otimes\C)^2\; |\; \mu-\nu\in X_*(\hat T)\right\}.$$

   For $G=\GL_{n}\times \G_m$ with diagonal torus $T$ we associate to $\chi\in X^*(T)\otimes \Q$ with $2\chi\in X^*(T)^\eta\cap X^*(T)^+$  the following continuous homomorphism, if
   $\chi+\rho_G=(b_1,\ldots,b_n;b_0)\in (\frac 12 \Z)^{n+1}$:
   $$\xi_\chi:\C^*\to \hat G=\GL_n(\C)\times \C^*,\quad z \mapsto
   diag(z^{b_1}\bar z^{b_n},z^{b_2}\bar z^{b_{n-1}},\ldots, z^{b_n}\bar z^{b_1}; |z|^{2 b_0}).$$
   For a characteristic set $S$ we let $\epsilon_S=1$ if $S\subset \Z$ (i.e. if $S$ is of type $D_g$ or $C_g$) and $\epsilon_S=-1$ if $S\subset \frac 12 +Z$ (i.e. if $S$    is of type $B_g$).
   We have $b_i+b_{n+1-i}=b_0$ and $b_i=\frac{b_0}2+\beta_i$ with $\beta_{n+1-i}=-\beta_i\in S_\chi$. Thus
   $(-1)^{(b_i,b_{n+1-i})}=(-1)^{b_i-b_{n+1-i}}=(-1)^{2\beta_i}=\epsilon_{S_\chi}$ and
   $\xi_\chi(-1)=(\epsilon_{S_\chi}\cdot Id_n, 1)\in \GL_{n}\times \G_m$.
   \klabsatz
   Then $\xi_\chi$ may be extended to a homomorphism also called $\xi_\chi$ from $W_\R$ to $\hat G$ in such a way that
   $\xi_\chi(\sigma)$ is of the form $(J_\chi,1)\in \GL_n(\C)\times \C^*$ where $J_\chi=(a_i\cdot \delta_{i,n+1-j})$ is an antidiagonal
   matrix, whose entries satisfy $$a_i\cdot a_{n+1-i}=\epsilon_{S_\chi}.$$
   The conjugacy class of $\xi_\chi$ does not depend on the special choice of the $a_i$ in the case that $n$ is even.
   For odd $n=2g+1$  the conjugacy class only depends on $a_{g+1}$.

   \begin{lemma}\label{Unterscheidungslemma} Let $n=2g$ in the situation above. There exists a conjugate of $\xi_\chi$ which factorizes through
      $\hat H= GSp_{2g}(\C)$ if and only if $S_\chi$ is of type $B_g$.
      There exists a conjugate of $\xi_\chi$ which factorizes through
      $\hat H= GSO_{2g}(\C)$ if and only if $S_\chi$ is of type $D_g$.
    \end{lemma}
    Proof:  With respect to the standard basis
    $(e_i)_{1\le i \le 2g}$ the standard symplectic form $b_a$ satisfies $b_a(e_i,e_j)=\varepsilon_i\cdot \delta_{i+j,2g+1}$
    with $\epsilon_i=-\epsilon_{2g+1-i}\in \set{\pm 1}$, i.e. $\epsilon_i\cdot \epsilon_{2g+1-i}=-1$ for all $i$.
    If $S_\chi$ is of type $B_g$ then $b_0$ is odd and $\xi_\chi$ respects this form:
    \begin{align*} b_a(\xi_\chi(z)e_i,\xi_\chi(z)e_j)\;&=\; 0 \;=\; |z|^{2b_0}\cdot b_a(e_i, e_j)\quad\text{ for } i+j\ne 2g+1\\
    b_a(\xi_\chi(z) e_i,\xi_\chi(z) e_{2g+1-i})\;&= \; b_a(z^{b_i}\bar z^{b_{2g+1-i}} e_i, z^{b_{2g+1-i}}\bar z^{b_i} e_{2g+1-i})
    \,=\,|z|^{2 b_0} b_a(e_i,e_{2g+1-i}),\\
    b_a(\xi_\chi(\sigma)e_i,\xi_\chi(\sigma)e_j)\;&=\; 0 \;=\;  b_a(e_i, e_j)\quad\text{ for } i+j\ne 2g+1\\
    b_a(\xi_\chi(\sigma) e_i,\xi_\chi(\sigma) e_{2g+1-i})\;&= \; b_a(a_{2g+1-i} e_{2g+1-i}, a_i  e_{i})
    \;=\;b_a(e_i,e_{2g+1-i}) ,\\
        \end{align*}
    since $a_{2g+1-i}a_{i}=-1$ and $b_a(e_{2g+1-i},   e_{i})=-b_a(e_i, e_{2g+1-i})$.
    Similarly, if $S_\chi$ is of type $D_g$, then  $\xi_\chi$ respects the standard symmetric form $q_s(e_i,e_j)=\delta_{i+j,2g+1}$.
    \absatz
    Conversely, if there exists a conjugate of $\xi_\chi$ which factorizes through $\hat H= GSp_{2g}(\C)$, then
    $\xi_\chi$ respects some non degenerate alternating form $q'_a$. For $i+j\ne 2g+1$ we have
    $z^{b_i}\bar z^{b_{2g+1-i}} \cdot z^{b_{j}}\bar z^{b_{2g+1-j}}\ne |z|^{2 b_0} $ for some $z\in\C^*$ and therefore
    the invariance under $\xi_\chi(z)$ implies
    $q'_a(e_i,e_j)=0$. Since $q'_a$ is non degenerate we get $q'_a(e_i,e_{2g+1-i})=-q'_a(e_{2g+1-i},e_i)\ne 0$.
    The invariance under $\xi_\chi(\sigma)$ then implies $-1=a_{2g+1-i}a_{i}=\epsilon_{S_\chi}$,  and $S_\chi$ must be of
    type $B_g$.
    Finally if there exists a conjugate of $\xi_\chi$ which factorizes through
      $\hat H= GSO_{2g}(\C)$, then the same reasoning implies that $S_\chi$ is of type $D_g$.
    \qed

    \Absatz\Numerierung Let $\ZZZ(\gGg)$ be the center of the universal enveloping algebra
    of $\gGg$, which is isomorphic to the $W$-invariant part of the symmetric algebra $\SsS(\ttt)$  by Harish-Chandra's isomorphism
    \cite[23.3.]{Humphreys}.
    Here $\ttt=Lie(T)\otimes  \C=X_*(T)\otimes \C.$ Thus every $\chi\in X^*(T)\otimes \C$ gives rise to a
    homomorphism from $\ZZZ(\gGg)=\SsS(\ttt)^W$ to $\C$. Recall that each irreducible $\gGg$-module gives rise to a character
    of $\ZZZ(\gGg)$.

    \Absatz For simplicity put $\tilde K_\infty=K_\infty\cdot Z_\infty$, where $K_\infty\subset G(\R)$ is maximal under the connected compact subgroups and where $Z_\infty$ denotes the real connected component of the center of $G(\R)$.

    \begin{Def} Let  $V_\chi$ be the irreducible algebraic representation of $G$ of highest weight $\chi\in X^*(T)^+$. An irreducible admissible $(\gGg,\tilde K_\infty)$-module $\pi_\infty$ is called cohomological with respect to $V_\chi$ if and only if there exists
    some $i\ge 0$ such that $H^{i}(\gGg,\tilde K_\infty, V_\chi\otimes \pi_\infty)\ne 0$.
    An irreducible admissible representation of $G(\R)$ is called cohomological with respect to $V_\chi$ if and only if its associated
    $(\gGg,\tilde K_\infty)$-module of $\tilde K_\infty$-finite vectors is cohomological with respect to $V_\chi$.
    An irreducible automorphic representation $\pi=\hat\otimes \pi_v$ of $G(\A)$ is called cohomological with respect to $V_\chi$ if and only if  $\pi_\infty$ is cohomological with respect to $V_\chi$.
    \end{Def}

    \begin{lemma} Let  $V_\chi$ be the irreducible algebraic representation of $G$ of highest weight $\chi\in X^*(T)^+$.
    Let $\pi_\infty$ be an irreducible admissible  $(\gGg,K_\infty)$-module which is cohomological with respect to $V_\chi$.
    \begin{itemize}
        \item[(a)]  $\ZZZ(\gGg)$ acts by the character $\chi+\delta_G$ on $V_\chi$;
        \item[(b)] $\ZZZ(\gGg)$ acts by the character $-(\chi+\delta_G)$ on $\pi_\infty$;
        \end{itemize}
    \end{lemma}
    Proof: (a) This is well known (\cite[23.3.]{Humphreys}).
    \klabsatz
    (b) This is a consequence of  Wigners Lemma (\cite[Cor. I.4.2]{BorelWallach}).
    \qed
       \begin{lemma} For an irreducible admissible representation $\pi_\infty$ of $G(\R)$ let $\phi=\LLL(\pi_\infty)$ be some associated
    homomorphism $\phi: W_\R\to \hat G$ in the sense of \cite{Langlandsreell}.
    \begin{itemize}
        \item[(a)]  If $G=\G_m$ and $\pi_\infty=\sigma_{s,\epsilon}$ with $\sigma_{s,\epsilon}(r)=|r|^s\cdot (sign(r))^{\epsilon}$ with $s\in \C$ and $\epsilon\in \Z/2\Z$
        then we have $\phi(z)=(z\bar z)^s$ and $\phi(\sigma)=(-1)^\epsilon$.
        \item[(b)] For $G=\GL_n$ we have $\LLL(\omega_{\pi_\infty})= \det\circ \LLL(\pi_\infty)$.
        \item[(c)] For $G=\GL_n\times\G_m$ let $\LLL(\pi_\infty\times\rho_\infty)= \xi_\chi$ in the notations of \ref{Weilgruppenhom}.
            Then we have $\rho_\infty(r)=|r|^{b_0}$ and $\omega_{\pi_\infty}(-1)=(-\epsilon_{S_\chi})^{g}$ for even $n=2g$ and
            $\omega_{\pi_\infty}(-1)=(-1)^{g}\cdot a_{g+1}$ for odd $n=2g+1$.
        \item[(d)] If $\phi(z)=z^\mu\cdot \bar z^\nu$ for $z\in \C^*\subset W_\R$ and
        $\mu,\nu\in X_*(\hat T)\otimes \C=X^*(T)\otimes \C$ then $\pi_\infty$ has infinitesimal character $\mu$.
        \end{itemize}
    \end{lemma}
    Proof: (a) may be deduced from \cite[2.8.]{Langlandsreell} and (b) from \cite[p.30 condition (ii)]{Langlandsreell}. (c) is a consequence of  (a) and (b),
    the explicit description of $\xi_\chi$ in \ref{Weilgruppenhom} and the fact that the determinant of the standard antidiagonal matrix of size $2g$ resp. $2g+1$ is $(-1)^g$.
    \medskip
    (d) is known for discrete series representations for all $G$: see \cite[Theorem 9.20]{Knapp} and observe that the characterization
    of discrete series representations via characters of maximal anisotropic tori in  \cite[Theorem 12.7]{Knapp} is the same as in \cite{Langlandsreell}.
        Then one checks immediately that the statement is invariant under parabolic induction and under forming Langlands quotients.
    \qed

\Absatz\Numerierung For an irreducible admissible representation $\rho$ of $\GL_m(F_v)$ resp. a cuspidal
    automorphic representation $\pi$ of $\GL_m(\A)$ and a natural number $n\ge 1$ we denote by $\pi=MW(\rho,n)$ the irreducible admissible representation of $\GL_{mn}(F_v)$ resp. of $\GL_{mn}(\A)$, which is the unique irreducible (Langlands) quotient of
    $$\pi|\cdot|^{\frac{n-1}2}\times \pi|\cdot|^{\frac{n-3}2}\times\ldots \times \pi|\cdot|^{\frac{1-n}2}.$$
    Thus $MW(\rho,1)=\rho$. Recall the main result of \cite{MoglinWald} that the discrete spectrum
    $L_{2,disc}(\GL_N(\Q)\bs \GL_N(\A))$ is the Hilbert direct sum of all $MW(\rho,n)$ where
    $n$ runs over all divisors of $N$ and $\rho$ over all cuspidal automorphic representations of $\GL_m(\A)$ where $N=mn$.

    \begin{lemma} Let $m,n\ge 1$ be natural numbers. Let $G_m=\GL_m\times \G_m$ and $\gGg_m$ the Lie algebra of $G_m(\R)$.
    Let $T_m\subset G_m$ denote the standard maximal split torus.
     \begin{itemize}
        \item[(a)]  For a cuspidal automorphic representation $\rho=\hat\otimes_v \rho_v$ of $\GL_m(\A)$ we have
        $MW(\rho,n)=\hat\otimes_v MW(\rho_v,n)$.
        \item[(b)] If $ \rho_\infty$ is cohomological with respect to $V_\chi$ for some $\chi\in X^*(T_m)^\eta\cap X^*(T_m)^+$, then we have:
        $ MW(\rho_\infty,n)$ is cohomological with respect to to $V_\sigma$ for
        $\sigma\in  X^*(T_{mn})^+$ if and only if $\sigma$ has the same weight as $\chi$ and if we have the identity
        $S_\sigma=MW(S_\chi,n)$ between the characteristic sets.
        \end{itemize}
    \end{lemma}
    \qed
   \Absatz\Numerierung {\sc More representations.}\label{tempered}
   For irreducible representations $\pi_1,\ldots,\pi_k$ of $\GL_{n_i}(F)$ we denote by
   $\pi_1\times\ldots\times \pi_k$ the normalized parabolically induced representation of $\GL_n(F)$, where
   $n=n_1+\ldots+n_k$.
   \absatz
   For $s\in \C$ and an integer $p>0$ let $D(s,p)$ be the discrete series representation of $\GL_2(\R)$, which is
   the infinite dimensional submodule of $\sigma_{\frac{s+p}2,p+1}\times \sigma_{\frac{s-p}2,0}$.
   Note that the central character of $D(s,p)$ is $\sigma_{s,p+1}$.
   Viewed as a representation of $\SL_2(\R)^{\pm}=\set{A\in \GL_2(\R)| \det(A)=\pm 1}$ it is independend of $s$ and will be denoted by $D(p)$. As a representation of $\SL_2(\R)$ id decomposes $D(p)=D(p)^+\oplus D(p)^-$, where in $D(p)^+$
   the lowest positive $\SO_2(\R)$-weight is $p+1$.
   For $w\in \Z$  we have $\LLL(D(w,p)\otimes \sigma_{w,0})=\xi_\chi$ for $\chi+\delta_{\GL_2}=(\frac{w+p}2,\frac{w-p}2;w)$.
   \absatz
   If $n=2g$ and $\chi+\delta_G=(b_1,\ldots,b_{2g};b_0)$ is as in \ref{Weilgruppenhom}, where $G=\GL_n\times \G_m$, then we may form the
   tempered representation
   $$ D_\chi \quad=\quad D(b_0, b_1-b_{2g})\times\ldots\times D(b_0,b_g-b_{g+1})\otimes \sigma_{b_0,0}
   \qquad\text { of } G(\R).$$
   We have $\LLL(D_\chi)=\xi_\chi$, and $D_\chi$ is the only tempered representation  of $\GL_{2g}(\R)$ which is cohomological with respect to $V_\chi$ (\cite[Thm 4.12.(3)]{Grobner}).

   \absatz
   If $n=2g+1$ and $\chi+\delta_G=(a_g,\ldots,a_1,0,-a_1,\ldots,-a_g)$, where $G=\PGL_{2g+1}$ then we may form
   the tempered representation
   \begin{align}\label{oddtempered}
    D_\chi \; = \; D(S_\chi)\quad :=\quad D(0,2 a_g)\times\ldots\times D(0,2 a_1)\times \sigma_{0,g}.
    \end{align}
   As a representation of $\GL_{2g+1}(\R)$ it has trivial central character and is thus a representation of $\PGL_{2g+1}(\R)$.
   Again $\LLL(D_\chi)=\xi_\chi$.

    \begin{prop}\label{LefschetzC} Let $G=\PGL_{2g+1}$  and $S$ be an $n$-admissible characteristic set of type $C_m$, where $2g+1=n\cdot (2m+1)$.
    Let $V_\chi$ be an irreducible representation of $G$ with  $S_\chi=MW(n,S)$ .
    Let $\pi_{\infty,1}=D(S)$ be the tempered representation
    of $\PGL_{2m+1}(\R)$ as in (\ref{oddtempered}) and $\pi_\infty= MW(n,\pi_{\infty,1})$. Then we have:
          $$ tr(\eta, H^*(\gGg,K_\infty,V_\chi\otimes \pi_\infty))\; =\; \pm 2^g .$$
    \end{prop}
Proof:
    Since the kernel $\mu_{2g+1}$ of the projection $\SL_{2g+1}\to G=\PGL_{2g+1}$
    satisfies $\mu_{2g+1}(\R)=\{1\}$ and $H^1(\R,\mu_{2g+1})=\{1\}$, we get an isomorphism for the real points
    $G(\R)\cong\SL_{2g+1}(\R)$,  and consequently  $K_\infty=\SO_{2g+1}$. Furthermore $\gGg= \mathfrak{sl}_{2g+1}$.
    Recall that $\eta(g)=J\cdot \theta(g)\cdot J^{-1}$, where $J=J_{2g+1}\in \SO_{2g+1}$ is anti-diagonal and
    satisfies ${}^tJ=J=J^{-1}$, and where
    $\theta: g\mapsto \tran g^{-1}$ is the Cartan involution.  Let
    $\kkk=\mathfrak{so}_{2g+1}=\left\{ A\in \mathfrak{sl}_{2g+1}\;| \; \tran A=- A\right\}$ and
    $\ppp= \left\{ A\in \mathfrak{sl}_{2g+1}\;| \; \tran A= A\right\}$ be the eigenspaces of $\theta$ in $\gGg$.
    Assume that a  maximal Cartan subalgebra $\ttt\subset\kkk$ is chosen in such a way,
    that $Ad(J)$ acts as identity on $\ttt$, i.e.
    $\ttt\subset \kkk^+=\{k\in \kkk| Ad(J)(k)=k\}$.
    In fact if we decompose  $\R^{2g+1}$  into the $\pm 1$ eigenspaces under the involution $J$ (of dimensions $g$ and $g+1$),
    then $A\in End(\R^{2g+1})$ is fixed under $Ad(J)$ if it respects this decomposition into eigenspaces. Thus
    $\kkk^+\cong \mathfrak{so}_{g}\times \mathfrak{so}_{g+1}$ and we can easily find a maximal torus $\ttt$ of rank $g$ inside
    $\kkk^+$.
    \klabsatz
    Then we have
    $\eta(x)=x$ for each $x\in\ttt$. Let $\gGg_\C=\gGg\otimes_\R \C= \mathfrak{sl}_{2g+1}(\C)$ be the complexification
    of $\gGg$. For a fixed $x\in \ttt$  let $\lLl=\left\{ A\in \gGg\;| \; [x,A]=0\right\}$, let $\lLl_\C=\lLl\otimes\C$  and let
    $\uuu_\C$ be the sum of all positive eigenspaces of $ad(ix)$ in $\gGg_\C$. Then $\qqq_\C=\lLl_\C\oplus \uuu_\C$ is a $\theta$-stable parabolic subalgebra of $\gGg_\C$ in the sense of \cite{Vogan}.
   Now let us assume that  $\pi_\infty= A_{\qqq}(\lambda)$ in the notations of \cite[Thm. 5.3.]{Vogan}, which will be proved later.
    Recall the isomorphism of \cite[Thm. 5.5.]{Vogan}:
    $$ H^{i}(\gGg,K_\infty,V_\chi\otimes \pi_\infty)\cong H^{i-R}(\lLl_\C,\lLl_\C\cap\kkk_\C,\C)
    \cong Hom_{\lLl_\C\cap\kkk_\C}(\Lambda^{i-R}(\lLl_\C\cap\ppp_\C),\C),$$
    where $R=\dim(\uuu_\C\cap\ppp_\C)$, which is given by a cap product with an element of $\Lambda^R(\uuu_\C\cap\ppp_\C)$
    (\cite[Proof of Thm 3.3]{Vogan}).
    Observe that $\kkk,\ppp,\lLl,\uuu_\C,\qqq_\C$ are all $\eta$-stable, since we have $\eta(x)=x$ and since $\eta=Ad(J)\circ \theta$ , where $J\in K_\infty$.
    If $\epsilon\in\set{\pm 1}$ denotes the eigenvalue of $\eta$ on $\Lambda^R(\uuu_\C\cap\ppp_\C)$, we immediately get
    $$tr(\eta, H^*(\gGg,K_\infty,V_\chi\otimes \pi_\infty))\; =\; \epsilon\cdot tr(\eta, H^*(\lLl_\C,\lLl_\C\cap\kkk_\C,\C)).$$
    Now let $L^K\subset K_\infty$ be the connected Lie subgroup with real Lie algebra $\lLl\cap \kkk$, and let
    $L^c$ be a compact Lie group containing $L^K$ with real Lie algebra $i(\lLl\cap \ppp)\oplus (\lLl\cap \kkk)\subset\gGg_\C$.
     Thus the complexified Lie algebra of $L^c$ is  $\lLl_\C$.
    We may define these groups as connected components of some stabilizers
    \begin{align*}
     L^K\quad&=\quad \left\{ A\in K_\infty | Ad(A)(x)=x\right\}^\circ\\
     L^c\quad&=\quad \left\{ A\in G(\C) | \bar A= \theta(A),\; Ad(A)(x)=x\right\}^\circ
    \end{align*}
    We will see later that the stabilizers are already connected in the examples we consider, and then $Ad(J)x=x$ implies
    $J\in L^K$.

    By \cite[Theorem V in 4.19.]{Greub2}  we have an isomorpism
    $$ H^{i}(\lLl_\C,\lLl_\C\cap\kkk_\C,\C) \cong H^{i}(L^c/L^K, \C).$$
    There exists involutions $\eta$ and  $\theta$ on the Lie group $L^c$, which induces the involutions with the same name on
    the Lie algebra, and the isomorphisms above are equivariant with respect to these actions.
    Since $\eta$ and the Cartan involution $\theta$ differ by conjugation with $J$,
    which lies in the path connected group $L^K$, they are
    homotopic as automorphisms of $L^c$ and thus induce the same map on cohomology.
    Since the Cartan involution acts by $-1$ on $\lLl\cap \ppp$, the Lefschetz number of $\eta$  is
    the sum of the dimensions of the cohomology groups $H^{i}(L^c/L^K, \C)$:
    $$  tr(\eta, H^*(\gGg,K_\infty,V_\chi\otimes \pi_\infty))\; =\; \epsilon\cdot Betti(L^c/L^K), $$
    where we denote by $Betti(X)=\sum_{i} \dim\left( H^{i}(X, \C)\right)$  the sum of all Betti numbers of a topological space $X $.
    \absatz
    Recall $2g+1=(2m+1)\cdot n$. Now we  fix $x\in \ttt$ in such a way that it is conjugate (under $K_\infty$)
     to a block diagonal $\tilde x$ consisting of blocks
    $$ \begin{pmatrix} 0_n& x_1 1_n\\ -x_1 1_n & 0_n\\ \end{pmatrix},\ldots ,
    \begin{pmatrix} 0_n& x_m 1_n\\ -x_m 1_n & 0_n\\ \end{pmatrix}, 0_n$$
    with $0<x_1<\ldots< x_m$.
    Then $\lLl$ is isomorphic to the algebra $\tilde\lLl$ of all block diagonal matrices of trace $0$ with blocks of the form
    \begin{align} \label{Blockmatrix}
     \begin{pmatrix} B_1& C_1\\ -C_1 & B_1\\ \end{pmatrix},\ldots ,
     \begin{pmatrix} B_m& C_m\\ -C_m & B_m\\ \end{pmatrix}, B_0,
    \end{align}
    where $B_0,\ldots,B_m,C_1,\ldots,C_m$ are real $n\times n$ matrices.
    We have $\lLl\cap\kkk\cong \tilde \lLl\cap \kkk$ and the latter is characterized by the conditions
    $\tran B_j=- B_j$ and $\tran C_j= C_j$ for all $j$ i.e. $B_0\in \mathfrak{so}_{n}$ and $B_\nu+iC_\nu\in \mathfrak{u}_n$ for $\nu=1,\ldots,m$.

     Thus we have $L^K= (U_n)^m\times \SO_n$, and the same computations imply that this $L^K$ is the full stabilizer
       of $x$ in $K_\infty$, so that $J\in L^K$. Similarly we get  $L^c= S((U_n\times U_n)^{m}\times U_n)$, where we write
    $S(H)=H\cap \SL_{2g+1}$ for a subgroup $H\subset \GL_{2g+1}$. The embedding $L^K\hookrightarrow L^c$ is of the form
    $U_n\ni A\mapsto (A,\bar A)\in U_n\times U_n$ in each factor.
    Thus $L^c/L^K = S(U_n^m\times (U_n/\SO_n))$.
    Using the finite covering $U_1\times S(U_n^m\times (U_n/\SO_n)) \to U_n^m\times (U_n/\SO_n)$
    and the Kunneth formula
    we get $Betti( L^c/L^K)=\frac 12\cdot Betti(U_n^m\times (U_n/\SO_n))$.
       Now $U_n$ has the same cohomology as $S^1\times S^3\times \ldots\times S^{2n-1}$, (\cite[Thm. IX in 6.27]{Greub3}),
     and since $n=2\nu+1$ is odd, the homogeneous space $U_n/\SO_n$ has the same cohomology as
     $S^1\times S^5\times S^9\times \ldots S^{4\nu+1}$ ( \cite[Prop. 31.4]{BorelKohom},
     compare \cite[Table I after 11.16]{Greub3}). Thus $Betti(U_n)=2^n$ and $Betti(U_{2\nu+1}/SO_{2\nu+1})=2^{\nu+1}$,
     and the Kunneth formula finally implies
     $tr(\eta, H^*(\gGg,K_\infty,V_\chi\otimes \pi_\infty)) = \epsilon\cdot \frac 12\cdot (2^n)^m\cdot 2^{\nu+1} = \epsilon\cdot 2^g$, since we have
    $g= mn+\nu.$

    \absatz
    We still have to prove  $\pi_\infty= A_{\qqq}(\lambda)$ in the sense of \cite[Thm. 5.3.]{Vogan}:
    $\tilde \lLl$ is the Lie algebra of the group $\tilde L=S(\Pi_{j=1}^m \GL_n(\C) \times \GL_n(\R))$:
    A block diagonal matrix as in (\ref{Blockmatrix})  belongs to $\tilde L$
     if and only if $\prod_{j=1}^m |\det(B_j+C_j)|^2\cdot \det(B_0)=1.$  Let $H$ be the Cartan subgroup of all such
    block diagonal matrices where each $B_j$ and each $C_j$ is diagonal and let
    $\hhh$ be its Lie algebra. Thus $H\cong S((\C^*)^{mn}\times (\R^*)^n)$. Write
    $$S=\{ -a_m,\ldots,-a_1,0,a_1,\ldots,a_m\} \text{ with }0=a_0<a_1<\ldots<a_m.$$ Since $S$ is $n$-admissible we have
    $a_j\ge a_{j-1}+n$ for all $j=1,\ldots,m$. Thus the numbers $d_j=a_j-j\cdot n$ satisfy $0=d_0\le d_1\le \ldots\le d_m$.
    Then we can write $\chi\in X^*(T)^\eta\subset\Z^{2g+1}$ in the form
        $$\chi =(d_m,\ldots,d_m,\ldots,d_1,\ldots,d_1,0,\ldots,0,-d_1,\ldots,-d_1,\ldots,-d_m,\ldots,-d_m),$$
    where each entry is repeated $n$ times.
    With the choice of a Cartan subgroup $H$ as above and a suitable choice of the positive roots for $H$ in $\gGg_\C$
    it is clear that $\chi$ (the highest weight of $V_\chi$) is the restriction of the following unitary character
    $\lambda$ of $\tilde L$:
    $$\lambda: \tilde L\ni (X_1,\ldots,X_m,B_0)\quad \mapsto \quad \prod_{j=1}^m (\det(X_j)/\overline{\det(X_j)})^{d_j}
    \qquad\text{ for }X_j\in \GL_{n}(\C).$$
    Now $A_{\qqq}(\lambda)$ is the representation defined in \cite[Thm. 5.3.]{Vogan}. To prove that it equals
    $\pi_\infty$ we use the description in \cite[Thm. 6.16.]{Vogan}: We decompose $H= T^+\cdot A^d$, where
    $A^d$ consists of all block matrices $B=B(b_{k,j})$ as in (\ref{Blockmatrix}) such that $C_j=0$ and $B_j=diag(b_{1,j},\ldots,b_{n,j})$ with
    $b_{k,j}>0$, and where $T^+\cong (S^1)^{nm}$ is compact.
    The unipotent subgroup $N^L$ consists of all block matrices (\ref{Blockmatrix}), where each $B_j$ is unipotent upper triangular
    and each $C_j$ is nilpotent upper triangular. Since the restriction of $\lambda$ to $A^d$ is trivial, the character
    $\nu^d$ of \cite{Vogan} is half the sum of all roots of $A^d$ on $\mathfrak{n}^L$ and this is the character sending
    the  $B=B(b_{k,j})$ to
    $$\prod_{j=1}^m\left(\prod_{k=1}^n\; b_{k,j}^{2k-n-1}\right)\cdot \prod_{k=1}^n b_{0,k}^{k-\frac{n+1}2}.$$
    The centralizer of $A^d$ is of the form $A^d\cdot M^d$ where $M^d\cong (\SL_2^{\pm}(\R))^{mn}\times \{\pm 1\}^{n}.$
    On $M^d$ we take the discrete series representation $\sigma$ which is the tensor product of the representations $D(2 d_j)$
    of $\SL_2^\pm (\R)$ on the factor with indices $j\in \set{1,\ldots,m}$ and $k\in\set{1,\ldots,n}$
    and of the representation $sign^{g}$ on $\set{\pm 1}$.
    \klabsatz
        The group $A$ of \cite[(6.8.)]{Vogan} consists of all block matrices $B$ as above, where $B_0=B_1=\ldots=B_m$ is a diagonal
    matrix with positive entries. Its centralizer in $\GL_{2g+1}$ is of the form $(\GL_{2m+1})^n$, and therefore
    it has Langlands decomposition $AM$ with $M\cong (\SL_{2m+1})^n$. The one dimensional representation
    $\tilde \nu^d: (\GL_{2m+1})^n \to \R_{>0},\; (X_1,\ldots,X_n)\;\mapsto \;\prod_{k=1}^n |\det(X_k)|^{k-\frac{n+1}2}$ restricts
    to the character $\nu^d$ on $A^d$.
    Now let $N^d\subset G$ be the unipotent radical of some parabolic $P^d$ which has $M^d\cdot A^d$ as  Levi factor and let
    $N$ be the unipotent radical of some parabolic $P\supset P^d$, which has $MA$ as Levi factor.
    Then  $A_\qqq(\lambda)$ is the Langlands quotient of the representation $Ind_{M^dA^dN^d}^G (\sigma\otimes \nu^d\otimes 1)$.
    But by induction by stages this is the Langlands quotient of $Ind_{MAN}^G \tilde\pi_{\infty}$, where
    $\tilde\pi_{\infty}=(\pi_{\infty,1}\times\ldots\times \pi_{\infty,1})\otimes \tilde \nu^d$ and
    $\pi_{\infty,1}=D(S)$ is the unitarily induced (tempered) representation as in (\ref{oddtempered}). Therefore
    $A_\qqq(\lambda)=MW(n,\pi_{\infty,1})=\pi_\infty$.

\qed

\Absatz\Numerierung \label{pm Eigenraum} Recall that in the case $G=\GL_{2g}$ we have $\tilde K_\infty=\SO_{2g}\cdot Z_\infty$ with
     $Z_\infty$ the connected component of the center of $G(\R)$. Put $\tilde K_\infty^m=O_{2g}\cdot Z_\infty$.
     Then $\tilde K_\infty^m$ acts on $Hom_{\tilde K_\infty}(\Lambda^*(\gGg/\tilde \kkk),V_\chi\otimes \pi_\infty)$ by $(k\phi)(l) = k(\phi(k^{-1}l))$.
     Since $\tilde K_\infty$ acts trivially we get an action of $\tilde K_\infty^m/\tilde K_\infty=O_{2g}/SO_{2g}=\Z/2\Z$,
      and we may decompose
     the cohomology into the $+1$ and $-1$ eigenspaces $H^*(\gGg,\tilde K_\infty,V_\chi\otimes \pi_\infty)^\epsilon)$ (here $\epsilon=\pm$)
     under the  action of the non trivial element of $\tilde K_\infty^m/\tilde K_\infty$.
          Since $\tilde K_\infty^m$ is stable under $\eta$ and since $\eta$ acts as identity on $\tilde K_\infty^m/\tilde K_\infty$, the action of
     $\eta$ respects the decomposition of $H^*(\gGg,\tilde K_\infty,V_\chi\otimes \pi_\infty)$ into these eigenspaces.

    \begin{prop}\label{LefschetzBD} Let $G=\GL_{2g}$  and let $S$ be an $n$-admissible characteristic set of type $B_m$ or $D_m$, where $g=n\cdot m$.
    Let $V_\chi$ be the irreducible representation of $G$  with describing set $MW(n,S)$ and of weight $w$. Let $\pi_{\infty,1}$ be a tempered representation of $\GL_{2m}$ of type $S$ with central character of weight $-w$,
    and let $\pi_\infty= MW(n,\pi_{\infty,1})$. Then we have:
          $$ tr(\eta, H^*(\gGg,K_\infty,V_\chi\otimes \pi_\infty)^\epsilon)\; =\; \pm 2^{g-1}\quad\text{ for }\epsilon=\pm .$$
    \end{prop}

    Proof: At first we note that $\pi_\infty$ is an induced representation:
    $$ \pi_\infty=Ind_{G(\R)^\circ}^{G(\R)} \pi_\infty^\circ,\quad \text{ where } G(\R)^\circ=\{g\in \GL_{2g}(\R)| \det(g)>0\}.$$
    This may be deduced from the analogous statement for the building constituents
    $D(p_i)$  of $\pi_{\infty,1}$.
    From this fact one concludes that we have $\eta$-equivariant isomorphisms for both choices of $\epsilon$:
    $$ H^*(\gGg,K_\infty,V_\chi\otimes \pi_\infty)^\epsilon\cong H^*(\gGg,K_\infty,V_\chi\otimes \pi_\infty^\circ).$$
    Now $Z_\infty$ acts trivially on $\Lambda^*(\gGg/\kkk)$ and by assumption also on $V_\chi\otimes \pi_\infty^\circ$, so that
    equivariance under $Z_\infty$ is automatically fulfilled. Therefore we may pass to the derived group situation:
    With $G^{der}=\SL_{2g}$, $\gGg^{der}=\mathfrak{sl}_{2g}$ and $K_\infty^{der}=\SO_{2g}$ we get:
    $$ tr(\eta, H^*(\gGg,\tilde K_\infty,V_\chi\otimes \pi_\infty)^\epsilon)=
    tr(\eta, H^*(\gGg^{der},K_\infty^{der},V_\chi\otimes \pi_\infty^\circ)).$$
     Now we can apply the results of \cite{Vogan}, and the rest of the proof is  analogous to the proof of the preceding proposition with the following modifications:
     Since $J_{2g}$ is skew symmetric, the Cartan subalgebra $\ttt$ can be chosen to be the algebra of all antidiagonal
     skew symmetric matrices. In the definition of $\tilde x$ we have to omit the last $0_n$ block, and therefore
     in (\ref{Blockmatrix}) the last entry $A_0$ disappears.
     Thus $L^K\cong (U_n)^m$ and $L^c\cong S((U_n\times U_n)^m)$ and therefore $L^K/L^c\cong S((U_n)^m)$ and
     $$tr(\eta, H^*(\gGg^{der},K_\infty^{der},V_\chi\otimes \pi_\infty^\circ))=\pm Betti(S((U_n)^m))=
     \frac {\pm Betti((U_n)^m)}2$$
     $$=\pm \frac{ (2^n)^m}2= \pm 2^{g-1}.$$
\qed

\meinchapter{Statement of the Main Theorem}\label{Ch Statement}

   \Absatz\Numerierung
   For $d\in \Q^*/(\Q^*)^2$ let $\SO_{2g}^d$ be the quasisplit form of $\SO_{2g}$ with splitting field $\Q(\sqrt d)$.
   Let $G_1$ be one of the groups $\SP_{2g}$ (case $C_g$), $\SO_{2g+1}$ (case $B_g$) or $\SO_{2g}^d$ (case $D_g$). Let $G=\PGL_{2g+1}$ in case $C_g$
   and $G=\GL_{2g}$ in the cases of orthogonal groups. Let $\iota: \hat G_1\hookrightarrow \hat G$ be the corresponding embedding of
   dual groups. Let $\chi_d:\A^*/\Q^*\to \set{\pm 1}$ be the quadratic character associated
   to the field extension $\Q(\sqrt d)/\Q$ by class field theory.

   \begin{Def} Let $\pi_1=\hat \otimes_v \pi_{1,v}$ be an irreducible automorphic representation of $G_1(\A)$ and $\pi=\hat \otimes_v \pi_v$
   be an irreducible automorphic representation of $G(\A)$.
   \begin{itemize}
   \item[(a)] $\pi$ is called a weak lift of $\pi_1$ if $\pi_v$ is the lift of $\pi_{1,v}$ for almost all finite places $v$, where $\pi_v$ and $\pi_{1,v}$ are unramified.
   \item[(b)]  $\pi$ is called a semi-weak lift of $\pi_1$, if $\pi$ is  a weak lift of $\pi_1$ and if
   $$\iota\circ\LLL(\pi_{1,\infty})=\LLL(\pi_\infty).$$
   \item[(c)] $(\pi,\pi_1)$ is called an elementary pair for $(G,G_1)$, if $\pi$ is a semi-weak lift of $\pi_1$, if $\pi_1$ is cuspidal and generic and if $\pi$ is cuspidal (and therefore generic) as well.
   \item[(d)] Let $S\subset \frac 12 \Z$ be a characteristic set of the same type as $(G,G_1)$.  An elementary pair $(\pi,\pi_1)$ is called cohomological of type $S$, if
        if in case $B_g$ or $C_g$ the representation $\pi_\infty$ is cohomological with respect to $V_\chi$ and if in case $D_g$ for every odd integer $w$ the representation $\pi_\infty\otimes |det|^{-w/2}$ is cohomological
       with respect to a coefficient system $V_\chi$ of weight $w$ with $S=S_\chi$.
   \end{itemize}
   \end{Def}
   Now we can state our main theorems
   \begin{thm} \label{Maintheorem}
   Let $\chi\in X^*(T_H)^+$ be a dominant weight for the group $H=\SP_{2g}\supset T_H$ and let
     $V_\chi$ be the associated highest weight module. Let $\tau=\tau_\infty\otimes \tau_f$ be an irreducible automorphic representation of $\Sp_{2g}(\A)$ such that
     $\tau_f$ appears with non trivial multiplicity in the alternating sum
     $$\sum_{i} (-1)^{i} H^{i}(\mathfrak{sp}_{2g}, U_g,L_{2,disc}(\SP_{2g}(\Q)\backslash \SP_{2g}(\A))\otimes V_\chi )$$
     considered as an element in the Grothendieck group. Then there exists a finite family of octuples
     $$(X_{\gamma_i}, G^{(i)}, G^{(i)}_1,n_i, d_i, S_i, \pi^{(i)},\pi_1^{(i)})_{i=1,\ldots,r}$$ such that for each $i=1,\ldots,r$ we have:
     \klabsatz
     $X_{\gamma_i}\in\set
     {B_{\gamma_i},C_{\gamma_i},D_{\gamma_i}}$ is a type,
     $(G^{(i)}, G^{(i)}_1)$ is an endoscopic pair of type $X_{\gamma_i}$, $n_i\in \N_{\ge 1}$,
      $d_i\in \Q^*/(\Q^*)^2$ a square class,  $S_i$ is an $n_i$-admissible characteristic set of type $X_{\gamma_i}$ and
     $(\pi^{(i)},\pi_1^{(i)})$ is an elementary pair for $(G^{(i)}, G^{(i)}_1)$, which is cohomological of type $S_i$,
     and such that we have:
     $$S_\chi\quad=\quad\bigcup_{i=1}^k MW(S_i,n_i)$$ is a partition in smaller characteristic sets and such that
     $$\tau \text{ lifts weakly to }\; MW(\tilde\pi^{(1)},n_1)\times\ldots\times MW(\tilde\pi^{(r)},n_r),\;\text{ where }$$
     \begin{itemize}
     \item
     If $X_{\gamma_i}= B_{\gamma_i}$, then $\tilde\pi^{(i)}=\pi^{(i)}$, $d_i=1$ and $n_i$ is even.
     \item If  $X_{\gamma_i}= C_{\gamma_i}$, then $\tilde\pi^{(i)}=\pi^{(i)}\otimes \chi_{d_i}$,  and $n_i$ is odd.
     \item If  $X_{\gamma_i}= D_{\gamma_i}$, then $\tilde\pi^{(i)}=\pi^{(i)}$, $(-1)^{\gamma_i}d_i>0$,
     $G_1^{(i)}= \SO_{2\gamma_i}^{d_i}$  and $n_i$ is odd.
     \end{itemize}
     We have furthermore: $g=\sum_{i=1}^r \gamma_i\cdot n_i$, $\Pi_{i=1}^r d_i=1$, $\omega_{\tilde \pi^{(i)}}=\chi_{d_i}$.
     \end{thm}
     \Absatz

     \begin{thm} \label{Umkehrtheorem}
   Let $\chi\in X^*(T_H)^+$ be a dominant weight for the group $H=\SP_{2g}\supset T_H$ and let
     $V_\chi$ be the associated highest weight module.
    Assume that there exists a finite family of octuples
     $$(X_{\gamma_i}, G^{(i)}, G^{(i)}_1,n_i, d_i, S_i, \pi^{(i)},\pi_1^{(i)})_{i=1,\ldots,r}$$
     satisfying the conditions stated in Theorem \ref{Maintheorem}.

     Then there exists an  irreducible automorphic representation $\tau=\tau_\infty\otimes \tau_f$ of $\Sp_{2g}(\A)$ such that
     $\tau_f$ appears with non trivial multiplicity in the alternating sum
     $$\sum_{i} (-1)^{i} H^{i}(\mathfrak{sp}_{2g}, U_g,L_{2,disc}(\SP_{2g}(\Q)\backslash \SP_{2g}(\A))\otimes V_\chi )$$
     and such that
     $$\tau \text{ lifts weakly to }\; MW(\tilde\pi^{(1)},n_1)\times\ldots\times MW(\tilde\pi^{(r)},n_r),$$
     where $\tilde\pi^{(i)}$ is related to $\pi^{(i)}$ as in theorem \ref{Maintheorem}.
     \end{thm}
     The proofs will be given in (\ref{ProofMT}) and (\ref{ProofUmkehr}).

\meinchapter{Comparison of Chevalley volumes }\label{Ch Chevalley}

\Absatz\Numerierung {\sc Chevalley elements and Chevalley volumes}
Let $\gGg$ be a semisimple split Lie algebra over $\R$   and
 $$ \gGg\quad =\quad \hhh\oplus\bigoplus_{\alpha} \gGg_\alpha$$ a decomposition
 into a Cartan subalgebra $\hhh\subset \gGg$ and the corresponding one dimensional root spaces $\gGg_\alpha$.
 We fix some Chevalley basis
 $Ch_{\gGg}=\{h_\alpha,\alpha\in\Delta\}\cup\{u_\alpha,\alpha\in\Phi\}$ of $\gGg$ in the sense of
 \cite[25.2]{Humphreys}, where $\Delta$ is the set of simple roots with respect to some Borel subalgebra $\bbb\supset\hhh$.
 Note that this expecially implies the existence of elements $h_\alpha\in \hhh$ for all $\alpha\in \Phi$
  such that $\gGg_\alpha=\R\cdot u_\alpha$ and
 $$ [h_\alpha,u_\alpha]=2 u_\alpha;\quad [h_\alpha,u_{-\alpha}]=-2 u_{-\alpha},\quad [u_\alpha,u_{-\alpha}]=h_\alpha$$
 and such that the $h_\alpha$ for $\alpha\in\Phi$ lie in the $\Z$-lattice spanned by the $h_\alpha$ for $\alpha\in \Delta$.
 It is clear that the elements $h_\alpha\in\hhh$ and $u_\alpha\wedge u_{-\alpha}\in \Lambda^2\gGg$ are uniquely characterized by these conditions.
 From this we deduce that the Chevalley element
 $$ \lambda_{Ch}=\bigwedge_{\alpha\in\Delta} h_\alpha \; \wedge \;\bigwedge_{\alpha\in\Phi}\quad\in\quad\Lambda^{\dim\gGg}\gGg$$
 (where we do not care about the exact order in the wedge product, so that it is only well defined up to sign) does not depend on
 the exact choice of the generators $u_\alpha\in\gGg_\alpha$. Since all pairs $(\hhh,\bbb)$ are conjugate under inner automorphisms of
 $\gGg$ and since inner automorphisms act as identity on the maximal exterior power $\Lambda^{\dim\gGg}\gGg$, it is clear that
 $\lambda_{Ch}$ is unique up to sign and does not depend on the choices in its construction.
 It is furthermore clear that to compute $\lambda_{Ch}$ we just have to check the commutation relations above and do not need to fix
 the $u_\alpha$ in such a way that the further conditions in \cite[25.2]{Humphreys} are satisfied (commutation relations between
 $u_\alpha$ and $u_\beta$ in the case that  $\beta\ne \pm\alpha$), since we only need  $u_\alpha\wedge u_{-\alpha}$
 for our computation.

 \Absatz
 If $G/\R$ is a semisimple algebraic group with Lie algebra $\gGg$ and if $G_K/\R$ is an inner form of $G$ then
 we have a canonical isomorphism $\Lambda^{\dim G}Lie(G_K)\cong\Lambda^{\dim G}\gGg$ (since $Ad$ acts trivially on
  $\Lambda^{\dim G}\gGg$) and thus get a Chevalley element
 $\lambda_{Ch}\in \Lambda^{\dim G}Lie(G_K)$, which again is unique up to sign. Then we define a Chevalley form to be an $G_K$-invariant
 differential form of top degree $\omega_{Ch}$ on $G_K$ which is the dual of the Chevalley element $\lambda_{Ch}$ at the
 identity of $G_K$.
 If $G_K$ is compact then we define the Chevalley volume of $G_K$ (and of its inner form $G$) to be
 $$vol_{Ch}(G_K)\quad:=\quad vol_{Ch}(G)\quad:=\quad \left| \int_{G_K}\omega_{Ch} \right|.$$
 Due to the absolute value this is completely well defined, i.e. does not depend on the sign of the Chevalley element and of the orientation of $G_K$.

 \Absatz
 \begin{prop}\label{BC Chevalleyvolume}
 For all $g\ge 1$ we have $$vol_{Ch}(Sp_{2g})\;=\; vol_{Ch}(Spin_{2g+1})\;=\; 2\cdot vol_{Ch}(SO_{2g+1}).$$
 \end{prop}
 Proof: The equation $vol_{Ch}(Spin_{2g+1})\;=\; 2\cdot vol_{Ch}(SO_{2g+1})$ is an immediate consequence of the isogeny of
 degree $2$ from $Spin_{2g+1}$ to $SO_{2g+1}$. The cases $g=1$ and $g=2$ are clear, since we have isomorphisms
 $Spin_3\cong Sp_2$ and $Spin_5\cong Sp_4$. We will prove the identity
 $vol_{Ch}(Sp_{2g})\;= 2\cdot vol_{Ch}(SO_{2g+1})$ by induction on $g$ and first give an outline of the argument:
 \absatz
 We will construct standard basis elements for $Lie(SO_{2g+1,K})$ and for $Lie(Sp_{2g,K})$ and denote their
 wedge product by $\lambda_{St}^{G}\in \Lambda^{\dim G}\gGg$.
  The invariant differential form $\omega_{St}^G$ of top degree which is dual to $\lambda_{St}^{G}$ can then be used to define
 the standard volume $vol_{St}(G_K)=|\int_{G_K}\omega_{St}^G|$.
 We will prove
 \begin{align}
 \label{SOcomprel}\omega_{St}^{SO_{2g+1}}\quad &=\quad \pm 2^{-(g+1)} \cdot \omega_{Ch}^{SO_{2g+1}}\qquad\text{ and}\\
 \label{Spcomprel} \omega_{St}^{Sp_{2g}}\quad &=\quad \pm 2^{-g^2} \cdot \omega_{Ch}^{Sp_{2g}}.\qquad\text{ }
 \end{align}
 The action of $SO_n$ on the unit sphere $S^{n-1}$ with stabilizer $SO_{n-1}$ will imply the
 relation
 \begin{align} \label{SOrek} vol_{St}(SO_{n})\quad&=\quad vol(S^{n-1})\cdot vol_{St}(SO_{n-1})
 \end{align} and the
 action of $Sp_{2g,K}$ on the unit sphere $S^{4g-1}$ in $\C^{2g}$ will imply the relation
 \begin{align} \label{Sprek} vol_{St}(Sp_{2g}) \quad =\quad vol(S^{4g-1})\cdot vol_{St}(Sp_{2(g-1)}),
 \end{align}
 where we use the Euklidean volume of the spheres, which is known to be
 $$ vol(S^{n-1})\quad =\quad \frac{2\cdot \pi^{n/2}}{\Gamma(n/2)}.$$
 Thus we get the recursive relations:
  $$ vol_{Ch}(Sp_{2g})=2^{g^2}\cdot vol_{St}(Sp_{2g})=
  2^{2g-1}\cdot 2^{(g-1)^2}\cdot vol(S^{4g-1})\cdot vol_{St}(Sp_{2(g-1)})$$
  $$= \frac{2^{2g}\cdot \pi^{2g}}{\Gamma(2g)}\cdot vol_{Ch}(Sp_{2(g-1)})$$
  and
  $$ vol_{Ch}(SO_{2g+1})=2^{-(g+1)}\cdot vol_{St}(SO_{2g+1})=
     \frac{2^{-g}}{2}\cdot vol(S^{2g})\cdot vol(S^{2g-1})\cdot vol_{St}(SO_{2g-1})$$
     $$=\quad \frac{2 \pi^{g+1/2}\cdot 2 \pi^{g}}{2\cdot \Gamma(g+1/2)\cdot \Gamma(g)}\cdot
     vol_{Ch}(SO_{2(g-1)+1}).$$
  By the doubling formula $\Gamma(g+1/2)\cdot \Gamma(g)=\frac{\sqrt\pi}{2^{2g-1}}\cdot \Gamma(2g)$
  the latter simplifies to
  $$ vol_{Ch}(SO_{2g+1})\quad =\quad \frac{2^{2g}\cdot \pi^{2g}}{ \Gamma(2g)}\cdot
     vol_{Ch}(SO_{2(g-1)+1})$$ and now the induction step from $g-1$ to $g$ is immediate.

  \Absatz\Numerierung \label{prepcalc} {\sc Preparation of explicit calculations:}
  For $N\in\N$ we denote by $E_{ij}$ the standard basis elements of the space of $n\times n$ matrices.
  Let $w_N=(\delta_{i,N+1-j})\in GL_N$ be the standard antidiagonal matrix. For elements $\epsilon_j=\pm 1$ we define
  $J=(\epsilon_j\delta_{i,N+1-j})$. Then $J^{-1}={}^tJ=(\epsilon_i\delta_{i,N+1-j})$.
  The algebraic group
  $G=\{A\in GL_N| \; {}^tA\cdot J\cdot A=J\}$ has the Lie algebra $\gGg=\{A\in Mat_{N\times N}| J^{-1}\cdot{}^tA\cdot J=-A\},$
  i.e. a matrix $A=(a_{ij})$ belongs to $\gGg$ iff we have
  $$ \epsilon_i\epsilon_j\cdot a_{N+1-j,N+1-i}=-a_{ij}\quad\text{ for all } 1\le i,j\le N.$$
  If we define
  \begin{align*} u_{ij}\quad &=\quad E_{ij}-\epsilon_i\epsilon_jE_{N+1-j,N+1-i}\in \gGg\quad\text{ for } i+j\ne N+1\\
    u_i\quad&=\qquad\qquad E_{i,N+1-i}\qquad\qquad\qquad\text{ for } 1\le i \le N
    \end{align*}
    then the elements $u_{ij}$ for $i+j<N+1$ and those $u_i$ for which $\epsilon_i\epsilon_{N+1-i}=-1$ form a basis
    of $\gGg$.
   Introducing the notation $t_i=u_{ii}$ (thus $t_{N+1-i}=-t_i$ and  $t_i=0$ for $2i=N+1$)
   and $t_{ij}= t_i-t_j$  for
   $i+j\ne N+1$ we get the commutation relations:
   \begin{align}\label{uu vertauschung}
   [u_{ij},u_{ji}]\quad&=\quad t_{ij}\\
      [u_i,u_{N+1-i}]\quad&=\quad t_i
      \end{align}
      For $i+j\ne N+1,\; i\ne j$ we  have furthermore
      \begin{align}\label{tu vertauschung}
      [t_{ij},u_{ij}]\quad &= \quad \begin{cases} 2u_{ij}\qquad\text{ if } 2j\ne N+1\ne 2i\\
                                    \; u_{ij}\qquad\text{ if $2j=N+1$ or $2i=N+1$.}
                                    \end{cases}
      \end{align}
      For $2i\ne N+1$ we finally get
      \begin{align*}
      [t_i,u_i]\quad &=\quad 2\cdot u_i .
      \end{align*}

   \Absatz\Numerierung\label{Sp2g} {\sc The case $Sp_{2g}$.}
   Putting $N=2g$ and $\epsilon_j=-1$ for $j\le g$ and $\epsilon_j=1$ for $g+1\le j\le 2g$ we have ${}^tJ=-J$ and get the group
   $G=Sp_{2g}$. A basis of the Lie algebra $\gGg$ consists of the elements $u_{ij}$ for $i+j<2g+1$ and $u_i$ for $1\le i\le 2g$.
   If we take $\hhh$ to be the subalgebra of diagonal matrices and $\bbb$ to be the subalgebra of upper triangular matrices we
    can use the following basis to compute the  Chevalley element (taking into account the commutation relations and $2i\ne 2g+1\ne 2j$ for all $i,j$):
   $$ u_{ij} \text{ for } i\ne j, i+j<2g+1;\quad u_i\text{ for } 1\le i\le 2g;\quad t_g;\quad t_{i,i+1}\text{ for } 1\le i\le g-1.$$
   From $t_{1,2}\wedge\ldots t_{g-1,g}\wedge t_g= u_{11}\wedge \ldots\wedge u_{g-1,g-1}\wedge u_{gg}$ we thus get the Chevalley element
   $$ \lambda_{Ch}\quad=\quad \pm \bigwedge_{i+j<2g+1} u_{ij}\; \wedge \; \bigwedge_{1\le i\le 2g} u_i.$$
   The inner forms of $G$ can be described as
   $$G_K(\R)\quad =\quad \left \{ A\in G(\C)|\; \overline A=BAB^{-1}\right\}, $$
   where $B\in G(\C)$ has to satisfy the cocycle relation that $\overline B\cdot B$ is trivial in the adjoint group, i.e.
   $\overline B\cdot B$ has to be a central matrix in $G(\C)$.
   In our case $G=Sp_{2g}$ we can take $B=J$  (observe $J^2=-id$) to get the compact inner form $G_K$.
    For symplectic matrices satisfying ${}^tA\cdot J\cdot A=J$  the condition
   $ \overline A=JAJ^{-1}$ is equivalent to the unitarity condition ${}^tA\cdot \overline A=E$ respectively
   ${}^t\overline A\cdot  A=E$. Thus
   \begin{align*}
   G_K(\C)\quad &=\quad \left\{ A\in \GL_{2g}(\C)\; | \; {}^tA\cdot J\cdot A=J, \qquad {}^t\overline A\cdot A=E\right\}\\
   Lie(G_K)\quad&=\quad \left\{ A\in Mat_{2g}(\C)\; | \; J^{-1}\cdot{}^t\!A\cdot J=-A, \qquad {}^t\overline A=-A\right\}.
   \end{align*}
   Then we get the following standard basis for the real Lie algebra $Lie(G_K)$:
   \begin{align*} i t_\nu\quad\qquad &\text{ for } 1\le \nu\le g\\
   u_{\nu\mu}-u_{\mu\nu}\;\quad &\text{ for } \nu<\mu,\; \nu+\mu<2g+1\\
   i(u_{\nu\mu}+u_{\mu\nu})\quad &\text{ for } \nu<\mu,\; \nu+\mu<2g+1\\
    u_\nu-u_{2g+1-\nu}  \;\quad &\text{ for } 1\le \nu\le g\\
    i( u_\nu+u_{2g+1-\nu})  \quad &\text{ for } 1\le \nu\le g
     \end{align*}
     From $(u_{\nu\mu}-u_{\mu\nu})\wedge i(u_{\nu\mu}+u_{\mu\nu})= 2i\cdot  u_{\nu\mu}\wedge u_{\mu\nu}$
     and the corresponding formula for the
     $u_\nu$ we immediately get (observing that we have $g^2$ positive roots):
     $$ \lambda_{St}= \pm (2i)^{g^2}\cdot i^g \lambda_{Ch}=\pm (-1)^{g(g+1)/2}\cdot 2^{g^2}\lambda_{Ch}.$$ The claim
     (\ref{Spcomprel})  is just the statement for the dual space.
     \klabsatz
     As subgroup of the standard unitary group the group $G_K(\R)=\Sp_{2g,K}(\R)$ acts on the unit sphere
     $S^{4g-1}=\{(a_\mu)\in \C^{2g}|\quad \sum_{\mu} |a_\mu|^2=1\}\subset \C^{2g}=\R^{4g}$ transitively in such a way that the
     stabilizer of $e_1=(1,0,\ldots,0)$ is isomorphic to $Sp_{2(g-1),K}(\R)$. By this action the $4g-1$ standard basis vectors of
     $Lie(Sp_{2g,K})$ with $\nu=1$ map to the euclidean base of the tangent space of $S^{4g-1}$ at $e_1$, while the
     other standard basis vectors descend to the standard basis of $Lie(Sp_{2(g-1),K})$. From these facts the recurrence formula
     (\ref{Sprek}) may be deduced easily.

    \Absatz\Numerierung {\sc The case $SO_{2g+1}$.} Now we put $N=2g+1$ and $\epsilon_j=1$ for $1\le j\le 2g+1$ in
    \ref{prepcalc} and get the split orthogonal group $G=SO_{2g+1}=\{A\in \GL_{2g+1}(\R)| {}^t\!A\cdot w_N\cdot A=w_N\}$
    with Lie algebra $\gGg=\mathfrak{so}_{2g+1}$. The elements $u_{ij}$ for $i+j<2g+2$ form a basis of $\gGg$.
    The commutation relations (\ref{uu vertauschung}) and (\ref{tu vertauschung}) now imply,
    that we can get the Chevalley element as the wedge product of the following basis elements:
    \begin{align*} u_{ij}\qquad\quad&\text{ for } i\ne j,\quad i+j<2g+2,\quad j\ne g+1\ne i\\
                   u_{i,g+1} \text{ and } 2\cdot u_{g+1,i}\quad&\text{ for } 1\le i\le g\\
                   t_{i,i+1}\qquad \; &\text{ for } 1\le i\le g-1\quad\text{ and finally}\\
                   2 t_g.\qquad\quad &{}\end{align*}
      From this we deduce immediately:
      $$\lambda_{Ch}\quad=\quad  \pm 2^{g+1}\cdot \bigwedge_{i+j<2g+2} u_{ij}.$$

      Let $SO_N^K\quad=\quad \{A\in \SL_N(\R)| {}^tA\cdot A= id\}$ be the compact form of $G$.
      Then the elements $s_{ij}=E_{ij}-E_{ji}$ for $1\le i< j\le N$ form a basis of its Lie algebra
      $\mathfrak{so}^K_N=\{A\in Mat_N(\R)| {}^t\! A=-A\}$.
      To get $SO^K_{2g+1}$ as an inner twist of the split group $SO_{2g+1}$ we introduce the block matrices
      $$ C=\begin{pmatrix} id_g & 0 & \frac 12 w_g\\ 0 & 1& 0\\ i\cdot w_g & 0 &-\frac i2 id_g\end{pmatrix}
      \quad \text{ with inverse }
         C^{-1}=\begin{pmatrix} \frac 12 id_g & 0 &  -\frac i2 w_g\\ 0 & 1& 0\\  w_g & 0 & i\cdot id_g\end{pmatrix}.$$
         Then we have ${}^t\! C\cdot C= J$, and the matrix
         $$ B\quad =\quad \overline C^{-1}\cdot C\quad=\quad
         \begin{pmatrix}0&0& \frac 12 w_g\\ 0& 1& 0\\ 2 w_g & 0 & 0\end{pmatrix}$$
         satisfies $\overline B\cdot B= id_{2g+1}$ and can thus be used to define the inner form
         $\widetilde{SO}_{2g+1,K}= \{ A\in SO_{2g+1}(\C)| \overline A =B\cdot A\cdot B^{-1}\}.$
         (One may replace $B$ by $-B$ if $\det(B)=-1$ to get an element in $SO_{2g+1}(\C)$.)
         With the notation $\phi(A)=C\cdot A\cdot C^{-1}$ we now have the equivalences:
         $$ A\in SO_{2g+1}(\C)\Leftrightarrow {}^t\!A\cdot ({}^t\! C\cdot C)\cdot A={}^t\! C\cdot C
          \Leftrightarrow {}^t\!(CAC^{-1})\cdot (CAC^{-1})= id_{2g+1}$$
          and
         $$ \overline A =B\cdot A\cdot B^{-1} \Leftrightarrow
         \overline A =\overline C^{-1}\cdot C\cdot A\cdot (\overline C^{-1}\cdot C)^{-1}
         \Leftrightarrow \overline{CAC^{-1}} = CAC^{-1} \Leftrightarrow \phi(A)=\overline{\phi(A)}.$$
         Thus we get an isomorphism
         $$ \phi: \widetilde{SO}_{2g+1,K}\quad \tilde  \to \quad SO_{2g+1}^K,\qquad A\mapsto CAC^{-1}.$$
         On the Lie algebra this isomorphism is given by the same formula.
         Its effect on the basis elements $u_{ij}$ is given by the following rules where $1\le \mu,\nu\le g$:
         \begin{align*}
          \phi(u_{\mu\nu})\; &=\; \frac 12\cdot\left( s_{\mu\nu}-s_{2g+2-\nu,2g+2-\mu}-i\cdot s_{\mu,2g+2-\nu}-i\cdot s_{\nu,2g+2-\mu}\right)\\
          \phi(u_{\mu,2g+2-\nu})\; &=\qquad  s_{\mu\nu}+s_{2g+2-\nu,2g+2-\mu}+i\cdot s_{\mu,2g+2-\nu}-i\cdot s_{\nu,2g+2-\mu}\qquad
          \text{ for } \mu<\nu\\
          \phi(u_{2g+2-\nu,\mu})\; &=\; \frac 14\cdot
          \left(- s_{\mu\nu}-s_{2g+2-\nu,2g+2-\mu}+i\cdot s_{\mu,2g+2-\nu}-i\cdot s_{\nu,2g+2-\mu}\right)
          \quad \text{for } \mu<\nu\\
          \phi(u_{\mu, g+1})\; &=\qquad s_{\mu,g+1}-i\cdot s_{g+1,2g+2-\mu}\\
          \phi(u_{g+1,\mu})\; &=\; \frac 12\cdot\left( - s_{\mu,g+1}-i\cdot s_{g+1,2g+2-\mu}\right).
         \end{align*}
         Here we put $s_{\mu\mu}=0$. Especially we have for $1\le \mu=\nu\le g$:
         $$\phi(t_\mu)=\phi(u_{\mu\mu})=(-i)\cdot s_{\mu,2g+2-\mu}.$$
         In the case $1\le\mu<\nu \le g$ we get by applying the first equation to $u_{\mu\nu}$ and to $u_{\nu\mu}$, using
         $s_{\nu\mu}=-s_{\mu\nu}$ and combining this with the next two equations:
         $$\phi\left( u_{\mu\nu}\wedge u_{\nu\mu}\wedge u_{\mu,2g+2-\nu}\wedge u_{2g+2-\nu,\mu}\right)
         \;=\; s_{\mu\nu}\wedge s_{2g+2-\nu,2g+2-\mu}\wedge s_{\mu,2g+2-\nu}\wedge  s_{\nu,2g+2-\mu}.$$
         Similarly we get from the last two equations:
         $$ \phi\left( u_{\mu,g+1}\wedge u_{g+1,\mu}\right)\quad=\quad (-i) \cdot s_{\mu,g+1}\wedge s_{g+1,2g+2-\mu}.$$
         Forming the wedge product of the last three families of equations we obtain:
         $$\phi\left(\bigwedge_{\mu+\nu<2g+2} u_{\mu\nu}\right)\; =\; \pm (-i)^g\cdot \bigwedge_{1\le i<j\le 2g+1} s_{ij}
         \qquad\text{ and finally }$$
         $$ \phi(\lambda_{Ch})\;=\;\pm 2^{g+1}\cdot \bigwedge_{1\le i<j\le 2g+1} s_{ij}.$$
          If we take the $s_{ij}$ as standard basis for the compact form $SO^K_{2g+1}$, we get the claim (\ref{SOcomprel}).
          One gets the recurrence relation (\ref{SOrek}) as in the case of the symplectic group. Now the proposition is proved.
         \qed

\meinchapter{The constants $\alpha$} \label{Ch alpha}

   \Absatz\Numerierung
   Let $G=\PGL_{2g+1}/\Q$ and $G_1=\SL_{2g}/\Q$ be its stable $\eta$-endoscopic group, where $\eta: g\mapsto
   J_{2g+1}^{-1}\cdot \tran g^{-1} J_{2g+1}.$ Let $\Delta$ be the set of simple roots of $G$ with respect to the Borel $P_\emptyset$
   of upper triangular matrices. Recall that the parabolic subgroups containing $P_\emptyset$ are parametrized by
   the subsets $I\subset \Delta$ and are denoted by $P_I$. We may identify the space of $\eta$-orbits $\Delta/\eta$ with the
   set of simple roots of $G_1$ in such a way that an $\eta$-stable $I\subset \Delta$ corresponds to a parabolic subgroup
   $P_{I,1}\subset G_1$ containing the Borel $P_{\emptyset,1}$ of upper triangular matrices.
   Since $G$ and $G_1$ are semisimple, we can arrange, that the open compact subgroups  $Z_f\subset Z_G(\A_f)$
   are small enough in the sense that the group
    $\zeta = Z_G(\QQ)\;\cap\;\left(K_\infty \cdot Z_\infty\;\times\; Z_f\right)$ is trivial in both cases.

    \Absatz
   The aim of this chapter is to prove:

   \begin{prop} \label{alpha} Let $I\subset \Delta$ be $\eta$-stable. If $\gamma\in P_I(\Q)$ and $\gamma_1\in P_{I,1}(\Q)$ are
   $\eta$-matching and contribute to the topological trace formula and if we use the Chevalley measures on the
    centralizers inside $G_1$ and $2$ times the Chevalley measures on the twisted centralizers inside $G$, then we have
   $$ \alpha_\infty(\gamma,1)\quad=\quad \alpha_\infty(\gamma_1,1).$$
   \end{prop}
  \Absatz\Numerierung Recall from \cite[Thm 4.9]{Uwe} the definition
    \begin{align*}\alpha_\infty(\gamma_0,h_\infty)\quad=
         \quad \frac{O^\infty_\eta\left(I,\gamma_0,h_\infty\right)}{d^{I}_{\zeta,\gamma_0}}
          \cdot(-1)^{\Delta(\gamma_0,\eta)}&\cdot
          \frac{\# H^1(\RR,T)}{vol_{db_\infty}\left((\overline{G_{\gamma_0,\eta}^{I}})'/\zeta\right)}.
   \end{align*}
   The factors on the right hand side will be explained in the following.  Thereby we will prove their coincidence
   for  $\gamma$ and $\gamma_1$:

   \Absatz\Numerierung For the definition of the sign factor $(-1)^{\Delta(\gamma_0,\eta)}$ we refer to \cite[Thm. 4.9]{Uwe}
   and \cite[3.9]{Uwe}. The coincidence of the sign factors has already be proved in \cite[5.26]{Uwe}.

   \Absatz The numbers $d^{I}_{\zeta,\gamma_0}$ are the cardinalities of certain non empty subsets
   of $H^1(\langle\eta\rangle,\zeta)$ as defined in \cite[2.24]{Uwe}. But since we assume that $\zeta$ is the trivial group,
   they all have to be $1$.

   \Absatz\Numerierung
   Recall that  $    O^\infty_\eta\left(I,\gamma,h_\infty\right)$
  is the order of the coset space
    \begin{align*}
       R_{\gamma,\eta}^{I}\quad&=\quad L_{\gamma,\eta}^m\bs \left(\tilde L^m/\tilde L\right)^{\eta_\gamma}
        \end{align*}
        where \begin{align*}
           \tilde L\quad&=\quad p_1\cdot K_\infty^{I} Z_\infty A_I \cdot p_1^{-1},\qquad\qquad
            \tilde L^m\quad=\quad p_1\cdot K_\infty^{I,m} Z_\infty A_I\cdot p_1^{-1},\\
                      L_{\gamma,\eta}^m\quad&=\quad\tilde L^m\cap G_{\gamma,\eta}^{I}(\RR),\qquad\qquad\qquad
            \eta_\gamma (x) \quad=\quad  (g_\eta\gamma)^{-1}\cdot \eta(x)\cdot g_\eta \gamma,\\
            K_\infty^{I}\quad&=\quad K_\infty\cap P_I(\R),\qquad\qquad\qquad
            K_\infty^{I,m}\quad=\quad K_\infty^m\cap P_I(\R).\\
        \end{align*}
     But for $G$ the group $K_\infty=SO_{2g+1}\subset G(\R)= \PGL_{2g+1}(\R)$ is already maximal compact, as is the
     subgroup $K_\infty=U_{g}\subset G_1(\R)=\Sp_{2g}(\R)$. Thus we have $K_\infty=K_\infty^m$ in both cases, which implies
     $\tilde L= \tilde L^m$ and therefore
     $O^\infty_\eta\left(I,\gamma,1\right)=1= O^\infty_\eta\left(I,\gamma_1,1\right)$.

      \Absatz\Numerierung
    To handle the two remaining factors we have to calculate the (twisted) centralizers:
      \begin{align*}
           G_{\gamma,\eta}^{I}(F)\quad&=\quad \set{A\in P_I(F)|\eta_\gamma(A)=A}\quad=\quad
           \set{A\in P_I(F)| \eta(A)^{-1}\cdot \gamma\cdot A= \gamma}.
       \end{align*}
       Here $F$ denotes a field of characteristic $0$. For $I=\Delta$, i.e. $P_\Delta=G$,  we omit the index: $G_{\gamma,\eta}=G_{\gamma,\eta}^\Delta$.
    Since $\eta:  A\mapsto J^{-1}\cdot {}^t\! A^{-1}\cdot J$ is formed with
    a symmetric and involutive antidiagonal matrix $J=J_{2g+1}$, the twisted centralizers may be written in the form:
    \begin{align*}
       G_{\gamma,\eta}(F)\;&=\;    \set{A\in G(F)|\; {}^t\! A\cdot J\gamma\cdot A\;=\; J\gamma}.
    \end{align*}

    \begin{lemma} (a) If $\chi_{I,\alpha}\left(\NNN(\gamma)\right)>1$ for all $\alpha\in\Delta -I$ and if
    $\NNN(\gamma)$ is conjugate to an element of $L_\infty^{I}=K_\infty^{I}Z_\infty A_I$, then we have
      $$G_{\gamma,\eta}=G_{\gamma,\eta}^{I}.$$
        \end{lemma}
    Proof:  Since the twisted centralizer $G_{\gamma,\eta}$ is contained in the usual centralizer
    $G_{\NNN(\gamma)}$ of the norm, it suffices to prove that $G_{\NNN(\gamma)}\subset P_I$. But this is a consequence of the two
    conditions: Since
     $\NNN(\gamma)$ is conjugate to an element of $L_\infty^{I}$, it is semisimple and may be conjugated inside $P_I(\C)$ to
      a diagonal matrix $\delta'$, where  in each block of the Levi $M_I$ the diagonal entries have the same absolute value.
      Then the condition $\chi_{I,\alpha}\left(\NNN(\gamma)\right)>1$ for all $\alpha\in\Delta -I$ implies, that the absolute values
      of these entries are different for different blocks of $M_I$. Thus the centralizer of $\delta'$ is contained in $M_I$, and
      the centralizer of $\NNN(\gamma)$ must be contained in $P_I$.\qed
        \klabsatz
     Since the two conditions are satisfied for $\gamma$ contributing to the topological trace formula \cite[Theorem 4.9.]{Uwe}, we
     will assume this in the sequel in proving proposition \ref{alpha}. Thus we may compute the twisted centralizers in the full group $G=P_\Delta$ resp. $G_1=P_{1,\Delta}$ and forget the index $I$.
     \begin{lemma}
     The canonical projection $\SL_{2g+1}\to \PGL_{2g+1}$ induces an isomorphism:
    $$ \set{A\in \SL_{2g+1}(F)|\; {}^t\! A\cdot J\gamma\cdot A\;=\; J\gamma}\quad \tilde\to\quad G_{\gamma,\eta}(F).$$
        \end{lemma}
    Proof:  If $A\in \GL_{2g+1}(F)$ represents an element of $G_{\gamma,\eta}(F)$, then
    ${}^t\! A\cdot J\gamma\cdot A\;=\; \lambda\cdot J\gamma$ for some $\lambda\in F^*$.
     Since $2$ and $2n+1$ are coprime there exists a unique $\mu\in F^*$, such that the two equations
    $\det(\mu A)= \mu^{2n+1}\det(A)=1$ and $\mu^2\lambda=1$ are satisfied i.e. such that $\mu A\in \SL_{2n+1}(F)$
    and ${}^t\! (\mu A)\cdot J\gamma\cdot (\mu A)\;=\;  J\gamma$. \qed
   \Absatz\Numerierung
    Recall  that $\NNN(\gamma)$ has a conjugate in $L_\infty^{I}$.
    But then $\NNN(\gamma)= \eta(\gamma)\cdot\gamma=  J^{-1}\cdot {}^t\gamma^{-1}\cdot J \cdot\gamma ={}^t(J\gamma)^{-1}\cdot J\gamma$ is semisimple, since every element in $L_\infty^{I}$ is semisimple.  Note that $\NNN(\gamma)$ is well defined as an element of $\SL_{2g+1}(F)\subset \GL_{2g+1}(F)$, since
    we have $\eta(z\gamma)\cdot z\gamma=\eta(\gamma)\gamma$ for $z$ in the center of $\GL_{2g+1}$ and since
    $\det(\eta(\gamma))=\det(\gamma)^{-1}$.
    For a matrix $M\in \GL_n(F)$ put $Aut(M)=\set{ A\in GL_n(F)|\; {}^t\!A\cdot M\cdot A=M}$. If $B\in GL_n(F)$ represents
    $J\gamma \in \PGL_n(F)$, we thus get
    $G_{\gamma,\eta}(F)\cong Aut(B)\cap \SL_{2g+1}(F)$ with semisimple ${}^t\!B^{-1}\cdot B$.
    We note that $A\in Aut(M)$ implies $A\in Aut({}^tM)$ by transposing the equation, and thus
    $A\in Aut(M+{}^tM)$ and $A\in Aut(M-{}^tM).$

    \begin{lemma} \label{Autlemma}
    Let $B\in \GL_n(F)$ with semisimple $S={}^t B^{-1}\cdot B$.
    Assume that $S$ has block diagonal form $S=diag(S_0, 1_m)$ such that $S_0-1_{n-m}\in\! \GL_{n-m}(F)$.
     \begin{itemize}\item[(a)]    The matrix $B$ is block diagonal:    $B=diag(B_0, B_s)$ where $B_s\in \GL_m(F)$ is symmetric, where $S_0={}^t B_0^{-1}\cdot  B_0$ and where
           the skew symmetric part $P_0=B_0-{}^tB_0$ lies in $\GL_{n-m}(F)$.
    \item[(b)]  We have $Aut(B)=Aut(B_0)\times Aut(B_s)$.
    \item[(c)] We have $S \in Aut(B-{}^tB)$ and $S_0\in Aut(P_0)$.
    \item[(d)] We have $Aut(B_0)=\{A\in Aut(P_0)| A\cdot S_0=S_0\cdot A\}$.
    \end{itemize}
    \end{lemma}
Proof:
      (a) If one writes $B$ as a block matrix, then the equation ${}^t BS=\;B$ easily implies, that $B$ has to be
      block diagonal, since $S_0-1_{n-m}$ is invertible, that $B_s$ is symmetric and $S_0={}^t B_0^{-1}\cdot  B_0$.
      But then $P_0=B_0\cdot (1_{n-m} -S_0^{-1})$ has to be invertible too.
      \klabsatz
      (b) The claim is a consequence of the fact that each  $A\in Aut(B)$ has to be block diagonal:
      From $Aut(B)\subset Aut(B-{}^tB)=Aut(diag(P_0,0))$ we deduce that $A$ has to be block lower triangular, since $P_0$ is invertible, and from
      $Aut(B)\subset Aut(B+{}^tB)$ we deduce that $A$ is block diagonal, since $B_s$ is invertible.
      \klabsatz
      (c) ${}^tS\cdot (B-{}^tB)\cdot S=
      {}^tB\cdot B^{-1}\cdot B\cdot {}^tB^{-1}\cdot B-{}^tB\cdot B^{-1}\cdot{}^t B\cdot {}^tB^{-1}\cdot B = B-{}^tB,$
      and the statement for $S_0$ is a consequence of this identity.
      \klabsatz
      (d) If $A\in Aut(B_0)$ then $A\in Aut(P_0)$ and $${}^tB_0^{-1}\cdot B_0={}^t({}^t\! AB_0A)^{-1}\cdot ({}^t\! AB_0A)
      =A^{-1}\cdot {}^tB_0^{-1}\cdot{}^t\! A^{-1}\cdot {}^t\!A\cdot B_0\cdot A= A^{-1}\cdot {}^tB_0^{-1}\cdot B_0\cdot  A$$
      implies $AS_0=S_0A$.
      Conversely $A\in Aut(P_0)$ means $${}^tA\cdot B_0(1_{n-m}-S_0^{-1})\cdot A = B_0(1_{n-m}-S_0^{-1}).$$
      If furthermore $A\cdot S_0=S_0\cdot A$ then
      we get ${}^t\! A\cdot B_0\cdot A\cdot (1_{n-m}-S_0^{-1})= B_0(1_{n-m}-S_0^{-1})$, and then the invertibility of
      $1_{n-m}-S_0^{-1}=S_0^{-1}\cdot(S_0-1_{n-m})$ implies
       $A\in Aut(B_0)$.
\qed
    \Absatz\Numerierung
     For $\beta\in \GL_{2g+1}(F)$ we denote its image in $\PGL_{2g+1}(F)$ by $\overline\beta$ and
      put $\tilde \gamma= \eta(\overline\beta)^{-1}\cdot\gamma\cdot \overline\beta\in \PGL_{2g+1}(F)$.
            Then we have:
      \begin{align*}
       G_{\tilde \gamma,\eta}\quad &=\quad \beta^{-1}\cdot G_{\gamma,\eta}\cdot \beta\\
       \NNN(\tilde\gamma)\quad&=\quad \beta^{-1}\cdot \NNN(\gamma) \cdot \beta \qquad\text{ (\cite[2.3 (14)]{Uwe})}\\
       \tilde B \quad :&= \quad {}^t\beta\cdot B\cdot \beta\quad\text{ represents } J\tilde\gamma,\\
       \tilde S \quad &=\quad \beta^{-1}\cdot S\cdot \beta \quad \text {where } \tilde S := \tilde {}^t B^{-1}\cdot  \tilde B.
      \end{align*}
       By varying  $\gamma$ in its $\eta$-conjugacy class  we may thus assume, that the semisimple matrix $S$ has block diagonal form
       $S=diag(S_0,1_m)$ with invertible $S_0-1_{2g+1-m}$, that $B=diag(B_0,B_s)$ is block diagonal with symmetric $B_s$ and
        that the antisymmetric part $P_0=B_0-{}^tB_0$ is the standard antidiagonal symplectic matrix $J_{2p}\in GL_{2p}(F)$ as in
        \ref{Sp2g}.
        Then we have $m=2(g-p)+1$.
        Now we put $\delta=diag(1_{g-p}, S_0,1_{g-p})\in \Sp_{2g}=: G_1$.

        \begin{lemma} \label{Zentralisatoren} With the notations introduced above we have:
        \begin{itemize}
        \item[(a)] $G_{\gamma,\eta}\quad\cong\quad SO(B_s)\times \{A\in \Sp_{2p}|\; AS_0=S_0A\}$.
        \item[(b)] $(G_1)_\delta\quad\cong\quad Sp_{2(g-p)}\times \{A\in \Sp_{2p}|\; AS_0=S_0A\}$.
        \item[(c)] The elements $\gamma\eta$ and $\delta$ are matching.
        \end{itemize}
        \end{lemma}
        Proof: (a) Since $Aut(B_0)\subset\Sp_{2p}\subset \SL_{2p}$ lemma \ref{Autlemma}(b) implies
        $$G_{\gamma,\eta}\quad = \quad Aut(B)\cap SL_{2g+1}\cong (Aut(B_s)\cap \SL_{2(g-p)+1})\times
                 Aut(B_0)$$ and the claim is a consequence of  lemma \ref{Autlemma}(d).
         \klabsatz
         (b) This is clear, since $S_0$ does not have the eigenvalue $1$.
         \klabsatz
         (c) Since this is essentially \cite[prop 6.3.(c)]{FLPGL5} we sketch the argument:
         Since the claim may be checked over an algebraic closed field $\bar F$,
         one can assume $\gamma$ to be a diagonal matrix. Then $B$ is antidiagonal and $S$ is diagonal, and an explicit calculation
          using the definitions of stable endoscopy implies the claim. \qed

\Absatz\Numerierung Now we can finish the proof of proposition \ref{alpha}: Since the elements in $P_{I,1}(\Q)$ matching with
       $\eta\gamma$ are all stable conjugate we may assume $\gamma_1=\delta$ in the notations of the preceding lemma.

       \absatz Recall that in the factor $\# H^1(\R,T)$ the torus $T$ can be chosen to be a maximal torus in the centralizer
       $G_{\gamma_0,\eta}$. Now we may chose some maximal torus $T_0$ inside the group $C=\{A\in \Sp_{2p}|\; AS_0=S_0A\}$,
        an $\R$-anisotropic maximal torus $T_1'$ in $Sp_{2(g-p)}$  (of rank $g-p$) and an $\\R$-anisotropic maximal torus $T'$ in
        $SO(B_s)$ again of rank $g-p$.
        Then we can chose the maximal torus $T$ in $G_{\gamma,\eta}$ in such a way that it maps to $T'\times T_0$ under the
        isomorphism of \ref{Zentralisatoren}(a) and the maximal torus $T_1$ in $(G_1)\delta$  in such a way that it maps
        to $T_1'\times T_0$ under \ref{Zentralisatoren}(b). Since we have $\# H^1(\R,T')=2^{g-p}=\# H^1(\R,T_1')$ we get immediately
        $\# H^1(\R,T)=\# H^1(\R,T_1)$.

  \Absatz\Numerierung For the remaining factor $vol_{db_\infty}\left((\overline{G_{\gamma_0,\eta}^{I}})'/\zeta\right)$
  we first recall the fact that $\zeta$ is the trivial group in both cases. For $G=G_{\gamma_0,\eta}^{I}=G_{\gamma_0,\eta}$ the group $\overline G/\R$ is
     the inner form of $G$ which is compact modulo the center of $G$ (\cite[3.9]{Uwe}), and then $\overline G'$ denotes
     the common kernel of all rational characters of $G$ (which may be viewed as characters of $\overline G$ using the
     isomorphism $G/G^{der}=\overline G/{\overline G}^{der}$).  But lemma \ref{Zentralisatoren} implies that we have
     $\overline{G_{\gamma,\eta}}\,'\;\cong SO_{2(g-p)+1} \times \overline C'$ and
     $\overline{(G_1)_{\delta}}\,'\;\cong Sp_{2(g-p),K}(\R) \times \overline C'$, where $Sp_{2(g-p),K}$ is the compact inner form
     of $Sp_{2(g-p)}/\R$.
     From prop \ref{BC Chevalleyvolume} and our assumptions on the measures we now get
     $$vol_{db_\infty}\left((\overline{G_{\gamma,\eta}^{I}})'/\zeta\right)=
     vol_{db_{\infty,1}}\left((\overline{(G_1)_{\gamma_1}^{I}})'/\zeta\right).$$
     This finishes the proof of proposition \ref{alpha}.
     \qed

\meinchapter{Restatement of the Topological Trace formula}\label{Ch Restatement}

\Absatz\Numerierung\label{Untersuchte Situationen}
    Now we assume that we are in one of the following situations:
    \begin{align*} (G,\eta,G_1)\quad&=\quad (\PGL_{2n+1},\eta,\Sp_{2n})\\
       (G,\eta,G_1)\quad&=\quad (\GL_{2n}\times\GL_1,\eta,\GSPIN_{2n+1})\text{ where $n=1$ or $n=2$. }\\
           \end{align*}
     We regret that the condition $n=1$ or $n=2$ is missing in \cite[Theorem 5.23]{Uwe} due to an formulation error:
     In fact for $n\ge 3$ the group $\GSPIN_{2n+1}$ is not cohomological trivial (so the stabilization of the trace formula will be more
      complicated), and the groups $K_\infty$ have not been specified, unless we can use the exceptional isomorphisms
      $\GSPIN_5\cong \GSp_4$ and $\GSPIN_3\cong\GL_2$. Now we may restate \cite[Theorem 5.23]{Uwe} and \cite[Corollary 5.27]{Uwe}:

\begin{thm}\label{Cohomology lifting theorem}
   (a)  Assume that $(G,\eta,G_1)$ is as in \ref{Untersuchte Situationen}.
   If the $ G$-module $M$  matches with
    the $ {G_1}$-module $M_1$, then
    \begin{align*}
       H^*_c\left(G(\QQ)\bs  G(\A)/K_\infty Z_\infty,\MMM\right)
        \in \Groth( G(\A_f)\rtimes \eta)\\
        \intertext{\qquad is the lift of }
        H^*_c\left( G_1(\QQ)\bs  G_1(\A)/K_{\infty,1} Z_\infty,\MMM_1\right)
        \in \Groth( G_1(\A_f));\\
        \end{align*}
     (b)   The same statement holds if one the cohomology with compact support $H^*_c$ is replaced by the usual cohomology
        $H^*$.
\end{thm}

      In the situation $(G,\eta,G_1)= (\PGL_{2n+1},\eta,\Sp_{2n})$ we remark that due to proposition \ref{alpha}
      the notion of lifting depends on a notion of matching test functions on the group of finite adeles which  is now defined
      intrinsically using the Haar measures on the centralizers which are defined using a Chevalley basis.

\meinchapter{Spectral decomposition of cohomology}\label{Ch Franke}
\Absatz\Numerierung {\sc Frankes Theorem.}\\
      The functor $V\mapsto H^*(G(\Q)\bs G(\A)/K_\infty Z_\infty,\VVV)$ induces a homomorphism
      between Grothendieck groups: $H^{*,G}:\Groth(G,alg)\to \Groth(G(\A_f))$,
      where $\Groth(G(\A_f))$ denotes the Grothendieck group of admissible representations of $G(\A_f)$
      in a complex vector space.

      For $I\subset \Delta$ the naive (not normalized) parabolic induction functor induces a homomorphism:
      \begin{align*}
          Ind_{P_I(\A_f)}^{G(\A_f)}:\quad \Groth(M_I(\A_f))\quad&\longrightarrow\quad \Groth(G(\A_f)).
      \end{align*}

      Let $\check{\aaa}_I= X^*(P_I)\otimes\RR\cong X^*(A_I)\otimes\RR$ and $\pi_I:X^*(A_0)\otimes\RR \to \check{\aaa}_I$
      be the canonical projection, which annihilates the simple roots in $I$. Let $\rho_I=\pi_I(\rho)$, where $\rho=\rho_G\in X^*(A_0)\otimes\RR$ is half the sum of the positive roots.
      Let $\overline{\check{\aaa}_I^+}$ be the set of those $\lambda\in X^*(P_I)\otimes\RR=\check{\aaa}_I$ for which
      $\langle\lambda,\hat\alpha\rangle\ge 0$ for all simple roots $\alpha\in \Delta-I$. Here $\hat \alpha=\alpha^\vee$ denotes the corresponding simple coroot.
      \absatz

      After tensoring the Grothendieck groups with the field $\Q$ we can define the following
      homomorphisms for $I\subset \Delta$:
      \begin{align*}
          J_{Eis,I}:\Groth(G,alg)_\Q\quad &\to \quad \Groth(M_I,alg)_\Q \\
               V\quad&\mapsto\quad
               \sum_{\lambda\in \overline{\check{\aaa}_I^+}}\;\frac{1}{n_I(\lambda)}\; V^*_{I,\lambda}
      \end{align*}
      where $n_I(\lambda)$ is the number of Weyl chambers to which $\lambda$ belongs and where
      $V^*_{I,\lambda}$ is  the part of the virtual $M_I$ module  $H^*(\nnn_I,V)$, on which
      $\aaa_I=X_*(A_I)\otimes \RR$ acts by $-\lambda$.

      If $G$ is a split group, let $W$ be its Weyl group. In view of Kostants theorem  $J_{Eis,I}$ may be rewritten in the form
      \begin{align*}
          J_{Eis,I}:\Groth(G,alg)_\Q\quad &\to \quad \Groth(M_I,alg)_\Q \\
               V_\chi\quad&\mapsto\quad
               \sum_{\substack{w\in W^{I},\\
               \pi_I(w(\chi+\rho))\;\in\; -\overline{\check{\aaa}_I^+}}}\;
               \frac{(-1)^{l(w)}}{n_I(w(\chi+\rho)-\rho)}\;
               V^{I}_{w(\chi+\rho)-\rho}
      \end{align*}
      where $V_\chi$ is an irreducible algebraic $G$ module with highest weight $\chi$, where
      $V^{I}_{w(\chi+\rho)-\rho}$ denotes the corresponding highest weight module for $M_I$ and where
      $$W^{I}=\{w\in W| w^{-1}\alpha>0\text{ for }\alpha\in I\}=\{w\in W| w^{-1}\hat \alpha>0\text{ for }\alpha\in I\}.$$

      \absatz
      If $\chi: Z_\infty A_I\to \C^*$ denotes a continuous character, we define
      $\CCC^\infty(M_I(\Q)\bs M_I(\A),\chi)$  to be
      $$\set{f\in\CCC^\infty(M_I(\Q)\bs M_I(\A))\; |\; f(zg)=\chi(z)f(g)\;\text{for}\; z\in Z_\infty \cdot A_I},$$
      and denote by  $L_{2}(M_I(\Q)\bs M_I(\A),\chi)$  its completion as a Hilbert space with respect to integration
      over $M_I(\Q)\bs M_I(\A)'$, where $M_I(\A)'$ is the common kernel of all homomorphisms of the form $m\mapsto |\phi(m)|$, where
      $\phi:M_I\to \G_m$ is a rational character.
     By $L_{2,disc}(M_I(\Q)\bs M_I(\A),\chi)$ we denote the intersection of $\CCC^\infty(M_I(\Q)\bs M_I(\A),\chi)$ with the discrete spectrum in  $L_{2}(M_I(\Q)\bs M_I(\A),\chi)$.
      Now we can define a homomorphism
      $H^*_{2,disc}=H^{*,M_I}_{2,disc}:\Groth(M_I,alg)\to  \Groth(M_I(\A_f))$, if we put for an
      irreducible $M_I(\Q)$-modul $V$, on which $Z_\infty\cdot A_I$ acts by the central character $\chi_V$:
      \begin{align*}
          H^{*,M_I}_{2,disc}(V)\quad &=\quad H^*(\mmm_I,K_\infty\cap M_I(\RR); L_{2,disc}(M_I(\Q)\bs M_I(\A),\chi_V^{-1})\otimes V).
      \end{align*}

\begin{thm}\label{Frankes theorem}
      The following identity between homomorphisms from $\Groth(G,alg)$ to $\Groth(G(\A_f))_\Q$ holds:
      \begin{align*}
           H^{*,G}\quad&=\quad \sum_{I\subset \Delta} Ind_{P_I(\A_f)}^{G(\A_f)}\circ
              H^{*,M_I}_{2,disc}\circ J_{Eis,I}.
       \end{align*}
\end{thm}
   This is a consequence of  the spectral sequence \cite[(7.4.1)]{Franke}, since the $E_1$ term coincides  with the $E_\infty$ term
   in the Grothendieck group. The statement is closely related to the considerations in
   \cite[7.7.(1) and (2)]{Franke}.
   \qed

 \begin{lemma}\label{I von w}
  For $w\in W^{I}$ and $\pi_I(w(\chi+\rho))\in -\overline{\check{\aaa}_I^+}$ we have
  $I=\set{ \alpha\in \Delta |\, w^{-1}\alpha>0}$.
   \end{lemma}
   Proof: For $\alpha\in\Delta$ and $w\in W^{I}$ write $w^{-1}\hat\alpha=\sum_{\beta\in \Delta} k_\beta \hat\beta$.
   Then we have
   $$ \langle w(\chi+\rho),\hat\alpha\rangle\;=\; \langle \chi+\rho,w^{-1}\hat\alpha\rangle\;=\; \sum_{\beta\in\Delta}
   k_\beta \langle \chi+\rho,\hat\beta\rangle.$$
   Now $\langle \chi+\rho,\hat\beta\rangle=\langle \chi,\hat\beta\rangle+ 1\ge 1$ since $\chi$ is a dominant weight.
    In the case $w^{-1}\hat\alpha>0$ (especially for $\alpha\in I$) all the $k_\beta$ are non negative with at least one of them being positive. Consequently
    $ \langle w(\chi+\rho),\hat\alpha\rangle>0$. For $w^{-1}\hat\alpha<0$ we get $ \langle w(\chi+\rho),\hat\alpha\rangle<0$
   by the analogous argument.
   \absatz
   Now write $\pi_I(w(\chi+\rho))= w(\chi+\rho)- \sum_{\beta\in I} n_\beta \beta$.
   The condition $\langle \pi_I (w(\chi+\rho)),\hat\alpha\rangle=0$ for $\alpha\in I$ implies that
   $n_\beta=\sum_{\gamma\in I} c_{\beta\gamma} \langle w(\chi+\rho),\hat\gamma\rangle$, where
   $C=(c_{\beta\gamma})_{\beta,\gamma\in I}$ is the inverse of the matrix $(\langle\beta,\hat\alpha\rangle)_{\alpha,\beta\in I}$.
   It is known that $c_{\beta\gamma}\ge 0$ for all $\beta,\gamma\in I$.
   For $\alpha\in \Delta-I$ we now deduce from $\pi_I(w(\chi+\rho))\;\in\; -\overline{\check{\aaa}_I^+}$
   \begin{align*} 0\quad &\ge\;  \langle \pi_I(w(\chi+\rho)),\hat\alpha\rangle\;=\;
       \langle w(\chi+\rho),\hat\alpha\rangle -\sum_{\beta,\gamma\in I} \langle\beta,\hat\alpha\rangle\cdot
       c_{\beta\gamma}\cdot \langle w(\chi+\rho),\hat\gamma\rangle \\
       \;&\ge \;\langle w(\chi+\rho),\hat\alpha\rangle
       \end{align*}
       since $\langle\beta,\hat\alpha\rangle\le 0$ for $\alpha\notin I,\beta\in I$ and since
        $c_{\beta\gamma}\ge 0$ and  $\langle w(\chi+\rho),\hat\gamma\rangle>0$ for all $\beta,\gamma\in I$.
    By the computation at the beginning we conclude $w^{-1}\hat\alpha<0$, which immediately implies the claim, since we have
    $w^{-1}\hat\alpha >0$ for $\alpha\in I$ by the definition of $W^{I}$.
    \qed
    \absatz
This implies that for each $w\in W$ there is at most one $I\subset \Delta$ such that
  the condition $w\in W^{I}$ and the cone condition
  $\pi_I(w(\chi+\rho))\;\in\; -\overline{\check{\aaa}_I^+}$
  are both satisfied. We write $I=I(w)$ if both conditions are satisfied.

\Absatz\Numerierung\label{Franke fuer lineare} {\sc The case of linear groups.}
  In the case of the automorphism $\eta$ of order $2$ on $G=\GL_n,\GL_n\times \GL_1,\PGL_n$
  we can make things more explicit:
  Recall that for $G=\GL_n$ we denote by $(a_1,a_2,\ldots,a_n)\in \ZZ^n\simeq X^*(T)$ the character
  $\chi:diag(t_1,\ldots,t_n)\mapsto t_1^{a_1}\cdots t_n^{a_n}$. Similarly we denote
 a character of the group $\G_m^n\times \G_m$ by $(a_1,a_2,\ldots,a_n;a_0)\in \ZZ^n\times \ZZ$.
  The characters of the diagonal torus in $\PGL_n$ may be described by the set $\set{(a_1,\ldots,a_n)\in \ZZ^n |\sum_{i=1}^n a_i=0}$.
  The simple roots are of the form $\alpha_i=e_i-e_{i+1}\in X^*(T)$ for the standard basis $e_i$ of $\Z^n$.
  The positive Weyl chamber is given by the inequalities $a_1\ge a_2\ge \ldots \ge a_n$. The
  automorphism $\eta$ acts by $(a_1,a_2,\ldots,a_n)\mapsto (a_n,a_{n-1},\ldots,a_1)$ (cases $\GL_n$ and $\PGL_n$)
  respectively by $(a_1,\ldots,a_n;a_0)\mapsto (a_0-a_n,\ldots,a_0-a_1;a_0)$
  in the case $G=\GL_n\times \GL_1$.
  We have $\rho=\left(\frac{n-1}{2},\frac{n-3}{2},\ldots,\frac{1-n}{2}\right)$.
  This implies for the Weyl-group action of $w\in W=S_n$:
  \begin{align*}
     w(\chi+\rho)-\rho\quad&=\quad ( b_1,\ldots, b_n)\quad\text{
  with } b_i=a_{w^{-1}(i)}+i-w^{-1}(i)
  \end{align*}
  If the complement of $I\subset \Delta$ is written in the form
  $\Delta-I=\set{\alpha_{i_1},\alpha_{i_2},\ldots,\alpha_{i_{r-1}}}$ with $i_1<i_2<\ldots<i_{r-1}$,  we introduce
  the notations  $i_0=0, i_r=n,\;\;j_\mu=i_\mu-i_{\mu-1}$ and
  $\tilde I=(i_1|i_2-i_1|\ldots|n-i_{r-1})=(j_1|\ldots|j_r)$.
  Then the Levi group $M_I$, which consists of  block diagonal
   matrices,  is isomorphic to
  $\GL_{j_1}\times\ldots\times GL_{j_r}\subset \GL_n$  (with obvious modifications
  for the other linear groups under consideration).The parabolic $P_I$ consists of upper triangular block matrices.

  If we consider $\chi=(a_1,a_2,\ldots,a_n)$  as a character on $T\subset M_I$ we write it in the
  form $(a_1,\ldots,a_{i_1}|a_{i_1+1},\ldots,a_{i_2}|\ldots|a_{i_{r-1}+1},\ldots, a_n)$ and call the
  sequence of numbers between two $|$ a block. For $\chi$ in the positive Weyl chamber
  the condition $w\in W^{I}$ then means that the numbers in each block of $w(\chi+\rho)-\rho=( b_1,\ldots, b_n)$
  with $ b_i=a_{w^{-1}(i)}+i-w^{-1}(i)$
  are in a semi decreasing order: $ b_{i_\mu+1}\ge\ldots\ge  b_{i_{\mu+1}}$. The numbers in the blocks of
  $w(\chi+\rho)=(\tilde b_1,\ldots,\tilde b_n)$ are then in a strictly decreasing order.
  We have $\tilde b_i=\tilde a_{w^{-1}(i)}$, where $\tilde a_i= a_i+\frac {n+1}2 - i$ is a strictly decreasing sequence.

  With the  arithmetic means of the numbers in each block
  $$m_\mu=\frac{1}{j_\mu}\cdot\sum_{\nu=1}^{j_\mu} \tilde b_{i_{\mu-1}+\nu}$$
   we can now describe $$\pi_I(w(\chi+\rho))= (m_1,\ldots,m_1| m_2,\ldots,m_2| \ldots| m_r\ldots,m_r)$$
   The condition $\pi_I(w(\chi+\rho))\;\in\; -\overline{\check{\aaa}_I^+}$ is equivalent to the statement
  that the sequence of rational numbers $m_\mu$
  is semi increasing: $m_1\le m_2\le \ldots \le m_r$.
\absatz
  The $\tilde b_i$ forming a strictly decreasing sequence in each block we get the relation
   $\tilde b_{i_\mu}\le m_\mu\le m_{\mu+1}\le \tilde b_{i_\mu+1}$ which implies the description
  $I=\{ \alpha_i\in \Delta| \tilde b_i> \tilde b_{i+1}\}= \set{\alpha_i\in \Delta| w^{-1}(i)< w^{-1}(i+1)} $ as in the preceeding lemma.

\meinchapter{Partitions}\label{Ch Partitions}

\Absatz\Numerierung {\sc numbered and unordered partitions.} By a numbered partition $J=(J_1,\ldots,J_r)$ of a finite set $T$ we mean an $r$-tuple of pairwise disjoint subsets
$J_i$ such that $T=\cup_{i=1}^r J_i$. The underlying unordered partition is the set $\set{J_1,\ldots,J_r}$. Two numbered partitions
$J=(J_1,\ldots,J_r)$ and $J'=(J_1',\ldots,J_r')$ are called equivalent (in symbols $J\sim J'$) if
$\set{J_1,\ldots,J_r}=\set{J_1',\ldots,J_r'}$.

  \Absatz\Numerierung\label{VariPartitionen} {\sc Totally ordered sets}.
  For a finite totally ordered set $T$ (e.g. a finite subset of $\Q$) let $\eta=\eta_T$ be the unique order reversing involution.

  If $T$ and $K$ are totally ordered sets of the same cardinality, then there exists
    a unique order preserving bijection $\iota$ between $T$ and $K$, and then $\iota\circ \eta_T=\eta_K\circ \iota$.
    Then $\iota$ induces a bijection between the numbered resp. unordered partitions of $T$ and those of $K$.

\begin{Def}  Let $J=(J_1, ...,J_r)$ be a numbered partition of the finite totally ordered set $T$.
      \begin{itemize}
      \item[(a)] $J$ is called chronological numbered if $i<j$ implies $t_i<t_j$ for all $t_i\in J_i$ and $t_j\in J_j$.
       \item[(b)] $\eta(J)=(\eta(J_r),\ldots,\eta(J_1))$.
       \item[(c)] $J$ is called $\eta$-stable, if $\eta(J)\sim J$
       \item[(d)] $J$ is called $\eta$-fixed, if $\eta(J)=J$ as a numbered partition.
       \item[(e)] $J$ is called $\eta$-invariant, if $\eta(J_\mu)=J_\mu$ for $\mu=1,\ldots,r$.
       \item[(f)] $J$ is called $\eta$-admissible,  if there is $l\le \frac r2$ such that
       $\eta(J_\mu)=J_{r+1-\mu}$ for $\mu=1,\ldots,l$ and $\eta(J_\mu)=J_\mu$ for $\mu=l+1,\ldots, r-l$.
       \end{itemize}
 \end{Def}

 \Absatz\Numerierung {\sc The permutation associated to a partition.} It is clear, that for each numbered partition $J=(J_1,\ldots,J_r)$ of a finite ordered set $T$ there exists a unique chronological numbered partition $J'=(J_1',\ldots,J_r')$ such that
  $\# J_\mu =\# J_\mu'$ for $\mu=1,\ldots,r$. For example if $T=\{1,\ldots,n\}$ we can put $i_\mu=\# J_1+\ldots+\# J_\mu$
  for $\mu=0,\ldots,r$ and then we have $J_\mu'=\set{i_{\mu-1}+1,\ldots,i_\mu}$. Now we define $w_J\in S_T$ to be the unique
  permutation of $T$ such that $w_J$ induces an order preserving bijection from $J_\mu$ to $J_\mu'$ for all $\mu=1,\ldots,r$.

\Absatz\Numerierung {\sc The partition associated to a permutation}.
Now we can continue the considerations of (\ref{Franke fuer lineare}):
  Let $I\subset \Delta$ be given. Writing
  $\Delta-I=\set{\alpha_{i_1},\alpha_{i_2},\ldots,\alpha_{i_{r-1}}}$ with $i_1<i_2<\ldots<i_{r-1}$ and
   $i_0=0, i_r=n$ we can consider the chronological numbered partition $J_\mu'=\set{i_{\mu-1}+1,\ldots,i_\mu}$ of
   $T=\set{1,\ldots,n}.$ Let $\chi=(a_1,\ldots,a_n)\in \Z^n$ be a fixed dominant weight and let $w\in W^{I}$ be
   such that $\pi_I(w(\chi+\rho))\;\in\; -\overline{\check{\aaa}_I^+}$.

   We define $$J_\mu\quad =\quad \set{ j\in\set{1,\ldots,n}\; |\quad i_{\mu-1} <w(j)\le i_\mu}\quad
   =\quad w^{-1}(J_\mu')$$
   and get a numbered partition $J=(J_1,\ldots,J_r)$. Since $\tilde a_i= a_i+\frac {n+1}2 - i$ is a strictly decreasing sequence and since $\tilde b_i=\tilde a_{w^{-1}(i)}$ is a strictly decreasing sequence if $i$ varies in some $J_\mu'$, we get that $w$ induces an order preserving bijection from $J_\mu$ to $J_\mu'$ for each $\mu=1,\ldots,r$. Consequently
   $w$ coincides with the element $w_J$ associated to the partition $J$. Therefore $I$ and $w$ are uniquely determined by $J$.
   \klabsatz
   We may rewrite the arithmetic means in the following form
   $$m_\mu=m(J_\mu)=\frac 1{\# J_\mu'} \sum_{j\in J_\mu'} \tilde b_j=\frac 1{\# J_\mu} \sum_{j\in J_\mu} \tilde a_j.$$
   Thus we get that the numbered partition $J$ is admissible numbered in the following sense:
   \begin{Def}
   A numbered partition $J=(J_1,\ldots,J_r)$ of $\{1,\ldots,n\}$ is called admissible numbered (with respect to a dominant
   $\chi=(a_1,\ldots,a_n)$) if the arithmetic means
   $$m(J_\mu)=\frac 1{\# J_\mu} \sum_{j\in J_\mu}  \left(a_j+\frac {n+1}2-j\right)$$
   form a semi increasing sequence of rational numbers:  $$m(J_1)\le m(J_2)\le \ldots\le m(J_r).$$
   \end{Def}
\Absatz\Numerierung {\sc Reformulation of Frankes theorem.}
It is clear that each numbered partition is equivalent to some admissible numbered partition.
    Now let $\PPP$ be a set of admissible numbered partitions of $\set{1,2,\ldots,n}$ which contains for every unordered
    partition exactly  one admissible numbered representative. To an admissible numbered partition $J$ we can associate
    the element $w=w_J\in S_n$ and using the numbers $i_\mu$ defining the chronological ordered partition $J'$ we get the subset
    $I(J)=\Delta-\{ \alpha_{i_1},\ldots,\alpha_{i_{r-1}}\}$ of $ \Delta$ associated to $J$.

\begin{prop}\label{FrankeGLn}
     For $G=\GL_n$ or $G=\PGL_n$ or $G=\GL_n\times \GL_1$ we have the following identity in $\Groth(G(\A_f))$:
    \begin{align*}
           H^{*,G}(V_\chi)\quad&=\quad \sum_{J\in \PPP} \; (-1)^{l(w_J)}\; Ind_{P_{I(J)}(\A_f)}^{G(\A_f)}
              H^{*,M_{I(J)}}_{2,disc}(V_{w_J(\chi+\rho)-\rho}^{I(J)}).
       \end{align*}
\end{prop}
    Proof:  The preceding combinatorial considerations show that the sum over $I$ and $w\in W^{I}$
    in  theorem \ref{Frankes theorem} may be replaced by a sum over all admissible numbered partitions. But the numbers
    $n_{I(J)}(w_J(\chi+\rho)-\rho)$ count the admissible numberings of the unordered partition
    underlying $J$. Since the usual intertwining operator induces an isomorphism
    $$ \iota_{J,\tilde J}:\quad Ind_{P_{I(J)}(\A_f)}^{G(\A_f)}
              H^{*,M_{I(J)}}_{2,disc}(V_{w_J(\chi+\rho)-\rho}^{I(J)})\quad  \simeq\quad
              Ind_{P_{I(\tilde J)}(\A_f)}^{G(\A_f)}
              H^{*,M_{I(\tilde J)}}_{2,disc}(V_{w_{\tilde J}(\chi+\rho)-\rho}^{I(\tilde J)}),$$
              if $J$ and $\tilde J$ are two different admissible numberings of the same unordered partition, we get the claim.\qed

 \begin{lemma} Let $\chi\in X^*(T)^\eta$ be $\eta$-invariant, and let $J$ be an admissible numbered partition of
  $\set{1,\ldots,n}$.
  \begin{itemize}
     \item[(a)]  $m(\eta(J_\mu))=a_0-m(J_\mu)$ where $a_0=0$ in the cases $G=\GL_n$ and $G=\PGL_n$.
     \item[(b)]  $\eta(J)$ is an admissible numbered partition.
     \item[(c)] If $J$ is $\eta$-admissible, then it is $\eta$-stable.
     \item[(d)] Each $\eta$-stable admissible numbered partition is equivalent to an $\eta$-admissible numbered partition.
     \item[(e)] An $\eta$-stable $J$ is $\eta$-fixed if and only if it is $\eta$-admissible and there is at most
     one index $\mu$ with $\eta(J_\mu)=J_\mu$.
  \end{itemize}
   \end{lemma}

 Proof: (a) $\chi=(a_1,\ldots,a_n;a_0)$ resp. $\chi=(a_1,\ldots,a_n)$ being $\eta$-invariant means
      $a_{\eta(j)}=a_0- a_j$ for $j=1,\ldots,n$. This implies $$\tilde a_{\eta(j)}= a_{\eta(j)}+\frac {n+1}2 -\eta(j)
      = a_0-a_j+\frac {n+1}2-(n+1-j)=a_0-\tilde a_j$$
      $$ m(\eta(J_\mu))=\frac 1{\# \eta(J_\mu)} \sum_{j\in \eta(J_\mu)} \tilde a_j
       =\frac 1{\# J_\mu} \sum_{j\in J_\mu} \tilde a_{\eta(j)}=a_0-\frac 1{\# J_\mu} \sum_{j\in J_\mu} \tilde a_j
       =a_0-m(J_\mu).$$
       (b) is an immediate consequence of (a).
       \klabsatz
       (c) is clear from the definitions since $\eta$ is an involution on the set $\set{1,\ldots,n}$.
       \klabsatz
       (d) If $J$ is $\eta$-stable let $l$ be the number of $\eta$-orbits  $\{J_\mu,\eta(J_\mu)\}$ containing two elements inside
       $\set{J_1,\ldots,J_r}$.
       From each of these orbits we pick an element $J_\mu$ with $m(J_\mu)\le m(\eta(J_\mu))$ and we renumber these elements
       by the integers $1,\ldots,l$ in such a way that we have $m(J_1)\le\ldots\le m(J_l)$.
       We put $J_{n+1-\mu}=\eta(J_\mu)$ for $\mu=1,\ldots, l$ and we let $J_{l+1},\ldots,J_{n-l}$ be any numbering of the
      subsets satisfying  $\eta(J_\mu)=J_\mu$. Since (a) implies
       $m(J_\mu)=\frac {a_0}2$ for $\mu=l+1,\ldots,n-l$ and furthermore $m(J_l)\le \frac {a_0}2\le m(J_{n+1-l})$
       we conclude that we get an admissible numbering, which is
       $\eta$-admissible by construction.
       \klabsatz
       (e) is an immediate consequence of the definitions.\qed

  \Absatz\Numerierung
        We  associate to every $\eta$-admissible partition
        $$J=(J_1,\ldots,J_l,J_{l+1},\ldots, J_{r-l},J_{r-l+1},\ldots,J_r)$$
        an $\eta$-fixed partition $\pi(J)=(\tilde J_1,\ldots,\tilde J_\lambda)$. If $2l=r$ we put $\pi(J)=J$ and $\lambda=2l=r$.
        In the case $2l<r$ we put  $\tilde J_\mu=J_\mu$ for $\mu=1,\ldots,l$, \;   $\tilde J_{l+1}=\bigcup_{\mu=l+1}^{r-l} J_\mu $ and
        finally $\tilde J_\mu = J_{\mu+r-2l-1}$ for $\mu=l+2,\ldots, 2l+1=:\lambda $.

        \klabsatz
        If $\pi(J)$ is replaced by an equivalent $\eta$-fixed partition $\tilde J'$ then it is easy to see that there
        exists an $\eta$-admissible partition $J'$ equivalent to $J$ such that $\pi(J')=\tilde J'$.
        \klabsatz
        By eventually replacing an admissible partition by an equivalent one we may assume that
        the set of partitions $\PPP$ is formed in such a way, that
        the subset $\PPP^{stable}$ of $\eta$-stable partitions in $\PPP$ consists of $\eta$-admissible
       partitions and that $\pi$ maps $\PPP^{stable}$ to the subset $\PPP^{fix}$ of $\eta$-fixed partitions in $\PPP$.
       The subset of $\eta$-invariant partitions in $\PPP$ is denoted by $\PPP^{inv}$.
       Observe that $\PPP^{inv}=\set{ J\in \PPP^{stable} | \pi(J)= J_{trivial}}$ where $J_{trivial}$ is the trivial partition
       $(\{1,\ldots,n\})$.

       \klabsatz
       If $J$ is an $\eta$-admissible partition as above with $2l<r$, then $w_{\pi(J)}$ is an order preserving bijection from
       $\tilde J_{l+1}$ to $\set{i_l+1,\ldots,n-i_l}$. Then the shift $\sigma: i\mapsto i-i_l$  is an order preserving bijection
       of this latter set to $\set{1,\ldots,\tilde n}$ where $\tilde n= n -2i_l$.
       Now $K_\nu=\sigma(w_{\pi(J)}(J_{l+\nu}))$ for $\nu=1,\ldots, r-2l$ defines a
       numbered partition of $\set{1,\ldots,\tilde n}$. We write $K(J)=(K_1,\ldots,K_{r-2l})$.
       It is easy to see that $K(J)$ is admissible numbered with respect to the character $\tilde \chi$ which is obtained by
        subtracting the $\rho$ of $\GL_{\tilde n}$ from the middle block of the character
        $w_{\pi(J)}(\chi +\rho)$.

   \begin{lemma} \label{etalemma} Let $J=(J_1,\ldots,J_r)$ be an admissible numbered partition.
    \begin{itemize}
       \item[(a)] We have  $I(\eta(J))=\eta(I(J))$ and $\eta(w_J)=w_{\eta(J)}$.
       \item[(b)] If $J$ is $\eta$-fixed then we have $w(J)\in W^{\eta}$ and $I(J)$ is $\eta$-invariant.
       \item[(c)] $\eta(K_\nu)=K_\nu$ for $\nu=1,\ldots,r-2l$.
       \item[(d)] If $J$ is $\eta$-admissible then $w(J)\cdot w(\pi(J))^{-1}$ is the identity on
       $\set {1,\ldots,i_l}$ and on $\set{n-i_l+1,\ldots,n}$, and we have
       $$ \sigma\circ w(J)\cdot w(\pi(J))^{-1}\circ \sigma^{-1}\quad=\quad w_{K(J)}.$$
       \end{itemize}
    \end{lemma}
    Proof: (a) To  $\eta(J)=(\tilde J_1,\ldots,\tilde J_r)$ with $\tilde J_\mu=\eta(J_{r+1-\mu})$
      one associates the numbers $\tilde i_\mu=\# \tilde J_1+\ldots +\#\tilde J_\mu = \#J_r+\ldots+\#J_{r+1-\mu}
      = n-\left( \# J_1+\ldots +\# J_{r-\mu}\right) =n -i_{r-\mu}$.
      Then we have $\eta(I(J))=\eta\left(\Delta-\set{\alpha_{i_1},\ldots,\alpha_{i_{r-1}}}\right)
      =\Delta -\set{\alpha_{n-i_1},\ldots,\alpha_{n-i_{r-1}}}=\Delta-\set{\alpha_{\tilde i_{r-1}},\ldots,\alpha_{\tilde i_{1}}}
      =I(\eta(J))$
      since $\eta(\alpha_i)=\alpha_{n-i}$.
      With the involution $w_n:i\mapsto n+1-i$ of $\set{1,\ldots,n}$ we have
      $\eta(w)= w_n\circ w\circ w_n$ in $W=S_n$.
      Now $w_n$ maps $\tilde J_\mu$ bijectively and order reversing to $J_{r+1-\mu}$. Then $w_J$ maps by definition
      $J_{r+1-\mu}$ bijectively and order preserving to $\set{i_{r-\mu}+1,\ldots, i_{r-\mu+1}}$ and finally
      $w_n$ maps this set bijectively and order reversing to $\set{n+1-i_{r-\mu+1},\ldots,n+1-(i_{r-\mu}+1)}
      =\set{\tilde i_{\mu-1}+1,\ldots,\tilde i_\mu}$. Thus $\eta(w_J)$ has the defining properties of $w_{\eta(J)}$ and thus
      must be equal to this permutation.
    \klabsatz
    (b) is a consequence of (a).
    \klabsatz
    (c) Since $w_{\pi(J)}$ and $J_{l+\nu}$ are fixed by $\eta=\eta_n$ we have $\eta_n(w_{\pi(J)}(J_{l+\nu}))= w_{\pi(J)}(J_{l+\nu})$
    and from this one concludes $\eta_{\tilde n}(K_\nu)=K_\nu$.
    \klabsatz
    (d) again is a consequence of  the defining properties of $w_j,w_{\pi(J)}$ and $w_{K(J)}$.\qed

    \begin{lemma} If  $\tilde J$ is an $\eta$-fixed partition, then the map $J\mapsto K(J)$ defines a bijection
    from the inverse image of $\tilde J$ under $\pi$ in $\PPP^{stable}$ to the set $\PPP^{inv}_{\tilde n}$ of $\eta$-invariant partitions
    of $\set{1,\ldots,\tilde n}$.
    \end{lemma}
    \qed

    \Absatz\Numerierung
    Now let $G=\PGL_{2g+1}$ and $G_1=\Sp_{2g}$ be its stable $\eta$-twisted endoscopic group.
    We have $X^*(T_1)=X^*(T)^{\eta}$ and $W(G_1)=W(G)^{\eta}$. The set of simple roots $\Delta_1$ of
    $G_1$ may be identified with the set of $\eta$-orbits in $\Delta$. Therefore the
    subsets $I_1$ of $\Delta_1$ (which parametrize the standard parabolics  of $G$)
    may be identified with the $\eta$-stable subsets $I$ of $\Delta$. Denote by $P_{1,I}$ resp. $M_{1,I}$ the corresponding
     parabolic resp. Levi subgroups of $G_1$.

     \begin{prop}\label{Frankesp2g}
     For $G_1=\Sp_{2g}$  we have the following identity in $\Groth(G_1(\A_f))$, if $\chi\in X^*(T_1)=X^*(T)^{\eta}$ is dominant:
    \begin{align*}
           H^{*,G_1}(V_\chi)\quad&=\quad \sum_{J\in \PPP^{fix}} \; Ind_{P_{1,I(J)}(\A_f)}^{G_1(\A_f)}
              H^{*,M_{1,I(J)}}_{2,disc}(V_{w_J(\chi+\rho)-\rho}^{I(J)}).
       \end{align*}
\end{prop}
  Proof: This again is a consequence of theorem \ref{Frankes theorem}, if we take $W(G_1)=W(G)^\eta$ and
  lemma \ref{etalemma} into account. \qed

   \begin{thm}\label{diskretes Lifting theorem}
   For dominant $\chi\in X^*(T_1)=X^*(T)^{\eta}$ we have
   \begin{align*}
     H^{*,G_1}_{2,disc}(V_\chi)\qquad\text{lifts to}\qquad \sum_{J\in \PPP^{inv}} \;
     Ind_{P_{I(J)}(\A_f)}^{G(\A_f)}H^{*,M_{I(J)}}_{2,disc}(V_{w_J(\chi+\rho)-\rho})
     \end{align*}
   \end{thm}

   The proof is by induction on $g$, the case $g=0$ being trivial. So we assume that the theorem is valid
   for all $\gamma<g$.
   Now  \cite[Cor 5.27]{Uwe} as restated in \ref{Cohomology lifting theorem}(b) implies, that
      $H^{*,G_1}(V_\chi)$ lifts to $H^{*,G}(V_\chi)$. Write
      \begin{align*}
      R_J\quad &=\quad (-1)^{l(w_J)}\; Ind_{P_{I(J)}(\A_f)}^{G(\A_f)}H^{*,M_{I(J)}}_{2,disc}(V_{w_J(\chi+\rho)-\rho}^{I(J)}) \qquad\text { and }\\
      R_J^1\quad &= \quad (-1)^{l_1(w_J)}\; Ind_{P_{1,I(J)}(\A_f)}^{G_1(\A_f)}
              H^{*,M_{1,I(J)}}_{2,disc}(V_{w_J(\chi+\rho)-\rho}^{I(J)})
      \end{align*}
      By  \ref{Frankesp2g} and \ref{FrankeGLn} we may reformulate:
      $$\sum_{ J\in \PPP^{fix}} R_J^1\quad \text{ lifts to }  \quad \sum_{J\in \PPP} R_J.$$
      Here we have to observe that the filtration which induces Frankes spectral sequence is $\eta$-invariant.
      Thus we get an $\eta$-action on the initial term of the spectral sequence and therefore on $\sum_{J\in \PPP} R_J$ in such a way that it is compatible with the $\eta$-action on the limit.
      \absatz
      Since $R_J$ appears together with $R_{\eta(J)}$ in the sum for $G$, if $J$ is not $\eta$-stable (thus
      $R_J+R_{\eta(J)}$ has trivial $\eta$-traces), we
      may replace the sum on the right hand side by
      $$  \sum_{J\in \PPP^{stable}} R_J = \sum_{ J\in \PPP^{fix}}  \sum_{J'\in \PPP^{stable}, \pi(J')=J} R_{J'}.$$
      
      We will prove that the induction assumption implies that for $J\in \PPP^{fix}, J\ne J_{trivial}$ we have:
      $ R_J^1 $ lifts to $\sum_{J'\in \PPP^{stable}, \pi(J')=J} R_J$. Then it is an immediate consequence that
      $H^{*,G_1}_{2,disc}(V_\chi)= R_{J_{trivial}}^1$ lifts to $\sum_{J\in \PPP^{inv}} R_J$.

      \klabsatz
      Indeed let $J=(J_1,\ldots,J_r)$ be $\eta$-fixed and put $j_i=\# J_i$. Then $j_i=j_{r+1-i}$ implies that
      $r= 2l+1$ must be odd, since $\sum_{i=1}^r j_i= 2g+1$ is odd. Write $j_{l+1}= 2\gamma +1$. Then
      $M_{1,I(J)}$ is of the form $L_J\times \Sp_{2\gamma}$ with $L_J=\Pi_{i=1}^l \GL_{j_i}\subset Sp_{2g}$.
      Similarly the inverse image of $M_{I(J)}$ in $\GL_{2g+1}$ is of the form
      $\widetilde {M_{I(J)}}= L_J\times \GL_{2\gamma+1}\times \eta(L_J)$. For $J'\in \PPP^{stable}$ with $\pi(J')=J$
      we get $\widetilde {M_{I(J')}}= L_J\times M_{K(J')}\times \eta(L_J)$ in the notations introduced above.
      In the case $J\ne J_{trivial}$ we have $\gamma< g$ and $L_J\ne 1$ and we can apply the induction assumption
      to the lift from $\Sp_{2\gamma}$ to $\PGL_{2\gamma+1}$ with respect to the character $\tilde \chi$ introduced earlier.
      Observe that the $J'\in \PPP^{stable}$ with
      $\pi(J')=J$ are in one-to-one correspondence with the equivalence classes of $\eta$-invariant partitions
      $K$ of $\set{1,\ldots,,2\gamma+1}$.

      It is an easy exercise that the induction assumption implies that
      $ H^{*,M_{1,I(J)}}_{2,disc}(V_{w_J(\chi+\rho)-\rho}^{I(J)})$ lifts to
      $\sum_{\pi(J')=J} (-1)^{l(w_{J'})-l_1(w_J)} H^{*,M_{I(J')}}_{2,disc}(V_{w_{J'}(\chi+\rho)-\rho}^{I(J')})$, and then the claim is clear, since the lifting
      relation remains valid after parabolic induction.

\qed

\meinchapter{Proof of the Main Theorem}\label{Ch Proof}

  \Absatz\Numerierung \label{Blockcharaktere} {\sc Characters on the blocks.}
     Let  $\chi = (a_g,\ldots,a_1,0,-a_1,\ldots,-a_g)\in X^*(T)^\eta\subset \Z^{2g+1}$ be a dominant weight.
     Recall that $$S_\chi=\left\{ -a_g-g,\ldots,-a_1-1,0,a_1+1,\ldots,a_g+g\right\}$$ is the characteristic set of type $C_g$ associated to
     $\chi$.
     Let $J=(J_1,\ldots,J_r)$ be an $\eta$-invariant partition of $\{1,\ldots,2g+1\}$.  Via \ref{VariPartitionen} this gives rise to
     an $\eta$-invariant partition $S_\chi=\Sigma_1\cup\ldots\cup \Sigma_r$. We remark that every ordering of an $\eta$-invariant partition is admissible, but we have to fix one. For $i=1,\ldots,r$  let $j_i=\#(\Sigma_i)$ and
     $k_i=-j_1-\ldots-j_{i-1}+j_{i+1}+\ldots+j_r$. Let $T_i\subset \GL_{j_i}$ be the diagonal torus and let $\kappa_i:\G_m\to \G_m,\quad z\mapsto z^{k_i}$. Observe
     $ k_i\equiv \sum_{\nu\ne i} j_\nu =2g+1-j_i \equiv 1+\#(\Sigma_i) \mod 2$.
     \absatz
     Observe that $M_I(\A)\ni (m_1,\ldots,m_r)\mapsto \prod_{i=1}^r |det(m_i)|^{k_i}$ is the modulus of the action of
     $M_I(\A)$ on the unipotent radical $U_I(\A)$ of $P_I(\A)$ and thus its square root is the additional factor, which describes the normalized parabolic induction in terms of the naive induction. This square root is attached to the difference $\rho_G-\rho_{M_I}$.

     \absatz
     There is exactly one index $i$ such that $0\in \Sigma_i$. Then $j_i=2\gamma_i+1$ is odd and $\Sigma_i$ is of type $C_{\gamma_i}$.
     We have $\Sigma_i=S_{\tilde \chi_i}$ for some $\eta$-invariant dominant character $\tilde \chi_i$.
     Since $k_i$ is even in this case, we may put $\chi_i=\tilde \chi_i\cdot (\det)^{-k_i/2}$ and get an $\eta$-invariant dominant
     character $\chi_i\times \kappa_i^{-1}$ on $T_i\times \G_m$. If
     $\Sigma_i= \left\{ -b_{\gamma_i}-\gamma_i,\ldots,-b_1-1,0,b_1+1,\ldots,b_{\gamma_i}+\gamma_i\right\}$ with
     $b_{\gamma_i}\ge \ldots\ge b_1\ge 0$ then
     $$\chi_i\;=\;\left(b_{\gamma_i}-\frac{k_i}2,\ldots,b_1-\frac{k_i}2,
     -\frac{k_i}2,-b_1-\frac{k_i}2,\ldots,-b_{\gamma_i}-\frac{k_i}2\right)\in \Z^{2\gamma_i+1}.$$
     \absatz
     If $0\notin \Sigma_i$ then $j_i=2\gamma_i$ is even, $\Sigma_i$ is of type $D_{\gamma_i}$ and $k_i$ is odd.
     By lemma \ref{Klassifikation D}(b) there exists an $\eta$- invariant dominant character of the form $\chi_i\times \kappa_i^{-1}$ on $T_i\times \G_m$ such that its characteristic set is $\Sigma_i$.
     If $\Sigma_i=\left\{ -b_{\gamma_i}-\gamma_i,\ldots,-b_1-1,b_1+1,\ldots,b_{\gamma_i}+\gamma_i\right\}$
     with $b_{\gamma_i}\ge \ldots\ge b_1\ge 0$
     then
     $$\chi_i\;=\;\left(b_{\gamma_i}+\frac{-k_i+1}2,\ldots,b_1+\frac{-k_i+1}2,
     -b_1+\frac{-k_i-1}2,\ldots,-b_{\gamma_i}+\frac{-k_i-1}2\right)\in \Z^{2g_i}.$$
     In both cases we have $\eta(\chi_i)=\chi_i+(k_i,\ldots,k_i) =\chi_i\cdot \kappa_i\circ\det$.

\Absatz\Numerierung Now we want to relate the discrete spectrum cohomology of the Levi $M_I$ to the cohomology of
    its building blocks $\GL_{j_i}$: Observe that $K_\infty\cap M_I(\R)$ is not connected: We have
    $O_{2g+1}\cap M_I(\R)=\prod_{i=1}^r O_{j_i}$ and thus we may pick orthogonal 
    matrices with determinant $-1$ at an even number
    of places to get an element of $K_\infty^{I}=\SO_{2g+1}\cap M_I(\R)$.
    Therefore we have to take eigenspaces for the action of $O_{j_i}/SO_{j_i}$ on the cohomology of the building blocks as in
    \ref{pm Eigenraum}. If $j_i$ is odd, then $O_{j_i}/SO_{j_i}$ is represented by central matrices, i.e. the central character
    of a representation decides whether an eigenspace is empty or not.

   \begin{prop} In the notations of theorem \ref{diskretes Lifting theorem} and of the preceding sections we have

      $$H^{*,M_{I(J)}}_{2,disc}(V_{w_J(\chi+\rho)-\rho}^{I(J)})\quad=\quad \sum_{\epsilon=\pm} \bigotimes_{i=1}^r H^{*,\GL_{j_i}}_{2,disc}(V_{\chi_i})^\epsilon.$$

   \end{prop}
   Proof: It is a matter of book keeping, that we may write:
   $$w_J(\chi+\rho)-\rho\quad =\quad (\chi_1|\ldots| \chi_r).$$
    If one modifies the left hand side by replacing in  $H^*(\mmm_I,K_\infty^{I}A_I\Z_\infty,\ldots)$ the group 
    $K_\infty^{I}$ by its connected component $(K_\infty^{I})^\circ=\prod_i SO_{j_i}$, then the K\"unneth formula implies, that this is
    $$ \bigotimes_{i=1}^r H^{*,\GL_{j_i}}_{2,disc}(V_{\chi_i})\quad=\quad
     \bigotimes_{i=1}^r \sum_{\epsilon=\pm} H^{*,\GL_{j_i}}_{2,disc}(V_{\chi_i})^\epsilon.$$
     Now invariance under $\SO_{2g+1}\cap M_I(\R)$ means that only those tensor products of eigenspaces appear, which share the sign. 
   \qed

\begin{lemma}\label{Gewichte-Lemma} Let $\chi=(a_1,\ldots,a_{2g})\in X^*(T)=\Z^{2g}$ be a dominant character of the split torus $T$ in
      $\GL_{2g}$ satisfying $a_{2g+1-i}+a_i=w$ for some $w\in \Z$ and let $V_\chi$ be the associated
      algebraic representation with highest weight $\chi$.
      Then the center of $\GL_{2g}$ acts by the character $z\mapsto z^{gw}$.
\end{lemma}
Proof:  This is clear since we have $\sum_{i=1}^{2g} a_i= gw$ by our assumption.

\qed
\begin{prop}\label{Zentchar} Let $\pi_f=\pi_{1,f}\times\ldots\times \pi_{r,f}$ be a representation of $G(\A_f)$ contributing to
   $Ind_{P_{I(J)}(\A_f)}^{G(\A_f)}H^{*,M_{I(J)}}_{2,disc}(V_{w_J(\chi+\rho)-\rho})$ for a fixed $\eta$-invariant $J$ in the sense that for $i=1,\ldots,r$ we have an automorphic representation
    $\pi_i=\pi_{i\infty}\times \pi_{i,f}$  of $\GL_{j_i}(\A)$ in the discrete spectrum, and $\pi_i\otimes|\det|^{ k_i/2}$ contributes to $H^{*,\GL_{j_i}}_{2,disc}(V_{\chi_i})$.
   Assume $\eta(\pi_f)\cong\pi_f$.
   \begin{itemize}
   \item[(a)] We have $\eta(\pi_i)\cong \pi_i$ for $i=1,\ldots,r$.
    If $\pi_i=MW(\rho_i,n_i)$ with $j_i=m_i\cdot n_i$ and a cuspidal automorphic representation $\rho_i$ of
    $\GL_{m_i}(\A)$ we have $\eta(\rho_i) \cong \rho_i$.
   \item[(b)] The central character $\omega_i=\omega(\pi_i)=\omega_{i\infty}\cdot \omega(\pi_{i,f})$ is a quadratic character.
   \item[(c)] If $j_i=2\gamma_i+1$ is odd, then $\pi_i\otimes \omega_i$ is an $\eta$-invariant automorphic representation
   of $\PGL_{2g_i+1}(\A)$.
   \item[(d)] If $j_i=2\gamma_i$ is even, then $\omega_{i\infty}(-1)=(-1)^{\gamma_i}$.
   \end{itemize}
   \end{prop}
 Proof: (a) Let $\pi=\pi_1\times\ldots\times \pi_r$. By assumption
     $\pi$ is a quotient of
     $$\rho_1|\cdot|^{\frac{n_1-1}2}\times\ldots\times\rho_1|\cdot|^{\frac{1-n_1}2}\times\ldots\times
       \rho_r|\cdot|^{\frac{n_r-1}2}\times\ldots\times\rho_r|\cdot|^{\frac{1-n_r}2}.$$
      Similarly $\eta( \pi)$ is a quotient of
      $$\eta(\rho_r)|\cdot|^{\frac{n_r-1}2}\times\ldots\times\eta(\rho_r)|\cdot|^{\frac{1-n_r}2}\times\ldots\times
       \eta(\rho_1)|\cdot|^{\frac{n_1-1}2}\times\ldots\times\eta(\rho_1)|\cdot|^{\frac{1-n_1}2}.$$
     By the extended version of the strong multiplicity one theorem \cite[Theorem 4.4]{JacquetSmult} the assumption
     $\eta(\pi_f)\cong\pi_f$ implies that for every $i=1,\ldots r$ and every $\nu_i\le n_i-1$ there exists $j \le r$ and $\nu_j'\le n_j-1$ such that we have
     $ \eta(\rho_i)|\cdot|^{\frac{n_i-1}2-\nu_i} \cong \rho_j|\cdot|^{\frac{n_j-1}2-\nu_j'}$. This implies $m_i=m_j$.
     But $H^p(\gGg_{m_i},K_\infty^{m_i},\rho_{i,\infty}\otimes V_{\chi_i})\ne 0$ for some $p\ge 0$ implies
     $$0\ne H^{p}(\gGg_{m_i},K_\infty^{m_i},\eta(\rho_{i,\infty})\otimes {}^\eta V_{\chi_i}))=
     H^{p}(\gGg_{m_i},K_\infty^{m_i},\eta(\rho_{i,\infty})\otimes  (V_{\chi_i}\otimes \kappa_i\circ\det))$$ since ${}^\eta V_{\chi_i}=V_{\eta(\chi_i)}=V_{\chi_i}\otimes \kappa_i\circ\det$.
     This implies $0\ne 
     H^{p}(\gGg_{m_i},K_\infty^{m_i},(\eta(\rho_{i,\infty})\otimes |\det|^{k_i})\otimes  V_{\chi_i})$.
     Similarly $H^{p'}(\gGg_{m_i},K_\infty^{m_i},\rho_{j,\infty}\otimes V_{\chi_j})\ne 0$ for some $p'\ge 0$.
     But then the fact, that $\eta(\rho_i)$ and $\rho_j$ differ by a power of the modulus character, implies that $\chi_i$ and $\chi_j$ have the same characteristic set. Since $S_{\chi_i}\ne S_{\chi_j}$ for
     $i\ne j$ we get $i=j$. Then the bijectivity of the map $(i,\nu_i)\mapsto (j,\nu_j')$ implies immediately
     that we have $\nu_j'=0$ for $\nu_i=0$, and thus we get $\eta(\rho_i)\cong \rho_i$. From this we deduce  $\eta(\pi_i)\cong \pi_i$.

 \absatz
    (b) is a consequence of  (a) and $\omega(\eta(\pi_i))=\omega(\pi_i)^{-1}$.
 \klabsatz
    (c) We have $\omega(\pi_i\otimes \omega_i)=\omega(\pi_i)\cdot \omega_i^{2g_i+1}=\omega_i^{2(g_i+1)}=1$
    by (b).
 \klabsatz

    (d) Since $\pi$ contributes to the cohomology we have some non vanishing class in
    $H^{p}(\gGg_{j_i},K_\infty^{j_i},\pi_{i,\infty}\otimes V_\chi)$. Therefore we have some non vanishing
    $K_\infty$-equivariant homomorphism from $\Lambda^{p}(\gGg_{j_i}/\kkk^{j_i})$ to $\pi_{i,\infty}\otimes V_\chi$.
    Since  $-id\in \SO_{2\gamma_i}\subset K_\infty$ acts as identity on $\Lambda^{p}(\gGg/\kkk)$, as $\omega_{i\infty}(-1)$
    on $\pi_{i,\infty}$ and as $(-1)^{\gamma_i}$ on $V_\chi$ by \ref{Gewichte-Lemma} (since $k_i$ is odd), the claim is clear.

    \qed

\Absatz\Numerierung {\sc Proof of the Main Theorem \ref{Maintheorem}}. \label{ProofMT}

If $\tau$ is an irreducible  automorphic representation on $\SP_{2g}(\A)$ such that $\tau_f$ appears
 with non trivial multiplicity in the alternating sum
     $$\sum_{i} (-1)^{i} H^{i}(\mathfrak{sp}_{2g}, U_g,L_{2,disc}(\SP_{2g}(\Q)\backslash \SP_{2g}(\A))\otimes V_\chi ),$$
  then theorem \ref{diskretes Lifting theorem} tells us that there is a partition of
  $S_\chi$ into characteristic sets $\tilde S_i$ such that $\tau$ weakly lifts to a representation
  of the form $\pi=\pi_{1}\times\ldots\times \pi_{r}$ as in proposition \ref{Zentchar} with $\eta(\pi_f)=\pi_f$. Here $\pi_i$ are automorphic representations in
  the discrete spectrum, and by the main result of \cite{MoglinWald} they are of the form
  $\pi_i=MW(\rho_i, n_i)$ with positive integers $n_i$ and cuspidal representations $\rho_i$ on some $\GL_{m_i}$. Then
  $\tilde S_i=MW(S_i,n_i)$ for some $n_i$-admissible characteristic set $S_i$ of type $X_{\gamma_i}$,
  and the $\rho_i$ are up to twist by
  a power of the modulus character cohomological with respect to some coefficient system $V_{\chi_i}$ with $S_i=S_{\chi_i}$.
  By \ref{Zentchar}(a) we have $\eta(\rho_i)=\rho_i$. By \ref{Zentchar}(b) we can write $\omega_{\rho_i}=\chi_{d_i}$ for
  some $d_i\in \Q^*/(\Q^*)^2$. In the case
  that $S_i$ is of type $C_{\gamma_i}$ we put $\pi^{(i)}=\rho_i\otimes \chi_{d_i}$. In all other cases we put
  $\pi^{(i)}=\rho_i$. Again we have $\eta(\pi^{(i)})=\pi^{(i)}$.
   By the main result of \cite{Soudry Paris} (compare \cite[Thm. 6.1.]{Cogdell}) we get that $\pi^{(i)}$ is a weak lift of some irreducible generic cuspidal automorphic representation  $\pi_1^{(i)}$ on the group $\Sp_{2\gamma_i}$ in case that $S_i$ is of type $C_{\gamma_i}$,
  respectively on either the group $G_1^{(i)}=\SO_{2\gamma_i+1}$ or on the group $G_1^{(i)}=\SO_{2\gamma_i}^{d_i}$.
  Since $\pi_1^{(i)}$ has a semi weak lift to the $\GL_*$ groups in question by \cite{Cogdell} and since we have strong multiplicity one for $\GL_{2\gamma_i}$ we get, that $\pi^{(i)}$ is a semi weak lift of $\pi_1^{(i)}$.
  But now lemma \ref{Unterscheidungslemma} together with the description of
  $\LLL(\pi_{\infty}^{(i)})$ imply, that the group $G_1^{(i)}$ is of the correct type $X_{\gamma_i}$.
  The relation $\Pi_{i=1}^r d_i=1$ in $\Q^*/(\Q^*)^2$ is a consequence of $\Pi_{i=1}^r \chi_{d_i}=Pi_{i=1}^r \omega_i=\omega_\pi=1$.
  The remaining assertions are obvious or may be deduced  from \ref{Zentchar}(c) and (d).
  \qed
  \Absatz\Numerierung {\sc Proof of the Converse Theorem \ref{Umkehrtheorem}:}
  \label{ProofUmkehr} The proof consists in reversing the arguments of the proof of the Main Theorem. But one has to check that for each octupel satisfying the conditions the $\eta$-action on the corresponding part of the cohomology has a
  non trivial Lefschetz number:
  The cohomology is of the form
  $$ I= Ind_{P_I(\A_f)}^{G(\A_f)}\; H^*(\mmm_I,K_\infty\cap M_I(\RR);(\pi_1\times\ldots\pi_r)\otimes V_S)$$
  for a suitable coefficient system $V_S$. We can calculate in $\GL_{2g+1}$,
  such that $M_I=G^{(1)}\times\ldots\times G^{(r)}$ and such that $V_S=V_1\otimes\ldots \otimes V_r$. We may assume
  $G^{(1)}=\GL_{2\gamma_1+1}$ (i.e. type $C_{\gamma_1}$), so that $G^{(i)}=\GL_{2\gamma_i}$ for $i\ge 2$.
  For $i\ge 2$  there are elements  $k_i\in K_\infty\cap M_I(\RR)$ which are $-id$ in the first factor $G^{(1)}$, lie in $O_{2\gamma_i}-SO_{2\gamma_i}$
  in the component $i$ and are trivial in all other components. Then $K_\infty\cap M_I(\RR)$ is generated by
  the $k_i$ and by its connected component $(K_\infty\cap M_I(\RR))^\circ$.
  Let $\omega_1$ be the archimedean component of
  the central character of $\pi_1$. This is a quadratic character on $\R^*$ with $\epsilon=\omega(-1)=sign(d_1)\in \{\pm 1\}$.
  Since the center of $G^{(1)}$ acts trivially on $V_1$ we get that equivariance under $K_\infty\cap M_I(\RR)$ translates into
  $$ I=Ind_{P_I(\A_f)}^{G(\A_f)}\; \left(H^*(\gGg^{(1)},K_\infty^{(1)}, \pi_1\otimes V_1)\times H^*(G^{(2)})^\epsilon
   \times\ldots \times H^*(G^{(r)})^\epsilon\right),$$
   where $H^*(G^{(i)})^\epsilon=H^*(\gGg^{(i)},K_\infty^{(i)}, \pi_i\otimes V_i)^\epsilon$ are the $O_{2\gamma_i}/SO_{2\gamma_i}$
   eigenspaces as  in
   \ref{pm Eigenraum} for $i\ge 2$,  and where now $K_\infty^{(i)}$ are  connected subgroups.
   Then  propositions \ref{LefschetzC} and \ref{LefschetzBD} imply, that the Lefschetz number of
   $\eta\times h_f$ equals $\pm 2^{g-(r-1)}$ times the trace of $\eta h_f$ on $Ind_{P_I(\A_f)}^{G(\A_f)} \pi_{1,f}\times \pi_{r,f}$.
   Since this is not identically zero the usual trace formula arguments give the claim.

   \qed

   Remark: The number $\pm 2^{g-(r-1)}$ should be the multiplicity, with which an individual $\Sp_{2g}(\A_f)$ module
   contributes to the alternating sum of the cohomology groups.

\meinchapter{Examples and Complements}\label{Ch Examples}

\Absatz\Numerierung {\sc Type $C_1$.} For $k\ge 0$ the characteristic set $S=\set{-k-1,0,k+1}$
        is of type $C_1$. It corresponds to the coefficient system
          $Sym^k$ on $G_1=\SP_2=\SL_2$.
           $\pi_1$ is cohomological with respect to $S$ if $ \pi_{1,\infty}\subset D(k+1)$. This means that $\pi_1$ corresponds to holomorphic (resp. antiholomorphic) cusp forms
          on $\Sp_2=\SL_2$ of weight $k+2$. On $G=\PGL_3$ we get the adjoint lift
          $\pi=Ad(\pi_1)$  in the sense of \cite{GelbartJacquet} (compare \cite{FPGL3global} for a trace formula approach).
                    We have $\pi_\infty=D(0,2k+1)\times \sigma_{0,1}$.
          A  cuspidal automorphic representations
           $\pi_1$ of $G_1(\A)$ is a subrepresenations of some cuspidal automorphic representation $\tilde \pi_1$ of $\GL_2$,
           which is unique up to character twists (\cite{Ramakrishnan}), and we may write $\pi=Sym^2(\tilde\pi_1)\otimes \omega_{\tilde\pi_1}^{-1}$  (independent of the choice of $\tilde \pi_1$).
                      $\pi$ is irreducible cuspidal unless $\tilde\pi_1$ is a CM-form, i.e.
          unless $\tilde\pi_1$ is a Weil representation of a Gr\"o\ss encharacter.

 \Absatz\Numerierung {\sc Type $D_1$.} For $k\ge 0$ the characteristic set $S=\set{-k-1,k+1}$ is of type $D_1$. Since $(-1)^1d>0$ we get
      an imaginary quadratic extension $F=\Q(\sqrt d)$ of $\Q$. The group
        $G_1=\SO_2^d$ is the kernel of the norm map $Res_{F/\Q}\G_m\to \G_m$, so that $G_1(\R)=\set{z\in \C| z\bar z=1}$.
         Then $\pi_{1,\infty}(z)=z^{k+1}$. The character $\theta:\A_F^*/F^*\to \C^*,  a\mapsto \pi_1(a\cdot \bar a^{-1})$
         is a  Gr\"o\ss encharacter with archimedean component  $\theta_\infty(z)= (z/|z|)^{2k+2}$.
         Then $\pi=\WWW(\theta)$ is the Weil representation attached to this Gr\"o\ss encharacter. We have $\pi_\infty=D(0,2k+1)$.
         Since $\theta(a)=1$ for $a\in \A_\Q/\Q^*$ we get $\omega_\pi=\chi_d$.
         \absatz
         We remark that

\Absatz\Numerierung {\sc Type $B_1$.} For $k\ge 0$ the characteristic set $S=\set{-k-\frac 12,k+\frac 12}$ is of type $B_1$.
        The elementary pairs are  cuspidal automorphic representations $\pi_1$ of $G_1=\SO_3\cong \PGL_2$,
        which may be viewed as representations $\pi$ of $G=\GL_2$. They are cohomological with respect to $S$ if
         $\pi_\infty= D(0,2k+1)$, and then they may be viewed as classical cusp forms of weight $2k+2$ with trivial multiplier.

 \Absatz\Numerierung {\sc Type $D_2$.}\label{BspD2} For $a\ge b\ge 0$ the characteristic set $$S=\set{ -a-2,-b-1,b+1,a+2}$$
       of type $D_2$ gives rise to a coefficient system,
        which pulls back to the system $Sym^{a+b+2}\otimes Sym^{a-b}$   under the exceptional isogeny
       $i:\SL_2\times \SL_2\to \SO_4$ (for the coefficient systems one can consider this over $\bar\Q$).  In the case $d=1$
       we can write $i^*(\pi_1)=\sigma_1\otimes\sigma_2$, where the cuspidal representations $\sigma_i$ of $\SL_2(\A)$
       are contained in cuspidal representations $\tilde\sigma_i$ of $\GL_2(\A)$ and we may assume
        $\omega_{\tilde\sigma_1}\cdot \omega_{\tilde\sigma_2}=1$. Then the automorphic representation $\pi$ on
        $\GL_4(\A)$ is of the form $\pi=\tilde\sigma_1\boxtimes\tilde\sigma_2$ in the sense of \cite{Ramakrishnan}.
        We have $\sigma_{1,\infty}=D(a+b+3)$ and $\sigma_{2,\infty}=D(a-b+1)$ and thus
        $\pi_1$ corresponds to a pair of classical cusp forms $(f,g)$
       of weights $a+b+4$ and $a-b+2$ respectively.
       In view of the different archimedean types $\tilde\sigma_2$ is not a twist of $\tilde\sigma_1$ and therefore
       $\tilde\sigma_1\boxtimes\tilde\sigma_2$ is a cuspidal representation of $\GL_4$ unless $\tilde \sigma_1$ and $\tilde \sigma_2$
       are both CM-representations with respect to the same imaginary quadratic extension $E/\Q$. This is a consequence of the cuspidality criterion of \cite[ch. 3]{Ramakrishnan}: if $\tilde\sigma_1=\WWW(\theta)$ for a Gr\"o\ss encharacter $\theta$
       on $\A_E$, but $\tilde\sigma_2$ has no CM by $E$,
       then the base change $\tilde\sigma_{2,E}$ of $\tilde\sigma_2$ is cuspidal and it admits no self twist
       with $\bar\theta\cdot \theta^{-1}$, since this latter character has infinite order.
       \klabsatz
       In the case $d\ne 1$ we have $\SO_4^d\cong Res_{F/\Q}\SL_2/\{\pm1 \}$ with $F=\Q(\sqrt d)$ real quadratic, so that
       $\pi_1$ corresponds to a Hilbert modular form of weight $(a+b+4,a-b+2)$.
       \absatz
       In the case $a=b\ge 0$ we can write $$\set{ -a-2,-a-1,a+1,a+2}=MW(\set{-a-\frac 32, a+\frac 32}, 2)$$and we get
        classes in the discrete spectrum, which are described by a holomorphic cusp form of weight $2a+4$ and the residue of the weight $2$ Eisenstein series.
\Absatz\Numerierung
{\sc The Case $g=2$.} For $a\ge b\ge 0$  consider the characteristic set of type $C_2$:
      $$ S\quad=\quad \set{ -a-2,-b-1,0,b+1,a+2}.$$
 (a) The trivial partition $S=S$ corresponds to cuspidal representations of $\SP_4(\A)$, which lift to irreducible cuspidal
  representations of $\PGL_5(\A)$.
  \klabsatz
  (a') In the case $a=b=0$ we can write $S=MW(\{0\},5)$ (where $\set{0}$ is of type $C_0$) and get the one dimensional representations in the residual spectrum on $G_1$ and on $G$.

 \absatz
  (b) The partition $S=\set{0}\cup S_2$, where $S_2=\{ -a-2,-b-1,b+1,a+2\}$ is of type $D_2$, corresponds to cuspidal
  automorphic representations $\tau$ of $\SP_4(\A)$, which are endoscopic lifts of representations $\pi_1$ of
  $\SO_4^d$ with $d>0$. In case $d=1$  the representation $\pi_1$ is the restriction of
  a cuspidal representation $\tilde\sigma_1\times\tilde\sigma_2$ of $\GL_2\times\GL_2/\G_m=\GSO_4$ to $\SO_4$, where  $\omega_{\tilde\sigma_1}\cdot\omega_{\tilde\sigma_2}=1$.
  Then $\pi= 1\times (\tilde\sigma_1\boxtimes \tilde\sigma_2)$ on $\PGL_5(\A)$. Here $\tilde\sigma_1\boxtimes \tilde\sigma_2$ is cuspidal, if
   $\tilde\sigma_1$ and $\tilde\sigma_2$ are not $CM$-forms for the same imaginary quadratic extension of $\Q$.
   $\tau$ is the restriction to $\SP_4$ of an endoscopic lift of $\tilde\sigma_1\times\tilde\sigma_2$ to $\GSp_4$ in the sense of
   \cite{WeiBuch}.
  If $d\ne 1$ then $\tau$ is an endoscopic lift to $\SP_4$ of a cuspidal representation $\pi_1$ on the quasisplit $\SO_4^d\cong Res_{F/\Q}\SL_2/\{\pm1\}$, but not a restriction of an endoscopic representations of $\GSp_4$.
  \klabsatz
 (b') For $a=b$ one can write $S=\set{0}\cup MW(\{-a-\frac 32,a+\frac 32\},2)$ and gets $\pi=1\times MW(\pi^{(2)},2)$,
  which lies in the residual spectrum on the $\GL_4$ factor. Then $\tau$ is either a CAP-representation with respect to the Siegel parabolic i.e. of  Saito-Kurokawa type \cite{Saito-Kurokawa} or it is a residual representation of $\Sp_4$, whose contribution
  to cohomology is described in the work of Schwermer \cite{SchwSieg1}, \cite[Prop. 4.6.(2)]{SchwSieg2}. This case may be viewed as a version of (b), where $\tilde \sigma_1$ is a one dimensional (residual) representation.
 \absatz
 (c) The partition $S=S_1\cup S_2$ with $S_1=\set{-a-2,0,a+2}$ of type $C_1$ and $S_2=\set{-b-1,b+1}$ of type $D_1$ corresponds to
    $\tau$ on $\SP_4$, which are restrictions of automorphic representations $\tilde \tau$  of $\GSp_4$, whose
  associated rank 4 motive is a  tensor product of two  $\GL_2$ motives, one of them being $CM$: let $\pi_1^{(1)}$ be a
  cuspidal representation of $\GL_2$, which is not of $CM$ type,  with $\pi_{1,\infty}^{(1)}=D(a+2,0)$ and let $ \pi^{(1)}$ be its
  adjoint lift to a cuspidal representation of $\PGL_3$. Let
  $\sigma=\pi_1^{(2)}$ be an automorphic character on the norm 1 group of  an imaginary quadratic extension $F=\Q(\sqrt d)/\Q$
  with $\sigma_\infty(z)=z^{b+1}$ and let
  $\pi^{(2)}=\WWW(\theta)$ be the cuspidal Weil representation of  $\GL_2$, where
   $\theta(x)=\sigma(x\cdot \bar x^{-1})$. Then $\tau$ lifts to $(\pi^{(1)}\otimes \chi_d)\times \WWW(\theta)$ on $\PGL_5$.
   If one extends $\sigma$ to a Grossencharacter $\tilde\sigma$ of $\A_F^*/F^*$, one can arrange $\tilde \tau$ on $\GSp_4=\GSPIN_5$
   in such a way that it lifts to $\pi_1^{(1)}\boxtimes \WWW(\sigma)\otimes \chi$ on $\GL_4\times \G_m$, where $\chi$ is the product
   of the central character of $\pi_1^{(1)}$ and of the restriction of $\tilde \sigma$ to $\A_Q^*/\Q^*$.
   In special cases these representations have been constructed by Ramakrishnan and Shahidi \cite[Theorem B']{RamaSiegel}.
 \klabsatz
 (d) The partition $S=S_1\cup S_2$ with $S_1=\set{-b-1,0,b+1}$  of type $C_1$ and $S_2=\set{-a-2,a+2}$ of type $D_1$ is completely analogous to the case (c), but with different weights.
 \klabsatz
 (d') If $b=0$ in case (d) we may write $S_1=MW(\set{0},3)$. Then $\tau$ is a CAP-representations with respect to the Klingen parabolic
 (\cite{Soudry Crelle}).
 \absatz
 (e) The partition $S=\set{0}\cup\set{-a-2,a+2}\cup\set{-b-1,b+1}$ may be viewed as a degenerate case of the preceeding cases.
 If one has $\sigma_1=\WWW(\theta_1)$ and $\sigma_2=\WWW(\theta_2)$ in case (b), where $\theta_1$ and $\theta_2$ are Gr\"o\ss encharacters
 on $\A_F^*/F^*$ for the same imaginary quadratic extension $F/\Q$, such that $\theta_1\theta_2$ is trivial on
 $\A_\Q^*/\Q^*$ then we have $\sigma_1\boxtimes \sigma_2= \WWW(\theta_1\theta_2)\times \WWW(\theta_1\bar\theta_2)$.
 If $\pi_1^{(1)}$ in case (c) or (d) is a Weil representation, we get another way to interpret this degeneracy.

  \Absatz\Numerierung\label{g3} {\sc The case $g=3$.}
   Let $a\ge b\ge c\ge 0$ be integers. Then we get $3$ partitions of the characteristic set
      $$S=\set{-a-3,-b-2,-c-1,0,c+1,b+2,a+3}$$
      of type $C_3$ into a characteristic set $S_1$ of type $C_1$ and a set $S_2$ of type $D_2$.
      We get the following weights of classical modular forms (compare \cite[7.10.]{BFGeer}), which can be lifted to Siegel modular forms:
      \begin{align*}
         S_1=\set{-a-3,0,a+3} \text{ of weight } a+4,\quad &\text{ and }S_2 \text{ of weights } b+c+4, \; b-c+2\\
       S_1=\set{-b-2,0,b+2} \text{ of weight } b+3,\quad  &\text{ and }S_2 \text{ of weights } a+c+5, \; a-c+3\\
       S_1=\set{-c-1,0,c+1} \text{ of weight } c+2,\quad  &\text{ and }S_2 \text{ of weights } a+b+6,\; a-b+2.\\
      \end{align*}

\Absatz\Numerierung {\sc The Ikeda Lift.} For $k\ge 0$ classical scalar valued holomorphic Siegel modular forms
          of weight $k+g+1$ on $\Sp_{2g}$
          (i.e. with automorphy factor  $det(CZ+D)^{k+g+1}$) correspond to cohomology classes with the coefficient system
          $(k,\ldots,k,-k,\ldots,-k)\in X^*(T_H)=\Z^g$. The characteristic set (of type $C_g$) is
          $$ S=\set{-k-g,\ldots,-k-1,0,k+1,\ldots,k+g}.$$

       If $g=2\gamma$ is even then we have a partition
       $S\;=\;\set{0}\cup  MW\left(S_1,2\gamma\right),$
       where $S_1=\left\{-k-\gamma-\frac12,k+\gamma+\frac 12\right\}$ corresponds to classical holomorphic forms of weight
       $2k+2\gamma+2$. This gives an interpretation of  the Ikeda-Lift \cite{Ikeda}\cite{Kohnen},
       and in the special case $\gamma=1$ we again get the Saito-Kurokawa-lift \cite{Saito-Kurokawa}.

\Absatz\Numerierung {\sc Relations with elliptic endoscopy of $\GSp_{2g}$.} Next we try to understand the relation of
        theorems \ref{Maintheorem}, \ref{Umkehrtheorem} with the work of Morel \cite{Morel}:
        If $S=\bigcup_{i=1}^r MW(S_i,n_i)$ is a decomposition
        and  $d_i\in \Q^*/(\Q^*)^2$ a family as in theorem \ref{Maintheorem}, let $I_0=\{1,\ldots,r\}$, assume $S_1$ is of type
        $C_{\gamma_1}$ and let $(\pi^{(i)},\pi_1^{(i)})_{i\in I_0}$ be a family of elementary particles, coming from an automorphic representation $\tau$ of $\SP_{2g}(\A)$.

        The set $\Sigma$ of all subsets $I$ of $I_0$ with $1\notin I$ is a group under the symmetric difference $\Delta$ with neutral element $\emptyset$. The kernel of the homomorphism
        $\pi: \Sigma\to \Q^*/(\Q^*)^2,\; I\mapsto \prod_{i\in I} d_i$ is thus
        $$\Sigma_\Q\quad=\quad \left\{ I\subset I_0| 1\notin I,\; \Pi_{i\in I} d_i=1\right\}.$$ For $I\in \Sigma_\Q$
        we can put $S_I=\bigcup_{i\notin I} MW(S_i,n_i)$ and $S_I'=\bigcup_{i\in I} MW(S_i,n_i)$.
        Then $S_I$ is of type $C_{g_1}$ and $S_I'$ is either empty or of type $D_{g-g_1}$ if $g_1<g$.
        It is known that the elliptic endoscopic groups of $\GSp_{2g}$ are of the form $G_{g_1}=G(\SP_{2g_1}\times \SO_{2(g-g_1)})$
        (\cite[Prop. 2.1.1]{Morel}).
        Now for each $I\in \Sigma_\Q$ their exists by  \ref{Umkehrtheorem} an automorphic representation $\tau_I$ on
        $\SP_{2g_1}(\A)$, such that $\tau_I$ lifts to $\Pi_{i\notin I} MW(\pi^{(i)},n_i)$ and their should exist (comp. \cite{Arthur}) an automorphic representation
        $\tau_I'$ on $\SO_{2(g-g_1)}(\A)$, which lifts to $\Pi_{i\in I} MW(\pi^{(i)},n_i)$ . Then  $\tau_I\otimes \tau_I'$ is a subrepresentation of some
        $\pi_I$ on $G_{g_1}(\A)$, which should be related to $\tau$ in the sense of elliptic endoscopy.
        We remark that $g-g_1$ is even, since $d_i<0$ if $S_i$ is of type $D_{\gamma_i}$ with  odd $\gamma_i$.
        \Absatz
        Let $D\subset \Q^*/(\Q^*)^2$ be the image of $\pi$ and let $E:=\Q(\sqrt D):=\Q(\sqrt{d_2},\ldots,\sqrt{d_r})$ be the corresponding field extension of degree $\#(D)$.
        The Galois group $\Gamma_{E/\Q}=\Gamma_\Q/\Gamma_E$ is canonical dual to $D$ in the category of $\F_2$ vector spaces:
        $\sigma(d)=\sigma(\sqrt d)/\sqrt d\in \set{\pm 1}$. Now  $\Sigma$ is self dual with respect to the bilinear form
        $\beta(I,J)=(-1)^{\#(I\cap J)}$ and the canonical map $\pi^\vee: \sigma\mapsto \set{i\in I_0| i\ne 1 \text{ and }
        \sigma(\sqrt{d_i})=-\sqrt{d_i}}$ is dual to $\pi$ with respect to these pairings.
        \klabsatz
        One expects for each $i\in I_0$ the existence of a $\lambda$-adic representation $V_i$ of $\Gamma_\Q$ of dimension
        $2^{\gamma_i}$ such that the Frobenius eigenvalues are related to the Satake parameters of the automorphic representation $\pi_1^{(i)}$. For $i\ge 2$ the restriction of this representation to $\Gamma_{E_i}$ where $E_i=\Q(\sqrt{d_i})$ should
        split into two representations $V_i^{+}$ and $V_i^{-}$ of dimension $2^{\gamma_i-1}$ in such a way that we have
        $V_i=Ind_{\Gamma_{E_i}}^{\Gamma_\Q} V_i^{+}=Ind_{\Gamma_{E_i}}^{\Gamma_\Q} V_i^{-}$ if $d_i\ne 1$.
        For $I\in \Sigma$ we can form
        $$V_I\;=\; Ind_{\Gamma_{E}}^{\Gamma_\Q}\; (V_1\otimes V_2^{\epsilon_2}\otimes \ldots\otimes V_r^{\epsilon_r}),$$
        where $\epsilon_i=-$ for $i\in I$ and $\epsilon_i=+$ for $i\notin I$.
        These are representations of dimension $2^g/2^{r-1}\cdot [E:\Q]= 2^g/\#(\Sigma_\Q)$.
        Then we have $V_I\cong V_{I'}$ if $I \equiv I' \mod \pi^\vee(\Gamma_E)$ and finally
        $ V_1\otimes\ldots \otimes V_r\;=\; \bigoplus_{I\in \Sigma/\pi^\vee(\Gamma_E)} V_I$.
        Here the index set $\Sigma/\pi^\vee(\Gamma_E)$ is canonical dual to $\Sigma_\Q=\ker(\pi:\Sigma\to D)$.
        Now $\tau$ should be associated to one of the $\lambda$-adic representations  $\Sigma_I$, and the $\tau'$ in the packet of
        $\tau$ may correspond to $V_{I'}$ for different $I'$.

\Absatz\Numerierung
        In the notations of the example \ref{g3} ($g=3$) let $S=S_1\cup S_2\cup S_3\cup S_4$ with
        $S_1=\set{0}, \; S_2=\{a+3,-a-3\},\; S_3=\{b+2,-b-2\},\; S_4=\{c+1,-c-1\}$ and $d_1=d_2=d_3=d_4<0$.
        Then we have $\Sigma=\{\emptyset,\{2,3\},\set{2,4},\set{3,4}\}$, which means that we have to combine two elementary particles to get an elliptic endoscopic representation. This should lead to motives of rank 2 in the cohomology of the associated Shimura variety, but this phenomenon does not appear in the totally unramified case considered in \cite{BFGeer}.

\appendix
\renewcommand{\Numerierung}{\refstepcounter{meinzaehler}{\bf(\Alph{chapter}.\arabic{meinzaehler}) }}
\renewcommand{\meinchapter}[1]{\refstepcounter{chapter}\section*{\Alph{chapter} $\ $ #1}
    \addcontentsline{toc}{chapter}{\Alph{chapter}$\ $ #1}
    \markboth{}{\scshape #1}
    }
{ 
\setlength{\itemsep}{0pt}
\setlength{\parsep}{0pt}
\setlength{\topsep}{0pt}
\setlength{\partopsep}{0pt}
\setlength{\baselineskip}{0pt}
\setlength{\smallskipamount}{0pt}
\setlength{\medskipamount}{0pt}
\setlength{\bigskipamount}{0pt}

}

Dr. Uwe Weselmann

Mathematisches Institut

Im Neuenheimer Feld 288
\smallskip

D-69121 Heidelberg

\medskip
weselman@mathi.uni-heidelberg.de

\Seitenumbruch


\begin{thebibliography}{RohlfsSchw} 
%
\addcontentsline{toc}{chapter}{Bibliography}
%
\bibitem[Ar12]{Arthur} J. Arthur, {\it The Endoscopic Classification of Representations: Orthogonal and Symplectic Groups}
             preliminary version available on http://www.claymath.org/cw/arthur/pdf/arthur-book-2012b.pdf
\bibitem[BWW02]{FLPGL5} J. Ballmann, R. Weissauer, U. Weselmann, {\it Remarks on the fundamental lemma
            for stable twisted Endoscopy of classical groups},
            Manuskripte der Forschergruppe Arithmetik {\bf 7} (2002) Mannheim--Heidelberg.
            arXiv:1302.0034
\bibitem[BFG11]{BFGeer} J. Bergstr\"om, C. Faber, G. van der Geer,
               {\it Siegel modular forms of degree three and the cohomology of local systems}
              arXiv:1108.3731v2 To appear in Selecta Mathematica
\bibitem[Bo53]{BorelKohom} A. Borel, {\it Sur La Cohomologie des Espaces Fibres Principaux et des Espaces Homogenes
              de Groupes de Lie Compacts}, Annals of Math. {\bf 57} (1953), 115--207.
\bibitem[BoW80]{BorelWallach}A. Borel, N. Wallach, {\it Continous Cohomology, Discrete Subgroups and Representations of
              Reductive Groups}, Annals of Mathematics Studies {\bf 94}, Princeton University Press (1980)

\bibitem[Clo]{Clozel} L. Clozel, {\it Motifs et Formes Automorphes: Applications du Principe de Functorialit\'{e}},
            In: L. Clozel, J.S. Milne: Automorphic Forms, Shimura Varieties, and $L$-functions I, Perspectives in Mathematics
            {\bf 10}, Academic Press (1990) 77--159.
\bibitem[CKPSS01]{Cogdell01} J.W. Cogdell, H. Kim,  I.I. Piatetski-Shapiro, F. Shahidi, {\it On lifting from classical
             groups to $\GL_N$}, Publ. Math. IHES {\bf 93} (2001), 5--30.
\bibitem[CKPSS04]{Cogdell04} J.W. Cogdell, H. Kim,  I.I. Piatetski-Shapiro, F. Shahidi, {\it Functoriality for the classical groups},
          Publ. Math. IHES {\bf 99} (2004), 163--233.
\bibitem[CPSS11]{Cogdell} J.W. Cogdell, I.I. Piatetski-Shapiro, F. Shahidi, {\it Functoriality for the quasisplit classical groups}
            In: Arthur, James (ed.) et al., On certain L-functions. Conference in honor of Freydoon Shahidi on the occasion of his 60th birthday,  Clay Mathematics Proceedings {\bf 13} (2011) 117--140.

\bibitem[Fl94]{FPGL3global} Y. F. Flicker, {\it On the Symmetric Square: Total Global Comparison},
             J. Funct. Anal. {\bf 122} (1994), 255--278.
\bibitem[Fr98]{Franke} J. Franke, {\it Harmonic Analysis in weighted $L_2$-spaces}, Ann. scient. \'{E}cole
            Norm. Sup. {\bf 31} (1998), 181--279.
\bibitem[FS98]{FrankeSchwermer} J. Franke, J. Schwermer, {\it A decomposition of spaces of automorphic forms and the Eisenstein cohomology of arithmetic groups}, Math. Ann. {\bf 311}, (1998), 765--790.
\bibitem[GJ78]{GelbartJacquet} S. Gelbart, H. Jacquet, {\it A relation between automorphic representations
           of $GL(2)$ and $GL(3)$. },  Ann. Sci. \'{E}c. Norm. Sup\'{e}r. (4)
             {\bf 11},   (1978), 471--542.
\bibitem[GRS97]{Ginzburg} D. Ginzburg, St. Rallis, D. Soudry, {\it Periods, poles of $L$-functions and
            symplectic-orthogonal theta lifts} J. reine angew. Math. {\bf 487} (1997), 85--114.
\bibitem[GHV73]{Greub2} W. Greub, St. Halperin, R. Vanstone, {\it Connections, Curvature and Cohomology, Vol II:
              Cohomology of  Principal Bundles and Homogeneous Spaces}, Academic Press, New York London (1973).
\bibitem[GHV76]{Greub3} W. Greub, St. Halperin, R. Vanstone, {\it Connections, Curvature and Cohomology, Vol III:
              Lie Groups, Principal Bundles and Characteristic Classes}, Academic Press, New York London (1976).
\bibitem[GR11]{Grobner} H. Grobner, A. Raghuram, {\it On some arithmetic properties of automorphic forms
            on $\GL_m$ over a division algebra}, preprint 2011,  arXiv:1102.1872
\bibitem[Ik01]{Ikeda} T. Ikeda, {\it On the lifting of elliptic modular forms to Siegel cusp forms of degree $2n$}, Ann. of Math.
              (2) {\bf 154} (2001), 641--681.
\bibitem[Hum72]{Humphreys} J. E. Humphreys, {\it Introduction to Lie Algebras and Representation Theory},
             Grad. Texts in Math. {\bf 9}, New York (1972).
\bibitem[JS76]{JacquetS} H. Jacquet, J. A. Shalika, {\it An non-vanishing theorem for zeta functions of
             $\GL_n$}, Invent. Math. {\bf 38} (1976), 1--16.
\bibitem[JS81a]{JacquetSmult1} H. Jacquet, J. A. Shalika {\it On Euler Products and the Classification of Automorphic Forms I}
             Am. J. of Math. {\bf 103} (1981): 499–-558.
\bibitem[JS81b]{JacquetSmult} H. Jacquet, J. A. Shalika {\it On Euler Products and the Classification of Automorphic Forms II}
             Am. J. of Math. {\bf 103} (1981): 777–-815.
\bibitem[Kim04]{Kim} H.H. Kim, {\it Funcoriality for the exterior square of $\GL_4$ and the
             sym\-me\-tric fourth of $\GL_2$}, J. of the AMS {\bf 16} (2002), 139--183.
\bibitem[Kn86]{Knapp} W.A. Knapp, {\it Representation Theory of Semisimple Groups}, Princeton University Press,
             Princeton (1986).
\bibitem[KoKo05]{Kohnen} W. Kohnen and H. Kojima, {\it A Maass space in higher genus}, Comp. Math. {\bf 141} (2005),
             313--323.
\bibitem[Kos61]{Kostant} B. Kostant, {\it Lie Algebra Cohomology and the generalized Borel-Weil Theorem},
            Annals of Math. {\bf 74} (1961), 329--387.
\bibitem[KoS99]{KShelstad} R. Kottwitz and D. Shelstad, {\it Foundations of twisted endoscopy},
            Ast\'erisque {\bf 255}  (1999).
\bibitem[KRS92]{Soudry3} St. Kudla, St. Rallis, D. Soudry, {\it On the degree 5 $L$-function for $Sp(2)$},
            Invent. Math. {\bf 107}  (1992), 483--541.
\bibitem[Lan73]{Langlandsreell} R.P. Langlands, {\it On the Classification of
                    Irreducible Representations of Real Algebraic Groups}, Preprint, Institute for Advanced Study (1973).
\bibitem[MW89]  {MoglinWald} C. M\oe glin and J.-L. Waldspurger, {\it Le spectre r\'{e}siduel de GL(n)}
            Ann. Sci. \'{E}cole Norm. Sup. (4) {\bf 22} (1989), 605-674.
\bibitem[Mor11] {Morel}  S. Morel, {\it Cohomologie d'intersection des vari\'{e}t\'{e}s
                modulaires de Siegel, suite},
              Comp. Math. {\bf 147} (2011) 1671--1740.
\bibitem[PS83]{Saito-Kurokawa}I.I. Piatetski-Shapiro, {On the Saito-Kurokawa lifting}, Invent. math. {\bf 71}
            (1983), 309--338.
\bibitem[Ram00]{Ramakrishnan} D. Ramakrishnan, {\it Modularity of the Rankin-Selberg $L$-series, and
             multiplicity one for $\SL(2)$}, Ann. of Math. {\bf 152}(2000), 45--111.
\bibitem[RS07]{RamaSiegel} D. Ramakrishnan, F. Shahidi, {\it Siegel Modular Forms of Genus 2 attached to Elliptic Curves},
             Math. Res. Lett. {\bf 14}, (2007),  315–-332.
\bibitem[Schw86]{SchwSieg1} J. Schwermer, {\it On arithmetic quotients of the Siegel upper half space
            of degree two}, Compos. Math. {\bf 58} (1986), 233--258.
\bibitem[Schw95]{SchwSieg2} J. Schwermer, {\it On Euler products and residual Eisenstein cohomology
            classes for Siegel modular varieties}, Forum Math. {\bf 7} (1995), 1--28.
\bibitem[Sou88]{Soudry Crelle} D. Soudry, {\it The CAP representations of $\GSp(4,\A)$},
              J. reine angew. Math. {\bf 383} (1988), 87--108.
\bibitem[Sou05]{Soudry Paris} D. Soudry, {\it On Langlands functoriality from classical groups to $\GL_n$ },
          In: Tilouine, Jacques (ed.) et al.: Automorphic forms (I) Proceedings of the Semester of the \'{E}mile Borel Center, Paris.
        Ast\'{e}risque {\bf 298} (2005), 335--390.

\bibitem[VZ84]{Vogan} D. Vogan, G. Zuckerman {\it Unitary representations with nonzero cohomology},
            Compos. Math. {\bf 53}  (1984), 51--90.

\bibitem[Wei00]{WWhittaker} R. Weissauer, {\it A remark on the existence of Whittaker models for $L$-packets of
            automorphic representations of $\GSp(4)$}, Manuskripte der For\-scher\-grup\-pe Arithmetik
            {\bf 24} (2000), University of Mannheim
\bibitem[Wei05]{W4dimGal} R. Weissauer, {\it Four dimensional symplectic Galois Representations},
       In: Tilouine, Jacques (ed.) et al.: Automorphic forms (II). The case of the group $\GSp(4)$.
        Astérisque {\bf 302} (2005) 67-150.
\bibitem[Wei09]{WeiBuch} R. Weissauer, {\it Endoscopy for $\GSp(4)$ and the Cohomology of Siegel Modular
         Threefolds}, Springer Lect. Notes  Math. {\bf 1968}, Berlin-Heidelberg (2009).

\bibitem[Wes12]{Uwe} U. Weselmann, {\it A twisted topological trace formula for Hecke operators and liftings from symplectic to general
        linear groups}, Comp. Math. {\bf 148} (2012)  65 -- 120  .
\end{thebibliography}
\end{document}